\documentclass[preprint,3p]{elsarticle}

\usepackage{amssymb}
\usepackage{amsmath}
\usepackage{lineno}

\usepackage{graphicx}
\graphicspath{{./figure}}

\usepackage{booktabs}
\usepackage{threeparttable}
\usepackage{siunitx}
\usepackage{overpic}
\usepackage{amsmath}
\usepackage{amsthm}
\usepackage{amssymb}
\usepackage{mathrsfs}
\usepackage{enumitem}
\usepackage[ruled]{algorithm2e}

\usepackage{mathtools}

\mathtoolsset{showonlyrefs=true}

\usepackage{overpic}
\usepackage[colorlinks=true, pdfstartview=FitV, linkcolor=blue,
citecolor=blue, urlcolor=blue]{hyperref}

\def\p{\boldsymbol{p}}
\def\P{\boldsymbol\Phi}
\def\n{\mathbf{n}}
\def\g{\boldsymbol{g}}
\def\ba{\boldsymbol{a}}
\def\a{\alpha}
\def\b{\beta}
\def\bt{\beta}
\def\c{\theta}

\def\e{\mathbf{e}}
\def\s{\sigma}
\def\eps{{\epsilon}}

\def\toe{{\tilde{\Omega}^\epsilon}}
\def\u{\boldsymbol{u}}

\def\R{{\mathbb R}}
\def\U{\mathcal{ U}}
\def\V{\mathcal{V}}
\def\T{\mathbb T}

\def\t#1{{\tilde{#1}}}

\def\tO{\t{\Omega}}
\def\tw{{\t{\omega}}}
\def\k{\bs{k}}
\def\I{\mathbf{I}}
\def\ka{\kappa}
\def\nb{\nabla}

\def\itd{{\,\mathrm{d}}}

\def\wc{\rightharpoonup}

\def\aeon{\text{a.e. }}
\def\bs#1{\boldsymbol{#1}}

\def\pr#1{\left(#1\right)}

\def\r{\bs{r}}

\def\div{\text{\upshape div }}

\def\br#1{\left\{#1\right\}}
\def\pr#1{\left(#1\right)}
\def\wmk{\text{ W/(mK)}}
\def\nswmk{\text{W/(mK)}}
\def\bx{\bs{x}}
\def\by{\bs{y}}
\def\inspace{\text{\, in\, }}
\def\pd{\partial}

\def\dt#1{\dot{#1}}
\def\tr{\text{tr }}

\def\BB#1{\langle{#1}\rangle}

\def\pd{\partial}
\DeclareMathOperator*{\argmin}{arg\,min}

\newtheorem{remark}{Remark}
\newtheorem{problem}{Problem}
\newtheorem{lemma}{Lemma}
\newtheorem{theorem}{Theorem}
\newtheorem{proposition}{Proposition}

\journal{Computer Methods in Applied Mechanics and Engineering}

\begin{document}

\begin{frontmatter}

\title{Asymptotic analysis and design of shell-based thermal lattice metamaterials}

\author{Di Zhang, Ligang Liu} 

\affiliation{organization={Unversity of Science and Technology of China},
            city={Hefei},
            country={China}
            }

\begin{abstract}
We present a rigorous asymptotic analysis framework for investigating the thermal conductivity of shell lattice metamaterials, extending prior work from mechanical stiffness to heat transfer. Central to our analysis is a new metric, the asymptotic directional conductivity (ADC), which captures the leading-order influence of the middle surface geometry on the effective thermal conductivity in the vanishing-thickness limit. 
A convergence theorem is established for evaluating ADC, along with a sharp upper bound and the necessary and sufficient condition for achieving this bound.
These results provide the first theoretical justification for the optimal thermal conductivity of triply periodic minimal surfaces. Furthermore, we show that ADC yields a third-order  approximation to the effective conductivity of shell lattices at low volume fractions. To support practical design applications, we develop a discrete algorithm for computing and optimizing ADC over arbitrary periodic surfaces.
 Numerical results confirm the theoretical predictions and demonstrate the robustness and effectiveness of the proposed optimization algorithm.
\end{abstract}


\end{frontmatter}

\section{Introduction}

Shell lattice metamaterials are a class of periodic microstructures constructed by thickening smooth periodic surfaces. Among various surface-based designs, those derived from \emph{triply periodic minimal surfaces} (TPMS) have attracted significant attention due to their exceptional physical performance~\cite{tpms-survey}. Their advantageous properties—such as high stiffness and efficient heat transfer—have inspired applications across disciplines including mechanical engineering, thermal management, and materials science~\cite{tpms-survey}. However, the theoretical explanation behind these superior behaviors remains limited, until a recent study~\cite{ads} introduced the concept of \emph{asymptotic directional stiffness} (ADS) to analyze the mechanical properties of TPMS shell lattices. The study reveals that the asymptotic bulk modulus (i.e., ADS under hydrostatic loading), achieves its theoretical upper bound if and only if the surface is a TPMS. While this result explains the optimal mechanical behavior, whether similar optimality holds for other physical phenomena—such as thermal conduction—has remained an open question.

Thermal conductivity is a central property in the design of functional metamaterials. Artificial thermal metamaterials can manipulate heat flow beyond the capabilities of natural materials~\cite{WANG2020101637}, enabling applications such as thermal cloaking, heat concentrators, and compact heat exchangers~\cite{PhysRevLett.108.214303,nat-rev-mat}. Notably, numerical and experimental studies have observed that TPMS shell lattices exhibit high thermal conductivity~\cite{GAO2023106976,MIRABOLGHASEMI201961}, leading to their use in advanced heat transfer devices~\cite{Gado2024,OH2023103778}. Yet, the underlying mechanism responsible for this enhancement has not been theoretically established.

In this study, we aim to establish the theoretical foundation for the high conductivity observed in TPMS shell lattices. We introduce a new metric, termed \emph{asymptotic directional conductivity} (ADC), to investigate how the geometry of the middle surface influences the conductivity of shell lattices as their thickness approaches zero. 
We establish a convergence theorem  for evaluating ADC, along with a necessary and sufficient condition for achieving its upper bound, which reveals the optimality of TPMS in thermal conductivity.
In addition, we show that ADC provides third-order accuracy in approximating the effective conductivity at low volume fractions.

For the practical design of shell lattices, we further propose a discretization algorithm to evaluate and optimize ADC based on the theory of discrete differential geometry~\cite{pmp-book}.
This discretization scheme, together with our optimality condition, enables the numerical generation of TPMS by maximizing the average asymptotic conductivity (Section~\ref{sec:opt-adc-tpms}).
As we will demonstrate in Section~\ref{sec:gen-tpms-aac}, this approach is more robust than mean curvature flow (MCF)~\cite{surf-evolver,Dziuk1990} and more flexible than conjugate surface construction methods~\cite{Makatura-proce,Karcher1996} in practice.
Additionally, our derivation also extends to other physical properties governed by similar mathematical formulations, such as electrical conductivity and mass diffusivity.

The remainder of this paper is organized as follows. Section~\ref{sec:rel-work} reviews relevant background on metamaterial design, the generation and application of TPMS, and asymptotic analysis. Section~\ref{sec:shell-latt} introduces the mathematical model of shell lattice metamaterials and its effective thermal conductivity, with auxiliary details provided in~\ref{sec:eff-cond} and~\ref{sec:shell-repr}. In Section~\ref{sec:asym-analysis}, we present our key results through asymptotic analysis, where the rigorous justification is provided in~\ref{app:asym-ana}.
 Section~\ref{sec:discrete} illustrates the numerical discretization for evaluating and optimizing ADC. 
 Section~\ref{sec:valid} provides numerical validations of our theory, including convergence analysis and comparison with experimental data.
  Section~\ref{sec:app} demonstrates applications in conductivity optimization and TPMS generation. 
  Finally, Section~\ref{sec:concl} concludes the paper and discusses future research directions.

\section{Related work}
\label{sec:rel-work}
\subsection{Metamaterial design}
Metamaterials are engineered materials whose effective macroscopic behavior is governed by the architecture of their internal microstructure. Their unit cell, often referred to as a representative volume element (RVE), serves as their building block that governs the macroscopic behavior~\cite{Kadic2019-mh}.
By tailoring the geometry of RVEs, metamaterials can exhibit effective properties that exceed those of natural materials~\cite{SUN2025145}. As a result, they have found widespread applications across various fields, including mechanics, acoustics, optics, and thermal transport~\cite{ASKARI2020101562}.
A fundamental tool for metamaterial design is homogenization theory~\cite{Allaire2002,Milton_2002}, which establishes the computational relationship between the structural configuration of RVEs and their effective macroscopic properties.
Optimizing these effective properties naturally leads to a class of topology optimization problems, which can be addressed using classical methods such as solid isotropic material with penalization (SIMP)~\cite{simp-mat,Zhou2008-gm,RADMAN2014266,porus-heat-limit}, evolutionary structural optimization (ESO)~\cite{eso-top,beso-review}, and level set methods~\cite{LIU2020113154}, among others.
In the graphics community, parametric approaches are commonly used to design microstructures with prescribed elastic properties~\cite{meta-contr-el2015,elastic-texture,worst-case-stress,flex-material}, enabling downstream applications such as deformation control~\cite{meta-contr-el2015,elastic-texture} and multiscale topology optimization~\cite{Wu2021-oz,Zhubo2017}.

Existing studies of shell lattice metamaterials have primarily focused on their mechanical behavior.
It has been shown that such lattices can achieve extremely low volume fractions while maintaining high strength and stiffness~\cite{shellular-mtl,Han2017-d-shell}.
Further investigation into their mechanical performance, particularly for TPMS-based shell structures, has been carried out through comprehensive simulations~\cite{BONATTI2019301,WANG2020108340}.
Inverse design methodologies of shell lattices typically utilize Kirchhoff-Love shell model for finite element simulation and sensitivity analysis~\cite{MA2022110426,opencell2023}, where cubic symmetry is assumed to simplify the NURBS representation of the middle surface.
Although numerical studies suggest that TPMS shell lattices also exhibit superior heat conduction~\cite{GAO2023106976,MIRABOLGHASEMI201961}, theoretical justification and shape optimization algorithms for this property are, to the best of our knowledge, still lacking.
The most related approaches~\cite{hjb-tpms-heat} optimizes the heat dissipation of macroscopic shell structures instead of designing metamaterials.

In addition to numerical optimization,
several analytical constructions~\cite{Allaire2002,GIBIANSKY2000461} have been  demonstrated to achieve the Hashin-Shtrikman (HS) upper bound~\cite{hashin1962variational} for bulk modulus and conductivity.
These constructions were significantly extended by~\cite{ads}, which demonstrated that any TPMS shell lattice achieves HS upper bound for bulk modulus when thickness is sufficiently small.
In the present study, we further show that these shell lattices also attains the HS upper bound for thermal conductivity.

\subsection{Construction and application of TPMS}
Traditional methods for constructing TPMS rely either on analytical expressions of fundamental minimal surface patches combined with symmetric duplication~\cite{FischerKoch1987,Karcher1996,schoen1970infinite}, or on numerical techniques such as MCF coupled with conjugate surface construction~\cite{SE1992,Ulrich1993,Makatura-proce}.
Recently, \cite{Stephanie2021} uses geometric measure theory (GMT) to convert the Plateau problem into a convex optimization problem, which yields globally optimal solutions while avoiding issues related to mesh quality and topological changes.
This idea is further explored in DeepCurrents~\cite{Palmer9879635}, which employs a neural network to represent current, enabling more efficient and robust computation.
As we will show in Section~\ref{sec:gen-tpms-aac}, the established optimal condition of ADC allows us to generate TPMS from a new perspective, i.e., by optimizing the average asymptotic conductivity.
Compared to  MCF-based method, which are often prone to degeneration~\cite{Makatura-proce,Karcher1996}, our approach is more robust and avoids mesh collapse in all experiments.
Moreover, it deforms arbitrary periodic surfaces into TPMS without requiring prior information such as user-specified edge boundaries, as needed in conjugate surface construction methods~\cite{Makatura-proce,SE1992} and  GMT based approaches~\cite{Stephanie2021,Palmer9879635}.

TPMS are known to exhibit superior physical properties across various domains, including mechanics, thermodynamics, and optics~\cite{FENG2021110050}. Over the past decades, they have been widely applied in additive manufacturing and metamaterial design~\cite{tpms-survey,TORQUATO2024120142,tpms-shell-lvlin,shell-lattice}, with most studies focusing on their mechanical performance.
Some investigations have explored the thermal transport properties of TPMS shell lattices through numerical simulations~\cite{MIRABOLGHASEMI201961,QURESHI2021121001,Forced-convec,Flow-and-thermal} and physical experiments~\cite{Zhou2300359}.
In the graphics community, TPMS have been adopted for designing internal supporting structures in 3D models to enhance stiffness~\cite{hjb-tpms-mech,8703138,10105512} and improve heat dissipation~\cite{hjb-tpms-heat}.
These works either focus on designing macroscopic structures rather than developing metamaterials, or are limited to specific types of TPMS, resulting in a lack of general conclusions. 

\subsection{Asymptotic analysis}
Our analysis of thermal shell lattices builds upon the framework developed in~\cite{ads}, which presents an asymptotic analysis of the stiffness of shell lattice metamaterials grounded in Ciarlet’s shell theory~\cite{CiarMem1996,CiarletFlex1996,CiarletGenmem1996}.
That work introduces a metric, termed \emph{asymptotic directional stiffness} (ADS), to quantify the influence of the middle surface geometry on the effective stiffness of shell lattices.
It establishes both an upper bound for ADS and the necessary and sufficient condition for achieving this bound, which provides the first rigorous justification of the optimal bulk modulus of TPMS.
In contrast to the linear elasticity problem, where the physical field is a  displacement vector field, the thermal conduction problem involves a scalar temperature field. This distinction allows us to derive stronger results  than in the mechanical case.
As we will show in Section~\ref{sec:upp-bound-adc} and~\ref{sec:opt-adc-tpms}, the optimality of conductivity in TPMS shell lattices holds uniformly in all directions, unlike the elastic case where optimal stiffness is restricted to hydrostatic loading conditions.
Furthermore, the ADC exhibits third-order accuracy in approximating the effective conductivity of shell lattices at low volume fractions.
From a computational perspective, the numerical evaluation and optimization of the proposed ADC metric are naturally formulated within the framework of discrete differential geometry~\cite{pmp-book}, which greatly facilitates implementation.

\section{Shell lattice metamaterial}
\label{sec:shell-latt}
\begin{figure}[t]
	\centering
	\begin{overpic}[width=0.6\textwidth,keepaspectratio]{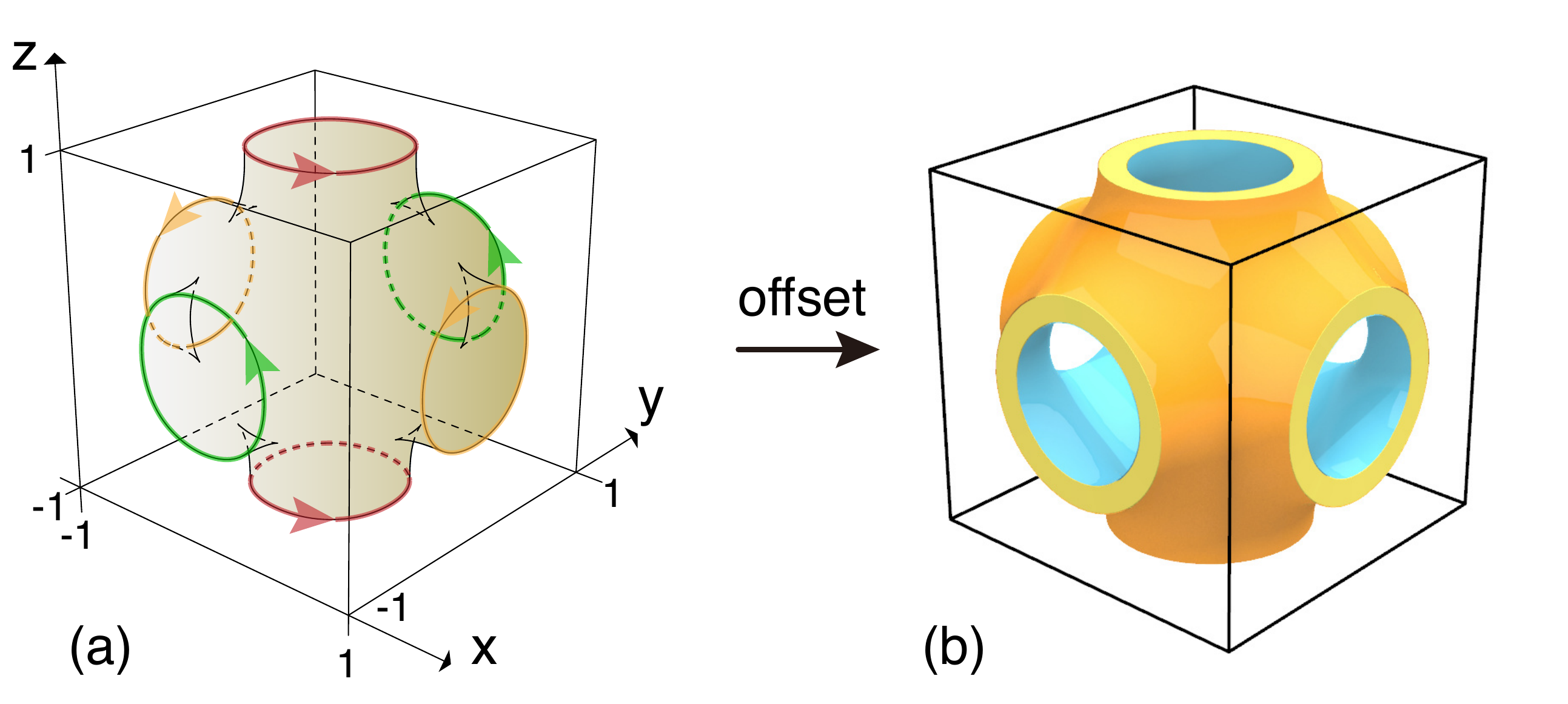}
		\put(8,38) {$\mathbb T^3$}
		\put(17,24) {$\tw$}
		\put(74,23) {$\toe$}
	\end{overpic}
	\caption{
		(a) a periodic surface $\tw$ embedded in $\T^3$, which is closed due to the the periodic boundaries (colored arrows).
		(b) a derived shell lattice $\toe$ by offsetting  $\tw$ normally on both sides.
	}
	\label{fig:thick-shell}
\end{figure}

The shell lattice metamaterial is constructed by thickening a smooth periodic surface within a RVE~\cite{ads}. For simplicity, we take $Y := [-1,1]^3$ as the RVE. Due to the periodic boundary conditions, $Y$ is identified with the 3-dimensional torus $\T^3 := \R^3 / (2\mathbb{Z})^3$.
The middle surface is a closed smooth  surface embedded in $\T^3$, denoted as $\tw$ (Figure~\ref{fig:thick-shell}(a)).
The shell lattice, denoted as $\toe$, is obtained by offsetting $\tw$ normally on both sides within $\T^3$ by a distance $\eps$ (Figure~\ref{fig:thick-shell}(b)).
We refer the reader to~\cite{ads} or \ref{sec:shell-repr} for more details, which introduces the necessary definitions and lemmas before proving the convergence theorem in Section~\ref{sec:conv-theorem}.

The effective conductivity of the shell lattice $\toe$ in the direction $\p \in \mathbb{R}^3$ is given by
\begin{equation}
	\label{eq:ka-eps-matrix-form}
	\ka_\eps(\tw;\p) = \k_\eps\cdot\p\cdot\p,
\end{equation}
where $\k_\eps$ denotes the effective conductivity matrix determined by homogenization theory~\cite{Allaire2002} (see~\ref{sec:eff-cond}).
The components of $\k_\eps$ are defined as
\begin{equation}
	\label{eq:k-eps-wise-expr}
	k_\eps^{ij} =\frac{1}{|Y|}\int_{\toe}\k\cdot\pr{\nabla u^i+\e^i}\cdot\pr{\nabla u^j+\e^j},
\end{equation}
where $\{\e^i\}_{i=1}^3$ are the canonical basis vectors in $\mathbb{R}^3$, and each $u^i$ is the solution to the following \emph{cell problem}:
\begin{equation}
	\label{eq:min-energy-main}
	\argmin_{u\in V_\#(\toe)} \frac{1}{|Y|}\int_{\toe}\k\cdot(\nb u +\p)\cdot(\nb u+\p)
\end{equation}
evaluated at $\p = \e^i$.
The function space is defined as $V_\#(\toe) := \{ u \in H^1(\toe) : \int_{\toe} u = 0 \}$, enforcing periodicity and zero average. The base material is assumed to be isotropic, with constant thermal conductivity matrix $\k = \ka \I$.

\begin{remark}
	\label{rem:p}
	Since $\ka_\eps(\tw; \p)$ is a quadratic form in $\p$, we consider only unit vectors (i.e., $\|\p\| = 1$) throughout this paper.
	We further denote the tangential and normal components of $\p$ on $\tw$ as $\p_{\tw}$ and $p_3$, respectively.
\end{remark}

\section{asymptotic analysis}
\label{sec:asym-analysis}
To investigate how the geometry of the middle surface influences the heat transfer properties of shell lattices,
we analyze the asymptotic behavior of the effective conductivity as the thickness tends to zero.
Since the effective conductivity vanishes in this limit, we introduce the \emph{asymptotic directional conductivity} (ADC) to extract the leading-order behavior, defined as
\begin{equation}
	\label{eq:dir-asymp-cond}
	{\ka}_A(\tw;\p):= \lim_{\eps\to0} \frac{{\ka}_\eps(\tw;\p)}{\rho_\eps(\tw)},
\end{equation}
where 
\begin{equation}
	\label{eq:rho-define}   
	\rho_\eps(\tw):=|\toe|/|Y|
\end{equation}
is the volume fraction of the solid region within the RVE. 
\subsection{Convergence theorem}
\label{sec:conv-theorem}
Our first key result is the  theorem below:
\begin{theorem}
	\label{thm:asym-cond}
	The limit in~\eqref{eq:dir-asymp-cond} admits the following explicit expression  
	\begin{equation}
		\begin{aligned}
			\label{eq:ka-dir-expr}
				{\ka}_A(\tw;\p)=\frac{1}{|\tw|}\int_{\tw}\ka\pr{\nabla\bar u+\bs p_{\tw}}\cdot\pr{\nabla\bar u +\bs p_{\tw}},
		\end{aligned}
	\end{equation}
	where $\bar{u}$ is the solution to the following equation on $\tw$:
	\begin{equation}
		\label{eq:u-div-eq-hp3}
		\Delta \bar{u} = -\text{\upshape div}\,\p_{\tw}.
	\end{equation}
\end{theorem}
See the detailed proof in \ref{sec:deri-limit-equation}.

\paragraph*{ADC matrix}
Since the solution $\bar{u}$ to~\eqref{eq:u-div-eq-hp3} is linearly dependent  on $\p$,   we can rewrite~\eqref{eq:ka-dir-expr} as the following quadratic form
\begin{equation}
	\label{eq:adc-matrix-form}
	\ka_A(\tw;\p)=\boldsymbol{k}_A\cdot \p\cdot\p.
\end{equation}
where the matrix $\k_A$ is referred to as the \emph{asymptotic conductivity matrix} (or ADC matrix), defined by
\begin{equation}
	\label{eq:bar-k-A-expr}
	k_A^{ij}:=\frac{1}{|\tw|}\int_{\tw}\ka\pr{\nabla\bar u^i+\e^i_{\tw}}\cdot\pr{\nabla\bar u^j+\e^j_{\tw}},
\end{equation}
in which $\bar{u}^i$ is the solution to~\eqref{eq:u-div-eq-hp3} when $\p = \e^i$, and $\e^i_{\tw}$ denotes its  tangential component on $\tw$.

By comparing the formulations in~\eqref{eq:adc-matrix-form},~\eqref{eq:ka-eps-matrix-form}, and~\eqref{eq:dir-asymp-cond}, we obtain the following relation:
\begin{equation}
	\label{eq:kmat-eps-to-kmat-A}
	\k_A=\lim_{\eps\to0}\frac{\k_\eps}{\rho_\eps}.
\end{equation}
In other words, the ADC matrix $\k_A$ represents the leading-order term of  $\k_\eps$.

\subsection{Upper bound of ADC}
\label{sec:upp-bound-adc}
The convergence theorem implies the following upper bound of $\ka_A$:

\begin{theorem}
	\label{thm:adc-upper-bound}
	The ADC of $\tw$ satisfies
	\begin{equation}
		\label{eq:ka-upp-bound-p3}
		\ka_A(\tw;\p)\le \frac{1}{|\tw|}\int_{\tw}\ka\pr{1-p_3^2}
	\end{equation}
	and  equality holds if and only if $\tw$ is minimal or a cylinder parallel to $\p$.
\end{theorem}
\begin{proof}
	From~\eqref{eq:ka-dir-expr}, we  can rewrite the integrand as:
	\begin{equation}
		\begin{aligned}
			\label{eq:ka-dir-hp-form}
			{\ka}_A(\tw;\p)&=\frac{1}{|\tw|}\int_{\tw}\ka\p_{\tw}\cdot\p_{\tw}+\ka(\nabla\bar u+\p_{\tw})\cdot\nabla\bar u+\ka\nabla\bar u\cdot\p_{\tw} \\
		\end{aligned}
	\end{equation}
	The weak form of~\eqref{eq:u-div-eq-hp3} is the following variational equation:
	\begin{equation}
		\int_{\tw}\pr{\nabla\bar u+\p_{\tw}}\cdot\nabla v=0,\quad\forall v\in H^1(\tw),
	\end{equation}
	which leads to the cancellation of the middle term in~\eqref{eq:ka-dir-hp-form} and replaces $\nabla\bar u\cdot\p_{\tw}$ with $-|\nabla\bar u|^2$.
	Consequently, we obtain
	\begin{equation}
		\label{eq:kap-ku2}
		{\ka}_A(\tw;\p)=\frac{1}{|\tw|}\int_{\tw}\ka\pr{1-p_3^2}-\frac{1}{|\tw|}\int_{\tw}\ka|\nabla\bar u|^2.
	\end{equation}
	where we have used the identity $\|\p_{\tw}\|^2 = 1 - p_3^2$ (Remark~\ref{rem:p}).
	Since the second term in~\eqref{eq:kap-ku2} is non-negative, the inequality~\eqref{eq:ka-upp-bound-p3} follows immediately.
	
	To determine the condition for equality, we analyze the second term: 
	\begin{equation}
		\label{eq:prove-equali-cond}
		\begin{aligned}
			\int_{\tw}\ka |\nabla\bar u|^2 =-\int_{\tw} \ka \bar u \Delta \bar{u}
			=\int_{\tw} \ka \bar u \div \p_{\tw} = \int_{\tw} 2\ka \bar u Hp_3,
		\end{aligned}
	\end{equation}
	where $H$ denotes the mean curvature of $\tw$, and the identity $\div \p_{\tw} = 2 H p_3$ is provided by Lemma~\ref{lem:div-tg-hp3}.
	The final integral in~\eqref{eq:prove-equali-cond} vanishes if $H p_3 = 0$ everywhere on $\tw$, which holds if either $H = 0$ (i.e., $\tw$ is a minimal surface) or $p_3 = 0$ (i.e., $\tw$ is a cylinder aligned with $\p$).
	The necessity of this condition is proved in Lemma~\ref{lem:nec-cond-upp}. 
\end{proof}

\subsection{The optimal ADC of TPMS}	
\label{sec:opt-adc-tpms}

From~\eqref{eq:adc-matrix-form}, we know $\ka_A(\tw;\p)$ has three principle directions, denoted as $\p_1$, $\p_2$, and $\p_3$.
The corresponding values of ADC in these directions are referred to as the \emph{principal asymptotic conductivities}, denoted $\ka_A(\tw; \p_i)$.
We define their average as the \emph{average asymptotic conductivity} (AAC), written as:
\begin{equation}
	{\ka}_A^{Avg}(\tw):=\frac{1}{3}\sum_i \ka_A(\tw;\p_i).
\end{equation}
We now state the following theorem.
\begin{theorem}
	\label{thm:tpms-opt-adc}
	The AAC of $\tw$ satisfies 
	\begin{equation}
		\label{eq:avg-const-upp-bound}
		\ka_A^{Avg}(\tw)\le\frac{2}{3}\ka =: \ka_A^*
	\end{equation}
	and equality holds if and only if $\tw$ is TPMS.
\end{theorem}
\begin{proof}
	Applying~\eqref{eq:ka-upp-bound-p3} to each $\p_i$, we obtain
	\begin{equation}
		\label{eq:avg-upp-sum}
		\begin{aligned}
			{\ka}_A^{Avg}(\tw)&\le \frac{1}{3}\sum_i\frac{1}{|\tw|}\int_{\tw}\ka\pr{1-\pr{\p_i\cdot\n}^2}\\
			&=\ka\pr{1-\frac{1}{|\tw|}\int_{\tw}\frac{1}{3}\sum_i\pr{\p_i\cdot \n}^2}
		\end{aligned}
	\end{equation}
	Since $\{\p_i\}$ forms an orthonormal basis, we have: $$\frac{1}{3}\sum_i(\p_i\cdot\n)^2=\frac{1}{3}\|\n\|^2=\frac{1}{3}.$$
	Substituting this into~\eqref{eq:avg-upp-sum} yields the upper bound~\eqref{eq:avg-const-upp-bound}.
	
	Now suppose equality holds in~\eqref{eq:avg-const-upp-bound}. Then equality must also hold in~\eqref{eq:ka-upp-bound-p3} for each $\p_i$.
	By Theorem~\ref{thm:adc-upper-bound}, this implies that $\tw$ must be either minimal ($H = 0$), or a cylinder parallel to each $\p_i$.
	The latter case is impossible, since the unit normal vector $\n$ cannot be simultaneously orthogonal to all three orthonormal directions $\{\p_i\}$.
	Therefore, $\tw$ must be a minimal surface.
\end{proof}
\begin{remark}
	Note that the inequality~\eqref{eq:avg-const-upp-bound} remains valid even when the surface is not closed or the shell thickness is non-uniform (\ref{sec:kavg-upp-general-shell}).
\end{remark}

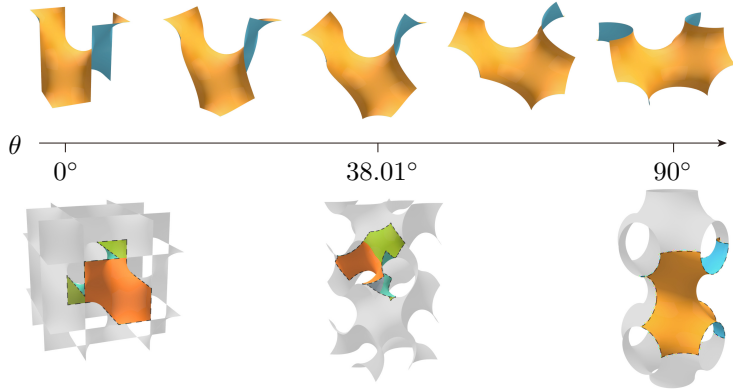
\begin{figure}[t]
	\centering
	\begin{overpic}[width=0.6\textwidth,keepaspectratio]{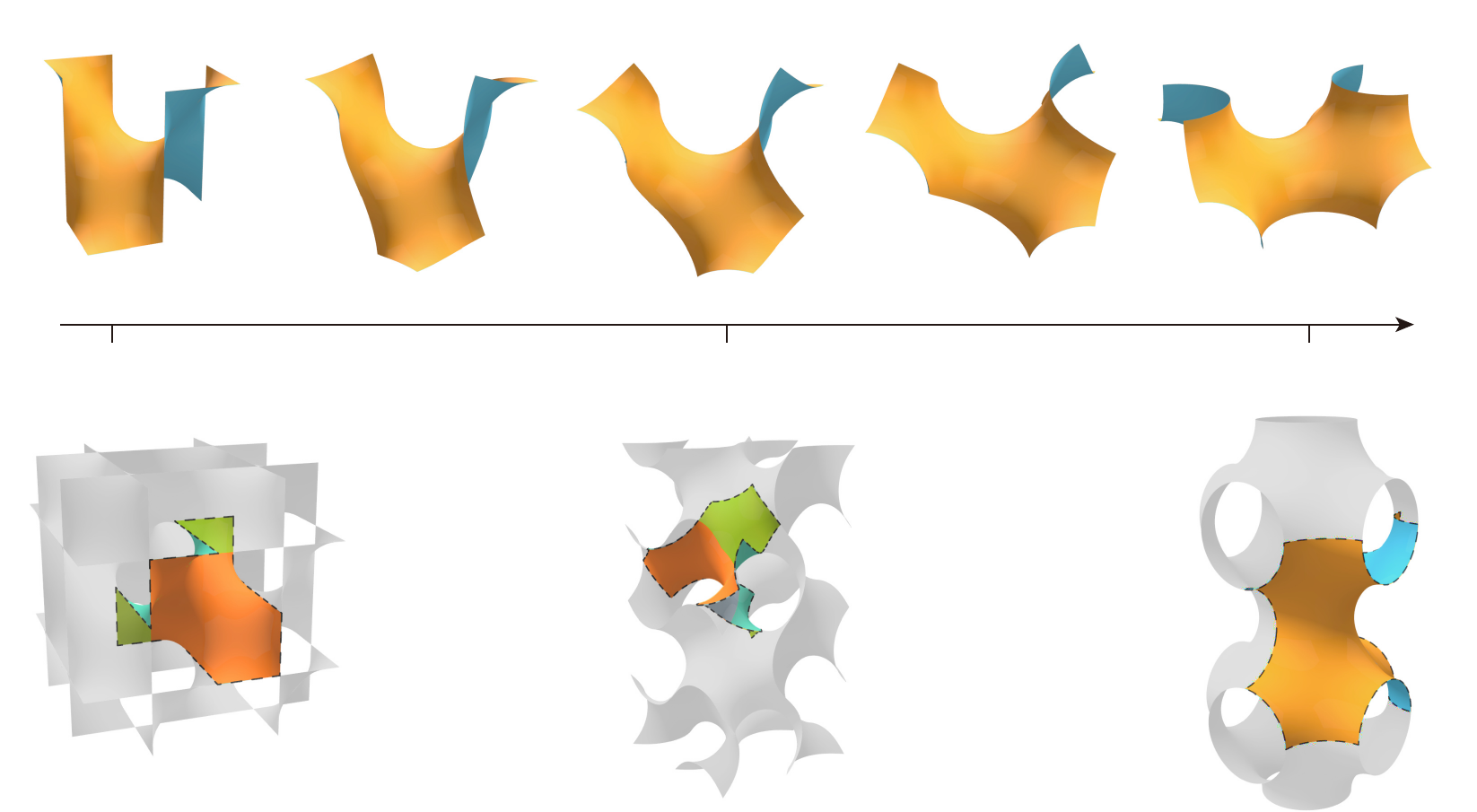}
		\put(0,31.5) {$\theta$}
		\put(6,28) {$\ang{0}$}
		\put(45,28) {$\ang{38.01}$}
		\put(86,28) {$\ang{90}$}
	\end{overpic}
	\vspace{-1mm}
	\caption{
		Associate family of Schwarz D surface.  
		$\theta$ denotes the Bonnet rotation angle.
		The surface is a patch of Gyroid and Schwarz P when  $\theta=\ang{38.01}$ and $\ang{90}$, respectively.
	}
	\label{fig:asso-family}
\end{figure}

\subsubsection{Isotropic ADC}
A surface $\tw$ is said to exhibit isotropic ADC if $\ka_A(\tw; \p)$ remains constant for all directions $\p$. In this case, the ADC equals its average value $\ka_A^{\mathrm{Avg}}$, thus the upper bound in Theorem~\ref{thm:tpms-opt-adc} applies uniformly in all directions:
\begin{theorem}
	If $\tw$ exhibits isotropic ADC, then
	\begin{equation}
		\label{eq:iso-cub-cons-bound}
		\ka_A(\tw)\le \frac{2}{3}\ka,
	\end{equation}
	and  equality holds if and only if $\tw$ is a TPMS (with isotropic ADC).
\end{theorem}

A typical class of surfaces with isotropic ADC are those with cubic symmetry.
This is claimed by Neumann’s principle~\cite{nye1985physical}, which states that cubic lattices always exhibit isotropic effective conductivity.
We can also show that  the RHS in~\eqref{eq:ka-upp-bound-p3}  simplifies to the constant $\frac{2}{3}\ka$ in this special case (see Lemma~\ref{lm:cubic-srf-int}).

\subsubsection{Associate family}
In addition to cubic-symmetric TPMS, Theorem~\ref{thm:adc-upper-bound} suggests a broader class of surfaces that achieve optimal isotropic ADC.
Specifically, it implies that the ADC of TPMS in a given direction $\p$ is determined solely by the distribution of its normal vectors.
It is well known that each minimal surface belongs to an associate family, comprising minimal surfaces that share  identical Weierstrass data and differ only by a Bonnet rotation~\cite{lawson-minsurf,conj-surf}.
A fundamental  property of such families is that all member surfaces  have the same normal vector distribution,  and thus exhibit the same ADC.
One typical example is the Schwarz D family, which includes Schwarz P and Gyroid surfaces (Figure~\ref{fig:asso-family}).
Despite the Gyroid and Schwarz D not being cubic-symmetric, they retain isotropic ADC due to the inclusion in the  family of Schwarz P, which is cubic-symmetric.

\subsection{Relation to Hashin-Shtrikman bound}
We stress that the upper bound in~\eqref{eq:iso-cub-cons-bound} aligns with the classical Hashin–Shtrikman (HS) upper bound~\cite{hashin1962variational,PABST2007479}, which provides a theoretical limit on the effective conductivity of isotropic composites. For a given volume fraction $\rho(\toe)$, the HS bound is given by:
\begin{equation}
	\label{eq:hs-bound-iso}
	\ka_\eps(\tw)  \le \ka^{HS}\pr{\rho(\toe)}:= \frac{2\rho(\toe)}{3-\rho(\toe)} \ka.
\end{equation}
Applying this bound to the definition of ADC yields:
\begin{equation}
	\kappa_A(\tilde\omega)\le \limsup_{\epsilon\to0}\frac{\kappa_\epsilon(\tilde\omega)}{\rho_\epsilon(\tilde\omega)}\le\lim_{\rho\to0}\frac{\ka^{HS}\pr{\rho}}{\rho}=\frac{2}{3}\kappa,
\end{equation}
which is exactly the asymptotic bound in~\eqref{eq:iso-cub-cons-bound}.
Therefore, our theorems suggest that \emph{any cubic TPMS shell lattice achieves the HS upper bound when the volume fraction is sufficiently small}.

\subsection{Third-order accuracy of ADC}
By definition~\eqref{eq:dir-asymp-cond}, the asymptotic directional conductivity $\ka_A(\tw; \p)$ captures the leading-order (first-order) behavior of the effective directional conductivity as the shell thickness $\eps$ tends to zero.
Interestingly, we show that this approximation is significantly more accurate, achieving third-order precision in $\eps$.
\begin{theorem}
	\label{thm:adc-3-order}
	The ADC of $\tw$ satisfies
	\begin{equation}
		\label{eq:adc-3-order}
		\ka_\eps(\tw;\p) -\frac{2\eps|\tw|}{|Y|}\ka_A(\tw;\p)\sim O(\eps^3) \quad as\,\,\eps\to0.
	\end{equation}
\end{theorem}
Here, the term ${2\eps |\tw|}/{|Y|}$ corresponds to the first-order approximation of the volume fraction $\rho_\eps(\tw)$.
The detailed proof of this result is provided in the \ref{app:third-oder-acc}, and numerical validation is presented in Section~\ref{sec:3order-acc}.

\begin{remark}
	The factor ${2\eps |\tw|}/{|Y|}$ in~\eqref{eq:adc-3-order} can be replaced by the full volume fraction $\rho_\eps(\tw)$, as the difference between them is also $O(\eps^3)$. See Lemma~\ref{lem:vol-remainder} for details.
\end{remark}

\section{Discretization}
\label{sec:discrete}
We discretize the continuous surface $\tw$ using a triangular mesh consisting of a set of vertices $\mathcal{V}$, edges $\mathcal{E}$, and faces $\mathcal{F}$. 
In Section~\ref{sec:dis-adc-form}, we first describe the numerical evaluation of the ADC matrix $\k_A$. 
Then, in Section~\ref{sec:shap-opt}, we derive the $L^2$ gradient of the ADC under normal deformation, which is used in a gradient descent flow for our shape optimization algorithm.

\subsection{Evaluation of ADC matrix}
\label{sec:dis-adc-form}
Following  the derivation of~\eqref{eq:kap-ku2},  we rewrite~\eqref{eq:bar-k-A-expr} in the form
\begin{equation}
	\begin{aligned}
		k_A^{ij}&= \frac{1}{|\tw|}\int_{\tw}\left(\e^i_{\tw}\cdot\e^j_{\tw} -\nabla\bar u^i\cdot\nabla\bar u^j\right).  \\
	\end{aligned}
\end{equation}
Using the identity
\begin{equation}
	\e^i_{\tw}\cdot\e^j_{\tw} = (\e^i\cdot\e^j) - (\e^i\cdot\n)(\e^j\cdot\n)=\delta_{ij}-n_in_j,
\end{equation}
we obtain the following matrix expression for the ADC:
\begin{equation}
	\begin{aligned}
		\k_A=\ka\pr{\I-\frac{1}{|\tw|}\int_{\tw}\n\n^\top-\mathbf{R}},
	\end{aligned}
\end{equation}
where $\mathbf I$ is the identity matrix, and $\mathbf{R}$ is given by
\begin{equation}
	\label{eq:asym-mat-L-def}
	[\mathbf R]^{ij} = -\frac{1}{|\tw|}\int_{\tw} \bar u^i\Delta \bar u^j,
\end{equation}
using the identity $\int_{\tw}\nabla\bar u^i\cdot\nabla\bar u^j=-\int_{\tw}\bar u^i\Delta \bar u^j$  through integration by parts.

The computation of $\k_A$ thus involves the following three steps:
\begin{itemize}
	\item Evaluate the integral $\int_{\tw}\n\n^\top$.
	We approximate this integral by 
	\begin{equation}
		\int_{\tw}\n\n^\top\approx \sum_{f\in \mathcal F} \n_f\n_f^\top A_f,
	\end{equation}
	where $\n_f$ denotes the normal vector of triangular face $f$.
	The notation $A_f$ denotes the area of $f$.
	
	\item Solve  for solution $\bar u^i$ from the equation $\Delta u=-\div\p_{\tw}$ on $\tw$ with $\p=\e_i$, i=1,2,3. We first represent $\bs p_{\tw}$ as its line integral on edges $\mathcal E$.
	Suppose $e_{ij}\in\mathcal E$ denotes the edge linking $\bs v_i\in\mathcal V$ to $\bs v_j\in\mathcal V$,  we assign a variable $p_{ij}$ to it, given by
	\begin{equation}
		p_{ij} := \int_{e_{ij}} \p_{\tw}\cdot d\r=\int_{e_{ij}}\p \cdot d\r \approx\p\cdot \pr{\bs v_j - \bs v_i}.
	\end{equation}
	The second equality holds because each edge is regarded as the geodesic path on $\tw$ and thus $d\r$ is tangent on $\tw$. 
	Then the divergence at vertex $v_i$ is approximated by its integral as
	\begin{equation}
		\rho_{v_i} :=-\int_{A_{v_i}} \div\p_{\tw}\itd\tw \approx -\sum_{v_j\in N(v_i)} w_{ij} p_{ij},
	\end{equation}
	where $N(v_i)$ signifies the set of adjacent vertices of $v_i$; the notation $w_{ij}$ denotes the cotangent weight and $A_{v_i}$ signifies the local averaging region~\cite{pmp-book} at vertex $\bs v_i$.
	Next, we stack all $\rho_{v_i}$ into a vector $\bs \rho=(\rho_{v_1},\rho_{v_2},\cdots,\rho_{v_n})^\top$, and represent the unknown as a vector $\bs u=\pr{u_{v_1}, u_{v_2},\cdots,u_{v_n}}^\top$, which stores the value of the solution $u$ at each vertex.  
	The Laplace operator $\Delta$ is discretized as the cotangent Laplacian matrix $\mathbf L$~\cite{pmp-book}.
	Then~\eqref{eq:u-div-eq-hp3} is discretized into the following linear equation
	\begin{equation}
		\label{eq:discr-poisson-eq}
		\mathbf L \bs u = \bs{ \rho}.
	\end{equation}
	After solving $\bs u$ from this equation, we obtain the solution's value at each vertices.
	We repeat this process with $\p=\e_i$, $i=1,2,3$ and denote the corresponding $\bs\rho$ and solution $\bs u$ as $\bs\rho^i$ and $\u^i$, respectively.
	\item Evaluate integral $\int_{\tw} \bar u^i\Delta \bar u^j$. We approximate it as 
	\begin{equation}
		\int_{\tw} \bar u^i\Delta \bar u^j \approx \u^i \cdot \bs\rho^j.
	\end{equation}
	Here, we have utilized the relation $\Delta \bar u^j = -\div\p_{\tw}$ where $\p=\e_j$.
\end{itemize}
The accuracy and convergence rate of this approximation are investigated in Section~\ref{sec:dis-conv-rate}. 

\subsection{Shape optimization}
\label{sec:shap-opt}
We adopt the shape optimization framework introduced in~\cite{ads} to optimize the ADC objective.
As tangential movement does not change the shape of $\tw$, we consider the evolution of $\tw$ under a normal velocity field $v_n: \tw \to \R$.

According to~\eqref{eq:ka-dir-expr},  ADC is computed as the ratio of an integral (denoted as $I_a^c$) to the surface area $A:=|\tw|$.
The time derivative of ADC is then expressed as the following form
\begin{equation}
	\label{eq:common-ad-rate}
	\dot{\kappa}_A=\frac{\dot{I}_a^c}{A}-\frac{I_a^c \dot{A}}{A^2}.
\end{equation}
To derive a gradient descent flow, we first identify the $L^2$-gradient of $\ka_A$ by requiring:
\begin{equation}
	\dot{\ka}_A[v_n]=\langle\text{grad }\ka_A, v_n \rangle_{L^2(\tw)}.
\end{equation}
Once the gradient $\text{grad }\ka_A$ is obtained, the optimization proceeds via an $L^2$-gradient flow:
\begin{equation}
	\label{eq:l2-grad-flow}
	\frac{d\tw}{dt}=-\text{grad }\ka_A.
\end{equation}

The change of area is easily computed by
\begin{equation}
	\label{eq:area-rate}
	\dot{A}:=-\int_{\tw}2v_n H.
\end{equation}
The remaining task is to compute the time derivative $\dot{I}_a^c$, which we address in the next section.

\subsubsection{Sensitivity analysis}
\label{sec:sens-ana}
To compute the $L^2$ gradient,  we first state the following result:
\begin{proposition}
	\label{prop:adc-sens}
	The time derivative $\dot{I}_a^c$ is given by
	\begin{equation}
		\label{eq:Ic-grad}
		\dot{I}_a^c=\int_{\tw}2v_n(\bs b-H\bs I)(\nabla \bar u+\p_{\tw})\cdot(\nabla \bar u+\p_{\tw}),
	\end{equation}
	where $\bar u$ is the solution of~\eqref{eq:u-div-eq-hp3}, $\bs{b}$ is the second fundamental form, $H$ is the mean curvature, and $\bs{I}$ is the identity matrix.
\end{proposition}
We provide the derivation details in \ref{app:dtime-flow}.

The normal velocity of $v_n$ is approximated by piecewise linear functions on $\mathcal M$.
Then time derivative $\dt{I}_a^c$ is discretized as
\begin{equation}
	\label{eq:Ic-grad-h}
	\small
	\begin{aligned}
		\dt{I}^c_a 
		&\approx\sum_{f\in\mathcal{F}}\sum_{i\in V(f)}2v_n^i\int_{f}\phi_i\pr{\nabla\bar u^h+P_f\bs p}^\top\pr{\bs b_f -\frac{\tr{\bs b_f}}{2}\mathbf I}\pr{\nabla\bar u^h+P_f\bs p}.
	\end{aligned}
\end{equation}
where

\textbf{-} $V(f)$ denotes the vertices of face $f$;

\textbf{-}  $\bar{u}^h=\sum_{i\in\mathcal{V}}\bar{u}_{i}\phi_i$ is the piecewise linear interpolation of the solution on the mesh vertices, computed in Section~\ref{sec:dis-adc-form};

\textbf{-} $P_f\in\mathbb R^{2\times3}$ is the projection matrix consisting of two orthonormal basis vectors on face $f$, i.e., $P_f = (\bs{g}_1, \bs{g}_2)^\top$;

\textbf{-}  $\bs{b}_f \in \mathbb{R}^{2 \times 2}$ is the second fundamental form of face $f$ under the basis $\{\bs{g}_1, \bs{g}_2\}$, approximated using a finite difference scheme~\cite{Rusinkiewicz2004,ads}.

The time derivative of the surface area is similarly discretized as:
\begin{equation}
	\label{eq:area-time-rate}
	\dt{A}\approx -\sum_{f\in\mathcal{F}}\sum_{i\in N(f)}v_n^i\int_{f}\phi_i\tr{\bs b_f}.
\end{equation}
where $N(f)$ denotes the vertices of face $f$.

Both integrals in~\eqref{eq:Ic-grad-h} and~\eqref{eq:area-time-rate} are evaluated using single-point quadrature per face.
We assemble the coefficients of $v_n^i$ in the discretized version of~\eqref{eq:common-ad-rate} into a vector $\bs{G}_{n}$, considered as the discrete $L^2$ gradient of $\ka_A$.

\begin{figure}[t]
	\centering
	\begin{overpic}[width=0.7\textwidth,keepaspectratio]{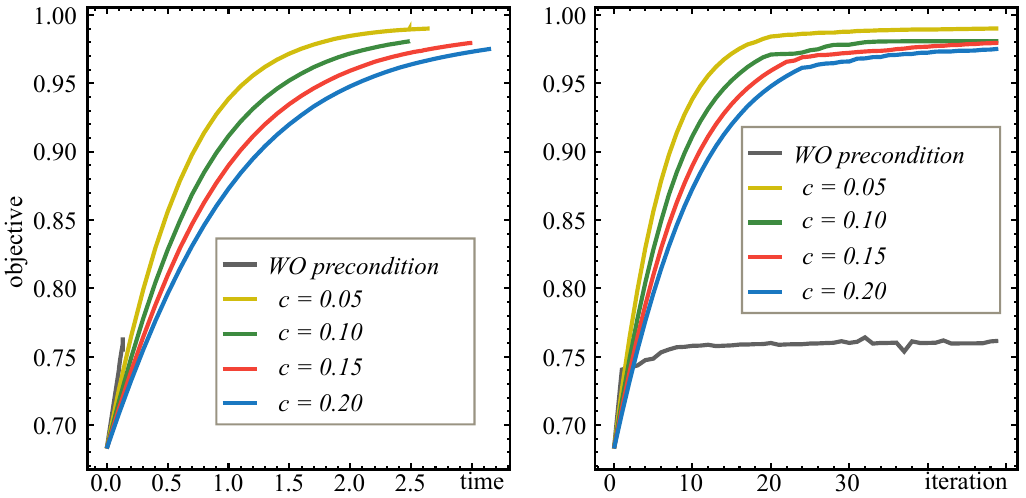}
	\end{overpic}
	\caption{
		Objective variation w.r.t. flow time (left) and iteration (right) when optimizing ADC along $x$-axis without precondition and with different precondition strength.
		The decrease of the objective arises from the failure of line search.
	}
	\label{fig:precon-c}
\end{figure}

\paragraph*{Custom objective}
To optimize a custom objective of the ADC matrix $f(\k_A)$, one must compute  its sensitivity.
Reusing the derivation from Proposition~\ref{prop:adc-sens}, we obtain the following result:
\begin{proposition}
	The time derivative $\dot{k}_A^{ij}$ is given by
	\begin{equation}
		\dot{k}_A^{ij}=\frac{\dot{I}_{ij}^c}{A}-\frac{I_{ij}^c\dot{A}}{A^2},
	\end{equation}
	where $I_{ij}^c$ denotes the integral in~\eqref{eq:bar-k-A-expr}, and its time derivative is 
	\begin{equation}
		\dot{I}^c_{ij}=\int_{\tw} 2v_n\pr{\bs b- H \bs I}\pr{\nabla \bar{u}^i+\e^i_{\tw}}\cdot\pr{\nabla \bar{u}^j+\e^j_{\tw}}.
	\end{equation}
\end{proposition}

The derivation is similar to Proposition~\ref{prop:adc-sens} with minor changes on the integrand, and is thus omitted here.
The corresponding discretization is given by
\begin{equation}
		\label{eq:Ic-ij-grad-h}
		\begin{aligned}
			\dot{I}^c_{ij} 
			&\approx\sum_{f\in\mathcal{F}}\sum_{q\in N(f)}2v_n^q\int_{f}\phi_q\pr{\nabla\bar u^h_i+P_f\e_i}^\top\pr{\bs b_f -\frac{\tr{\bs b_f}}{2}\mathbf I}\pr{\nabla\bar u^h_j+P_f\e_j},
		\end{aligned}
\end{equation}
where $\bar u_i^h$ is the solution of~\eqref{eq:u-div-eq-hp3} when $\p=\e^i$, computed in Section~\ref{sec:dis-adc-form}.

\begin{figure}[t]
	\centering
	\begin{overpic}[width=0.7\textwidth,keepaspectratio]{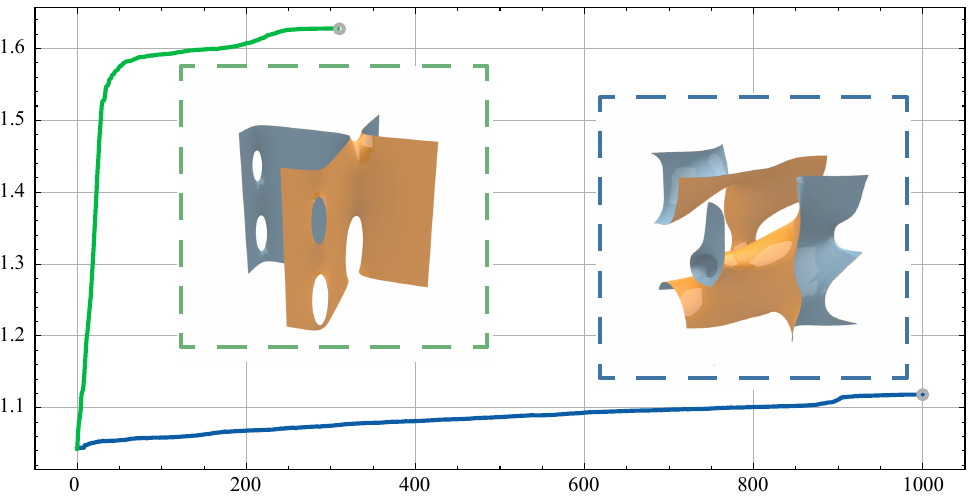}
		\put(25,40){\small W prec.}
		\put(70,37){\small WO prec.}
		\put(48,-2){\small iteration}
		\put(-4,20){\small \rotatebox{90}{objective}}
	\end{overpic}
	\vspace{2mm}
	\caption{
		Optimizing objective $k_A^{33}+\ka_A^{Avg}$ with/without precondition after remesh operation is introduced in each iteration.
	}
	\label{fig:precon-opt}
\end{figure}
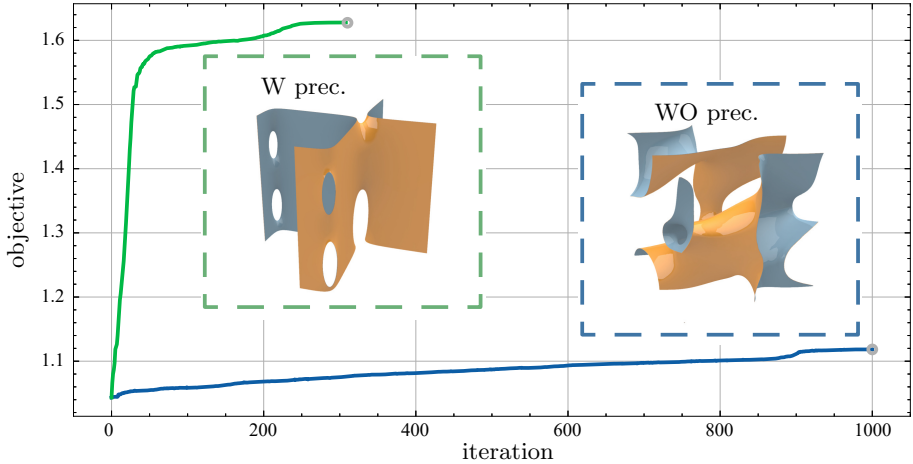

\subsubsection{Optimization}
As noted in~\cite{ads}, the $L^2$ gradient flow~\eqref{eq:l2-grad-flow} often leads to very small time steps during line search, even when initialized with a large step size.
This eventually stalls the optimization (Figure~\ref{fig:precon-c} right), and is difficult to recover from—even with remeshing (Figure~\ref{fig:precon-opt}).
To address this issue, we adopt the preconditioning strategy proposed in~\cite{ads}. Specifically, we use a preconditioned gradient $\bs{d}_n$ obtained by solving the following linear system:
\begin{equation}
	\label{eq:precon-eq}
	\pr{\mathbf M - c\mathbf L} \bs d_n =  (c+1)\bs g_n,
\end{equation}
where $\mathbf{M}$ is the mass matrix, $\mathbf{L}$ is the cotangent Laplacian, and $c \in [0, \infty)$ denotes the preconditioning strength.
While preconditioning slows down the objective change with respect to flow time  $t$ in \eqref{eq:l2-grad-flow} (Figure~\ref{fig:precon-c} left), it enables significantly larger time step in each iteration, which greatly improves convergence speed in terms of iteration (Figure~\ref{fig:precon-c} right).

With the preconditioned gradient,  we update the vertices' position ($\bs x_i$) by
\begin{equation}
	\label{eq:update-vertex}
	\bs{x}_i^{k+1} = \bs{x}_i^k +\Delta t d_i\n_i,
\end{equation}
where $d_i$ is the $i$-th component of $\bs{d}_n$, and time step $\Delta t$ is determined through Armijo line search.
In practice, we optionally add a diffusion term $w \mathbf{L} \bs{x}$ to regularize the step direction for mesh fairing, with $w = 0.1$ in our experiments.

Our optimization pipeline  proceeds as follows in each iteration:
\begin{enumerate}
	\item Detect singularities and perform numerical surgery;
	\item Apply dynamic remeshing;
	\item Evaluate the ADC matrix $\k_A$, objective $f^i = f(\k_A)$, and shape gradient $\bs{g}_n$ (Sections~\ref{sec:dis-adc-form} and~\ref{sec:sens-ana});
	\item Solve for the preconditioned gradient $\bs{d}_n$ using~\eqref{eq:precon-eq};
	\item Perform line search: $\Delta t \gets ArmijoLineSearch(f^i,\bs g_n,\bs d_n)$;
	\item Update vertex positions using~\eqref{eq:update-vertex};
	\item Check convergence; if satisfied, output the mesh; otherwise, return to step 1.
\end{enumerate}
The surgery and remeshing steps (steps 1 and 2) are adopted from~\cite{ads} and are essential for addressing necking singularities (Figure~\ref{fig:surj-seq}) and preserving mesh quality—both of which are known to impede optimization.
The optimization is considered converged (step 7) when either of the following criteria is met: (i) the time step $\Delta t$ drops below $10^{-4}$ for five consecutive iterations; or (ii) the time derivative of the objective function, estimated via linear regression over the last 50 iterations, falls below $10^{-3}$.

	\section{Validation}
\label{sec:valid}
\subsection{Preparation}
\label{sec:vali-prep}

We  process and generate periodic surface meshes using the same method as described in Section 7 of~\cite{ads}. 
A total of 24 types of TPMS are generated using Surface Evolver~\cite{surf-evolver}, and 1000 non-TPMS surfaces are created via a random perturbation process~\cite{ads}.
The source files of these surfaces, along with the associated experimental statistics, are available in the supplementary material~\cite{sm}. 

To evaluate the effective conductivity matrices of the resulting shell lattices, the surface meshes are imported into ABAQUS for finite element simulation.  
Each triangular facet of the surface mesh is assigned a DS3 shell element.
The thermal conductivity of the base material is set to $\ka=1$.
Periodic boundary conditions are enforced using the Micromechanics plugin, which also computes the effective conductivity matrix.

\subsection{Validating the convergence theorem}
\label{sec:val-conv-theorem}
We numerically approximate the limit in~\eqref{eq:kmat-eps-to-kmat-A} by the ratio $\k_\eps/\rho_\eps$ under sufficiently small volume fraction.
To this end, we set the shell thickness to $0.004$ in the simulation.
The effective conductivity matrix $\k_\eps$ is computed from the ABAQUS simulation, while the corresponding ADC matrix $\k_A$ is obtained following the procedure in Section~\ref{sec:dis-adc-form}.
We observe that the relative error
\begin{equation}
	\mathcal E^c:=\frac{\|\k_A-\k_\eps/\rho_\eps\|_F}{\|\k_A\|_F}
\end{equation}
for all the 1024 surfaces generated in Section~\ref{sec:vali-prep} satisfies $\mathcal E^c\le 0.3\%$ (see the statistic tables in the supplementary material~\cite{sm}).
Considering the discretization error inherent in both numerical methods, this result strongly supports the convergence predicted by~\eqref{eq:kmat-eps-to-kmat-A}.

\subsection{Convergence rate of the discretization}
\label{sec:dis-conv-rate}
\begin{figure}[t]
	\centering
	\begin{overpic}[width=0.5\textwidth,keepaspectratio]{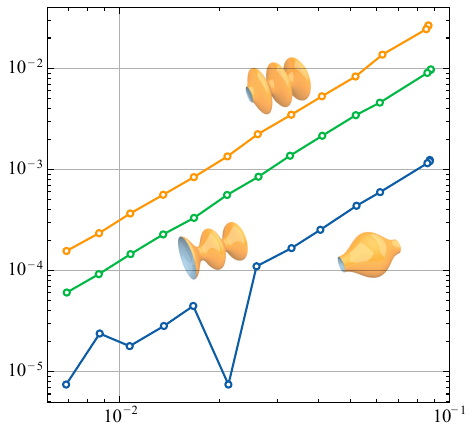}
		\put(54,0){$h$}
		\put(-7,50){\rotatebox{0}{$\mathcal{E}^h$}}
	\end{overpic}
	\caption{
		Discretization error ($\mathcal E^h$) with different element size ($h$).
	}
	\label{fig:convadc}
\end{figure}

We investigate the convergence rate of our discretization scheme in Section~\ref{sec:dis-adc-form}, using a set of revolution surfaces for which the ADC along the medial axis admits an analytical expression (\ref{app:exp-form-adc-revo}).
The surfaces are constructed by revolution around the $x$-axis with radius profile $R(x)=(2+\cos \pi x)/4$.
Following the setup in~\cite{ads}, we refine the mesh by remeshing with progressively decreasing target edge lengths, and compute the relative error with respect to the analytical value:
\begin{equation}
	\mathcal{E}^h = \frac{|\ka_A-\ka_A^h|}{\ka_A},
\end{equation}
where $\ka_A^h$ denotes the numerically computed ADC on a mesh with element size $h$, defined as the average circumradius of the faces.
The change of errors is plotted in Figure~\ref{fig:convadc}, which demonstrates a second-order convergence rate.

\subsection{Validating third-order accuracy of ADC}
\label{sec:3order-acc}
Although Section~\ref{sec:val-conv-theorem} validates the convergence theorem, these simulation results cannot be used to verify the third-order accuracy of ADC (Theorem~\ref{thm:adc-3-order}) directly, as the convergence order of the DS3 element discretization is not theoretically established.
To rigorously validate Theorem~\ref{thm:adc-3-order}, we instead employ a semi-analytical method (see~\ref{app:semi-any-eval}) to evaluate the effective conductivity of the shell lattices constructed from the revolution surfaces in Section~\ref{sec:dis-conv-rate}.

\begin{figure}[t]
	\centering
	\begin{overpic}[width=0.6\textwidth,keepaspectratio]{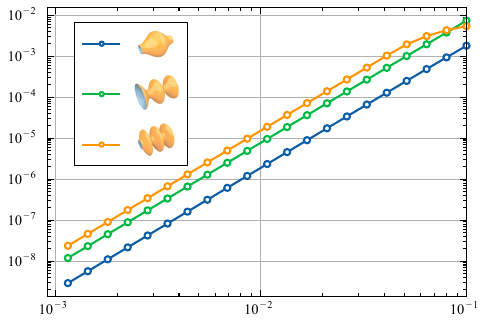}
		\put(30,1){$\epsilon$}
		\put(-3,37){$E_r$}
	\end{overpic}
	\caption{
		Convergence of the residual term in Theorem~\ref{thm:adc-3-order} as thickness approaches zero.
	}
	\label{fig:adc2order}
\end{figure}

We define the residual term as:
\begin{equation}
	E_r = \left|\ka_\eps -\frac{2\eps|\tw|}{|Y|}\ka_A\right|,
\end{equation}
and examine its behavior with respect to the variation of shell thickness, as shown in Figure~\ref{fig:adc2order}.
The results exhibit clear third-order convergence for all three tested surfaces, thereby confirming the prediction of Theorem~\ref{thm:adc-3-order}.

\begin{figure}[t]
	\centering
	\begin{overpic}[width=0.6\columnwidth,keepaspectratio]{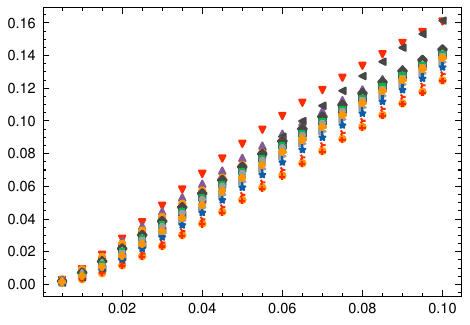}
		\put(50,-1) {$2\eps$}
		\put(-4,35) {\rotatebox{90}{\small $\mathcal{E}_a$}}
	\end{overpic}
	\caption{
		The errors of the AACs to the theoretical upper bound for the selected 24 types of TPMS as thickness approaches zero.
	}
	\label{fig:tpms-apac-conv}
\end{figure}

\subsection{Validating the optimal ADC of TPMS}
To verify the conclusion of Theorem~\ref{thm:tpms-opt-adc}, we evaluate the AAC ($\ka_A^{\mathrm{Avg}}$) for the TPMS generated in Section~\ref{sec:vali-prep}, using the identity $\ka_A^{\mathrm{Avg}} = \mathrm{tr}(\k_A)/3$, where $\k_A$ is approximated by $\k_\eps / \rho_\eps$.
For all generated TPMS, the relative errors to the upper bound, defined by 
\begin{equation}
	\mathcal E_a := 1-\frac{1}{3}\pr{\tr\k_\eps/\rho_\eps}/\ka_A^*, 
\end{equation}
consistently decrease to zero (Figure~\ref{fig:tpms-apac-conv}) as the thickness approaches zero.

\begin{figure}[t]
	\centering
	\begin{overpic}[width=0.6\textwidth,keepaspectratio]{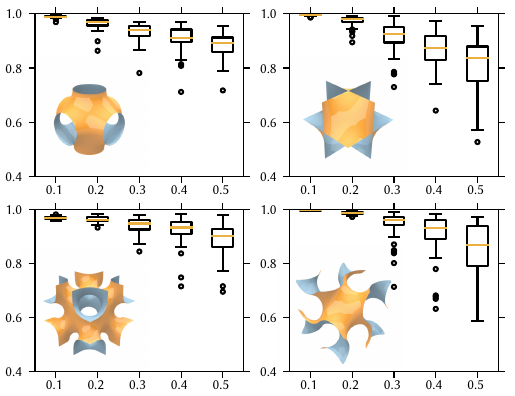}
		\put(-5,54) {\rotatebox{90}{\small$\bar\ka_A^{Avg}/ \bar\ka_A^*$}}
		\put(-5,16) {\rotatebox{90}{\small$\bar\ka_A^{Avg}/ \bar\ka_A^*$}}
		\put(27,-2) {$\bar s$}
		\put(77,-2) {$\bar s$}
	\end{overpic}
	\caption{
		The AAC of perturbed surfaces derived from Schwarz P, Schwarz D, IWP, and Gyroid. Each subplot corresponds to one TPMS type.
		The abscissa denotes the perturbation strength~\cite{ads}, with 50 surfaces sampled for each setting.
	}
	\label{fig:sample-surf}
\end{figure}

On the other hand, we plot the ratio $\mathrm{tr}(\k_A)/(3\ka_A^*)$, which is the approximation of $\ka_A^{Avg}/\ka_A^*$, for non-TPMS surfaces in Figure~\ref{fig:sample-surf}, with a fixed thickness $2\eps=0.004$.
The results indicate that this ratio remains strictly less than one for all such surfaces.
Moreover, as the strength of geometric perturbation increases—i.e., as the surfaces deviate further from the TPMS family—the ratio tends to decrease.
These numerical results are consistent with our theoretical prediction and further support the optimality of TPMS in terms of ADC.

\subsection{Validation from experimental data}
We further compare our theoretical predictions with experimental results reported in~\cite{Zhou2300359}, where the heat transfer properties of 3D printed TPMS structural zirconia  ceramics  are investigated.
The base material used in their study, $\text{ZrO}_2$, has a thermal conductivity of  $\ka = 2.249 \wmk$.
For the Schwarz P surface, they report that the derived shell lattice exhibits an effective thermal conductivity of $\ka_\eps(\tw) = 0.2262 \wmk$ at a volume fraction of $\rho_\eps(\tw)=14.7\%$.  
Notably, this measured value is very close to the corresponding Hashin–Shtrikman upper bound, given by $ \ka^{HS}_\eps(\tw) = 2\rho_\eps(\tw)\ka/(3-\rho_\eps(\tw))=0.2318\wmk$,
which is consistent with our theoretical prediction.
Similar agreement is observed for other TPMS structures reported in their study (see Table~\ref{tab:tpms-exp-data}).
Our analysis predicts that  $\ka_\eps / (f_\eps\ka)$ should approaches $2/3$ when $\rho_\eps$ is small, and the experimental results show that this ratio indeed remains close to $2/3$ across all tested samples.
\begin{table}[h]
	\centering
	\caption{
		\label{tab:tpms-exp-data}
		Volume fractions and conductivities of the shell lattices derived from 4 types of TPMS.
		The data in first three columns are from~\cite{Zhou2300359}.
		For each type of TPMS, three specimens are presented, corresponding to physical cell (RVE) sizes of 2 mm, 3 mm, and 4 mm, respectively. 
	}
	\vspace{2mm}
	\resizebox{0.6\columnwidth}{!}{
		\begin{threeparttable}
			\begin{tabular}{lcccr}
				\multicolumn{1}{c}
				{\textrm{Type}} &  $\rho_\eps(\%)$ & $ \ka_\eps (\nswmk)$ & $ \ka^{HS}_\eps (\nswmk)$ & $\ka_\eps/(\rho_\eps\ka)$\\
				\midrule
				Schwarz P  & 14.7 & 0.2262 & 0.2318 & 0.6842\\
				& 14.5 & 0.2241 & 0.2284 & 0.6872\\
				& 13.6 & 0.2194 & 0.2135 & 0.7173\\
				Schwarz D  & 17.3 & 0.2280 & 0.2753 & 0.5860\\
				& 15.9 & 0.2167 & 0.2517 & 0.6060\\
				& 14.2 & 0.2109 & 0.2235 & 0.6603\\
				IWP & 15.9 & 0.2148 & 0.2517 & 0.6007 \\
				& 14.7 & 0.2100 & 0.2318 & 0.6352\\
				& 14.5 & 0.2104 & 0.2284 & 0.6452\\
				Gyroid & 14.7 & 0.2179 & 0.2318 & 0.6591\\
				& 14.5 & 0.2106 & 0.2284 & 0.6458\\
				& 13.9 & 0.2070 & 0.2185 & 0.6621\\
				\bottomrule
			\end{tabular}
		\end{threeparttable}
	}
\end{table}

\section{Experimental results}
\label{sec:app}

\subsection{Tailoring ADC}
\begin{figure}[t]
	\centering
	\begin{overpic}[width=0.6\textwidth,keepaspectratio]{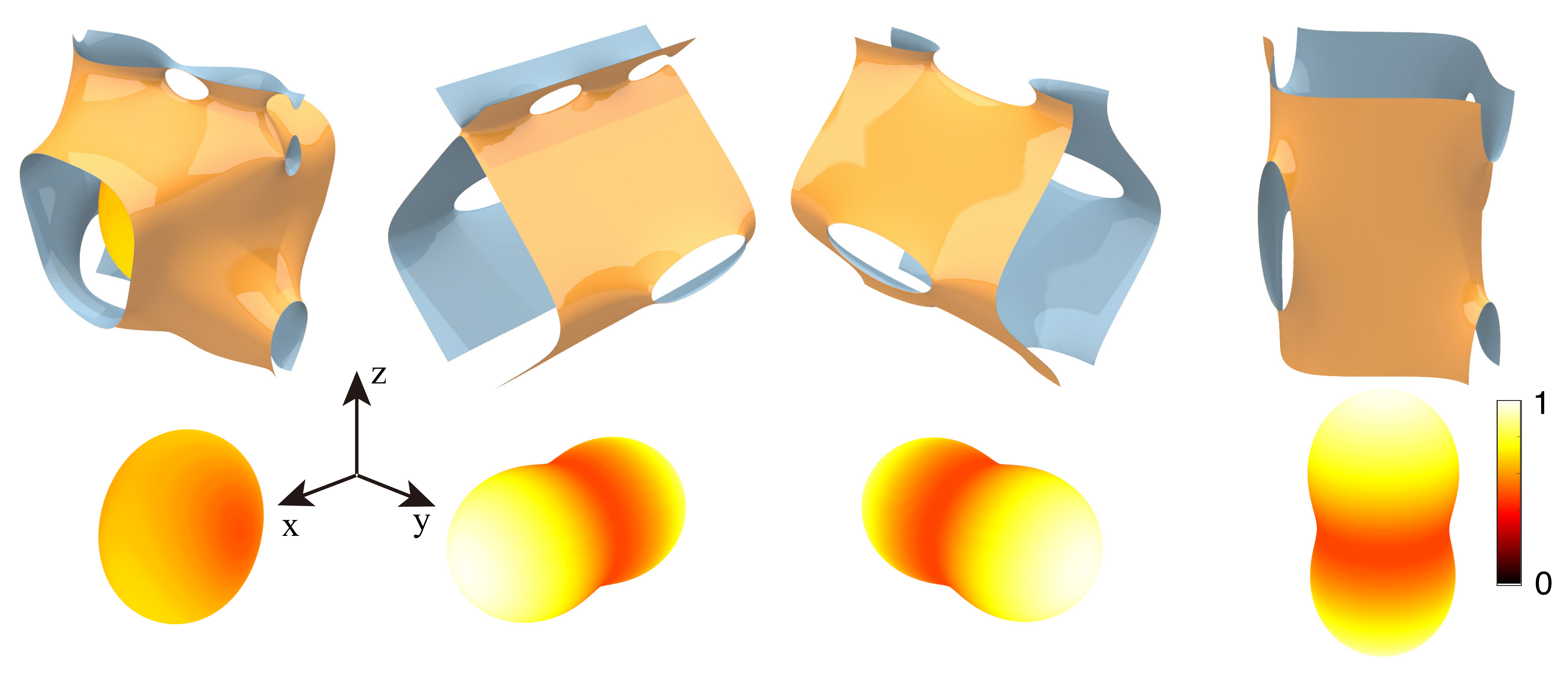}
		\put(7,-2){Input}
		\put(28,-2){$k_A^{11} + \ka_A^{Avg}$}
		\put(55,-2){$k_A^{22} + \ka_A^{Avg}$}
		\put(80,-2){$k_A^{33} + \ka_A^{Avg}$}
	\end{overpic}
	\vspace{2mm}
	\caption{
		Optimize ADC with same input surface (left) and different objectives (right 3 columns).
		The directional distribution of ADC is visualized below each surface.
	}
	\label{fig:adc-max-iwp}
\end{figure}
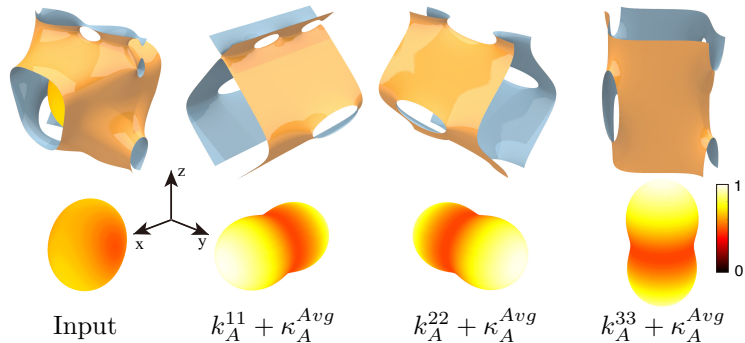

\paragraph*{Maximizing  ADC }
\label{sec:max-adc}

\begin{figure}[t]
	\centering
	\begin{overpic}[width=0.6\textwidth,keepaspectratio]{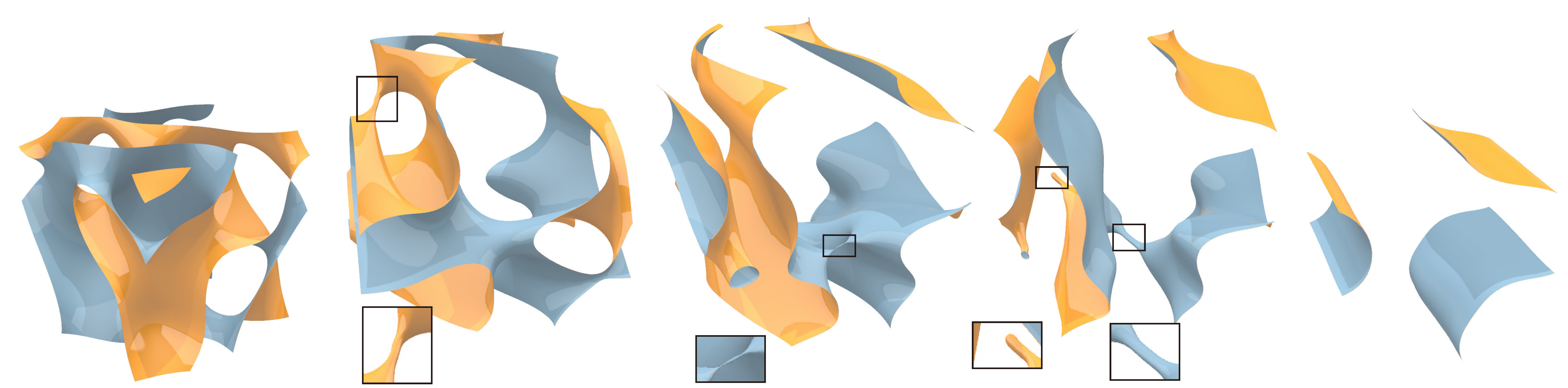}
		\put(7,-4) {Input}
		\put(85,-4){Final result}
	\end{overpic}
	\vspace{4mm}
	\caption{
		A surface evolves into a cylinder when optimizing $k_A^{22}$.
		The surface developed a sequence of singularities during the flow, which are all removed through surgery operation, and finally orients  to the objective direction.
	}
	\label{fig:surj-seq}
\end{figure}

We test our optimization pipeline using the following three objectives  
\begin{align}
	\label{eq:obj-f1}
	f_{i} =  k_A^{ii} + \ka_A^{Avg},\quad i=1,2,3,
\end{align}
where each objective emphasizes conductivity in one axis direction.
The resulting optimized surfaces exhibit clear directional bias aligned with the corresponding axis  (Figure~\ref{fig:adc-max-iwp}).
The additional term $\ka_A^{Avg}=\tr\,\k_A/3$ is used to prevent the surface from degenerating into a cylinder, which would compromise connectivity in directions orthogonal to the targeted axis (Figure~\ref{fig:surj-seq}).

\begin{figure}[t]
	\centering
	\begin{overpic}[width=0.6\textwidth,keepaspectratio]{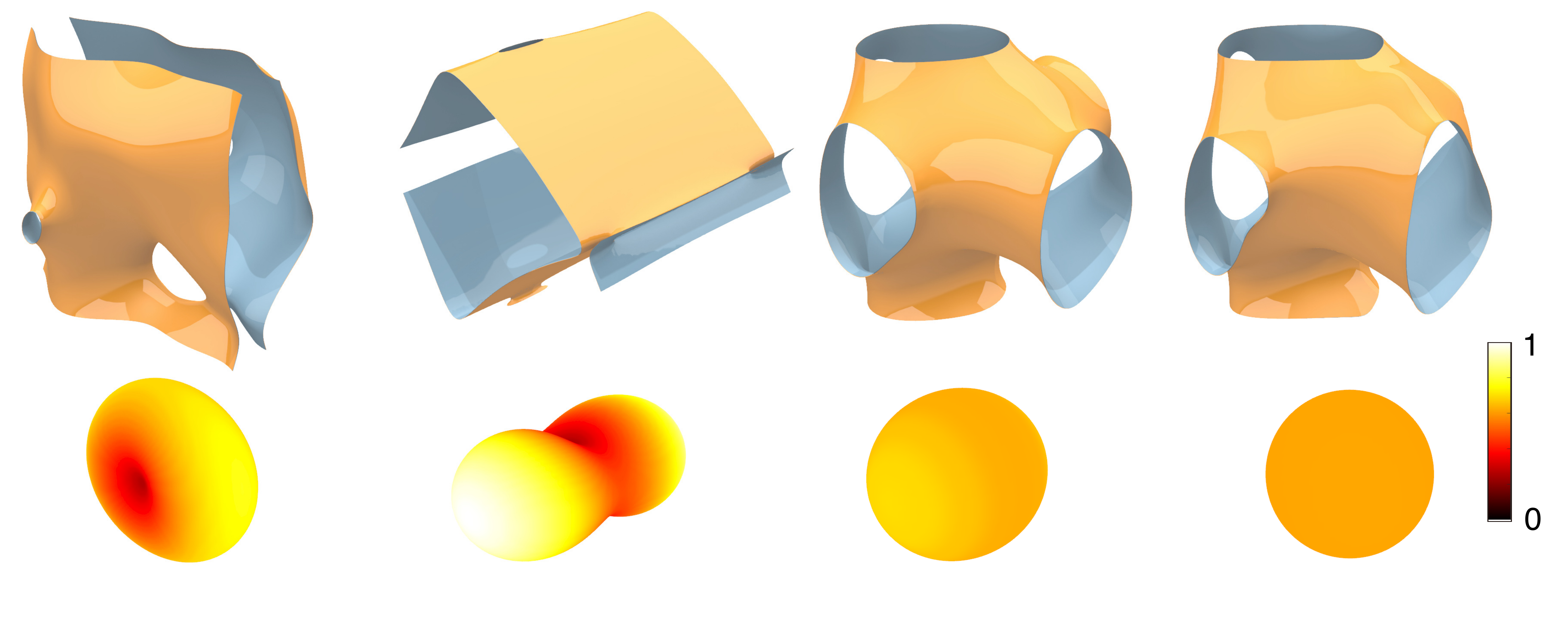}
		\put(7,-2){Input}
		\put(34,-2){$k_A^{11}$}
		\put(53,-2){$k_A^{11}-4f_{iso}^c$}
		\put(85,-2){$-f_{iso}^c$}
	\end{overpic}
	\vspace{2mm}
	\caption{
		Maximizing the conductivity with isotropic penalty.
		The objective is shown below each column.
		Top: the geometry of surfaces;
		Bottom: the directional distribution of ADC.
	}
	\label{fig:adc-iso}
\end{figure}

\paragraph*{Isotropic penalty}
An isotropic conductivity matrix is a scaled identity matrix, whose eigenvalues are equal.
To encourage isotropy in the optimized ADC matrix, we introduce the following penalty term:
\begin{equation}
	f_{iso}^c(\bs k_A) = \lambda_{max}(\bs k_A) -\lambda_{min}(\bs k_A)
\end{equation}
where $\lambda_{\max}$ and $\lambda_{\min}$ denote the largest and smallest eigenvalues of $\bs{k}_A$, respectively.
The time derivative of this function is given by
\begin{equation}
	\frac{d}{dt} f_{iso}^c(\bs k_A) = \bs p_{max}^\top \dt{\bs k}_A \bs p_{max}-\bs p_{min}^\top\dt{\bs k}_A \bs p_{min}
\end{equation}
where $\bs p_{max}$ and $\bs p_{min}$ are the normalized eigenvectors  associated with $\lambda_{\max}$ and $\lambda_{\min}$, respectively.
After adding this penalty term to the optimization objective, the ADC distribution of the optimized surface exhibit significantly improved isotropy (Figure~\ref{fig:adc-iso}).

\begin{figure}[t]
	\centering
	\begin{overpic}[width=0.8\textwidth,keepaspectratio]{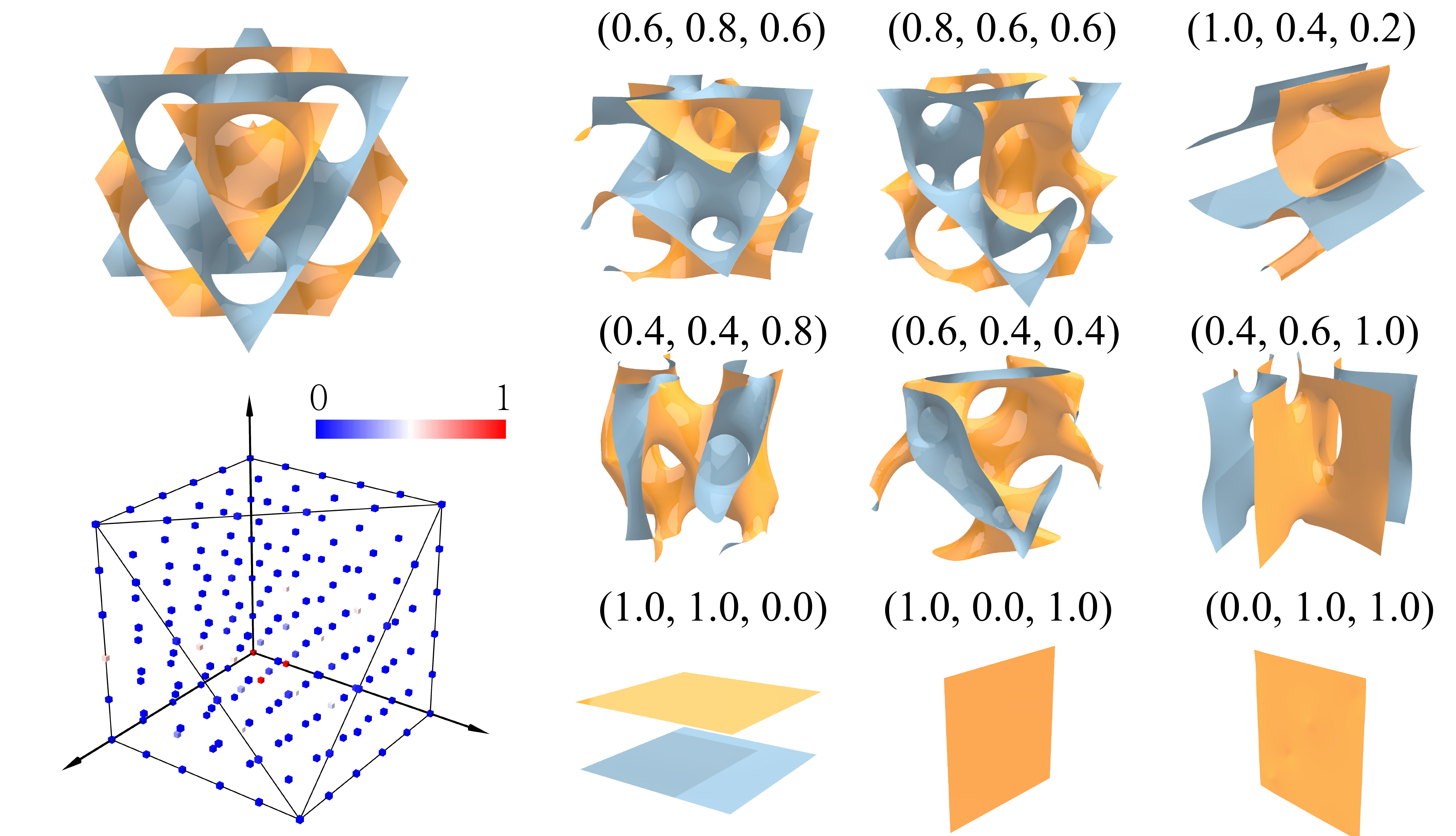}
	\put(5,32){(a)}
	\put(5,-2){(b)}
	\put(68,-2){(c)}
	\put(26,30){$\mathcal{E}_{tgt}$}
	\put(7,4){\small 1}
	\put(29,6){\small 1}
	\put(17.5,26){\small 1}
	\put(2,2){\small $\ka_{1}$}
	\put(34,6){\small $\ka_{2}$}
	\put(13,28){\small $\ka_{3}$}
	\end{overpic}
	\vspace{2mm}
	\caption{
		(a) The input surface for optimizing toward specified ADC matrices.
		(b) Sampling points (dots) within the feasible region of $\k_A$, where color indicates the relative error $\mathcal{E}_{\mathrm{tgt}}$ between the optimized and target ADC matrices.
		(c) The optimized surfaces with target principle conductivity shown above.
	}
	\label{fig:adctex}
\end{figure}

\paragraph{Specified ADC matrix}
\label{sec:speci-ADC}
 To design surfaces whose ADC matrix matches a prescribed target, we minimize the following objective function:
\begin{equation}
	\label{eq:tgt-ka-obj}
	f^c(\k_A) = \|\k_A-\k^*\|_F
\end{equation}
where $\k^*$ denotes the target ADC matrix and $\|\cdot\|_F$ is the Frobenius norm.
The feasible region for  $\k_A$, as implied by Theorem~\ref{thm:adc-upper-bound} and~\ref{thm:tpms-opt-adc}, is given by
\begin{equation}
	\label{eq:ka-feas-region}
	\left\{ 
	\begin{aligned}
		&0\le\ka_i \le 1,\, i=1,2,3\\
		&\ka_1+\ka_2+\ka_3\le 2
	\end{aligned}
	\right.,
\end{equation}
where $\ka_i$ are the eigenvalue of $\k_A$.
We sample $\k^*$ in this region by setting diagonal entries $k_A^{ii} = m \Delta \kappa$ and off-diagonal entries $k_A^{ij} = 0$ for $i \ne j$, where $\Delta\ka=0.2$ is a step and $m$ is a positive integer.
For each sampled target, we optimize the objective~\eqref{eq:tgt-ka-obj} starting from a common input surface (Figure~\ref{fig:adctex}(a)).
The relative errors between the optimized and target ADC matrices, given by
\begin{equation}
	\mathcal{E}_{tgt}=f^c(\k_A)/\|\k^*\|_F,
\end{equation}
 are plotted in Figure~\ref{fig:adctex}(b).
 Several optimized surfaces are visualized in Figure~\ref{fig:adctex}(c).
 The results indicate that, for most cases, the optimized ADC matrices closely match the targets, with $\mathcal{E}_{{tgt}} \approx 0$.

\begin{figure}[t]
	\centering
	\begin{overpic}[width=0.7\textwidth,keepaspectratio]{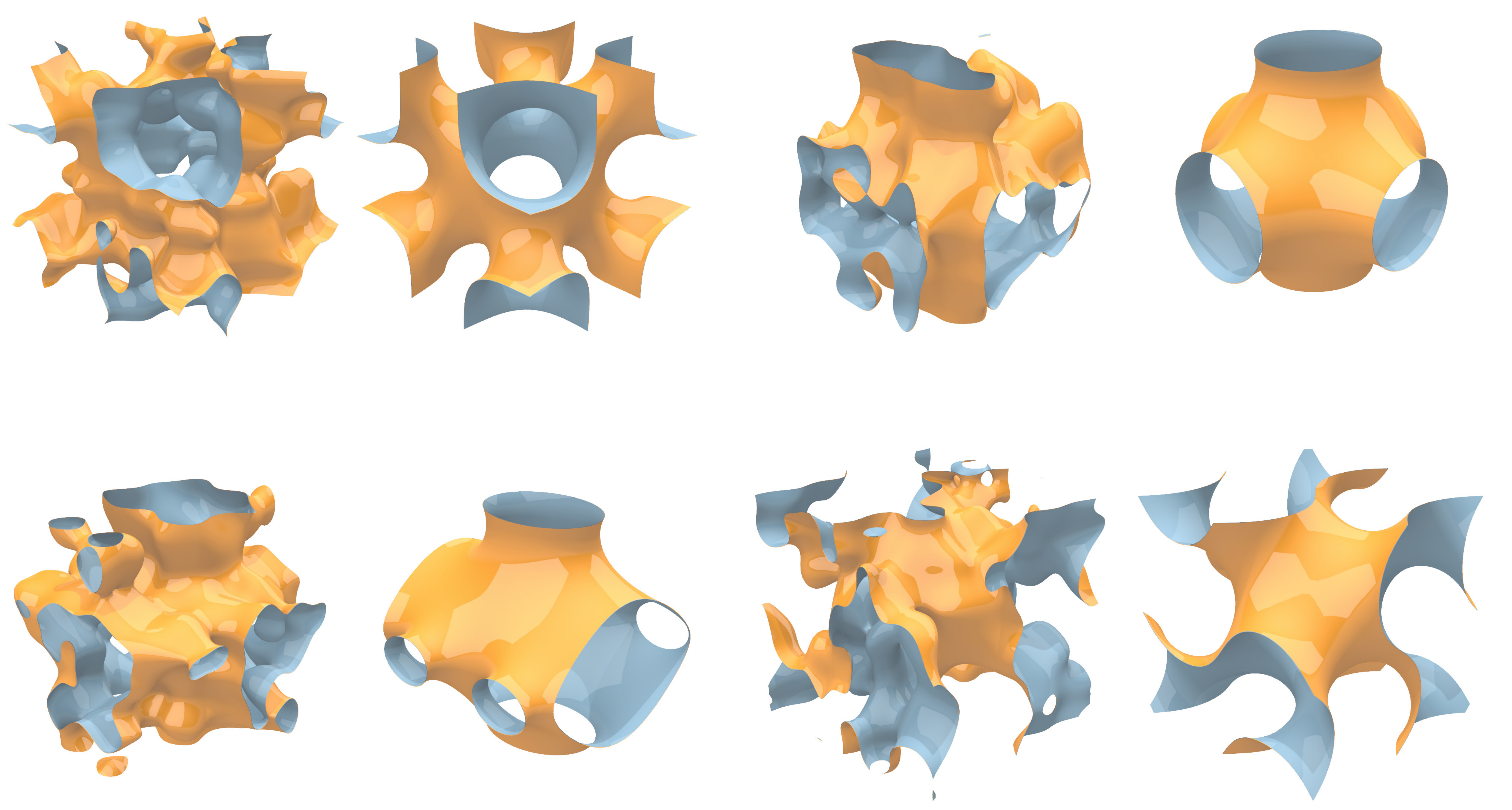}
		\put(52,-2){\small $\ka_A^{Avg}=0.4875$}
		\put(77,-2){\small$\ka_A^{Avg}=0.6664$}
		\put(52,27){\small $\ka_A^{Avg}=0.4763$}
		\put(77,27){\small$\ka_A^{Avg}=0.6666$}
		\put(2,-2){\small $\ka_A^{Avg}=0.4547$}
		\put(27,-2){\small$\ka_A^{Avg}=0.6666$}
		\put(2,27){\small $\ka_A^{Avg}=0.4870$}
		\put(27,27){\small$\ka_A^{Avg}=0.6666$}
		\put(8,53){\small  Input}
		\put(31,53){\small  Result}
		\put(58,53){\small  Input}
		\put(81,53){\small  Result}
	\end{overpic}
	\vspace{2mm}
	\caption{
		Optimizing $\ka_A^{Avg}$ yields TPMS shaped surfaces.
		The values of objective  are noted below  surfaces.
	}
	\label{fig:apac}
\end{figure}

		\begin{figure}[t]
				\centering
				\begin{overpic}[width=0.7\textwidth,keepaspectratio]{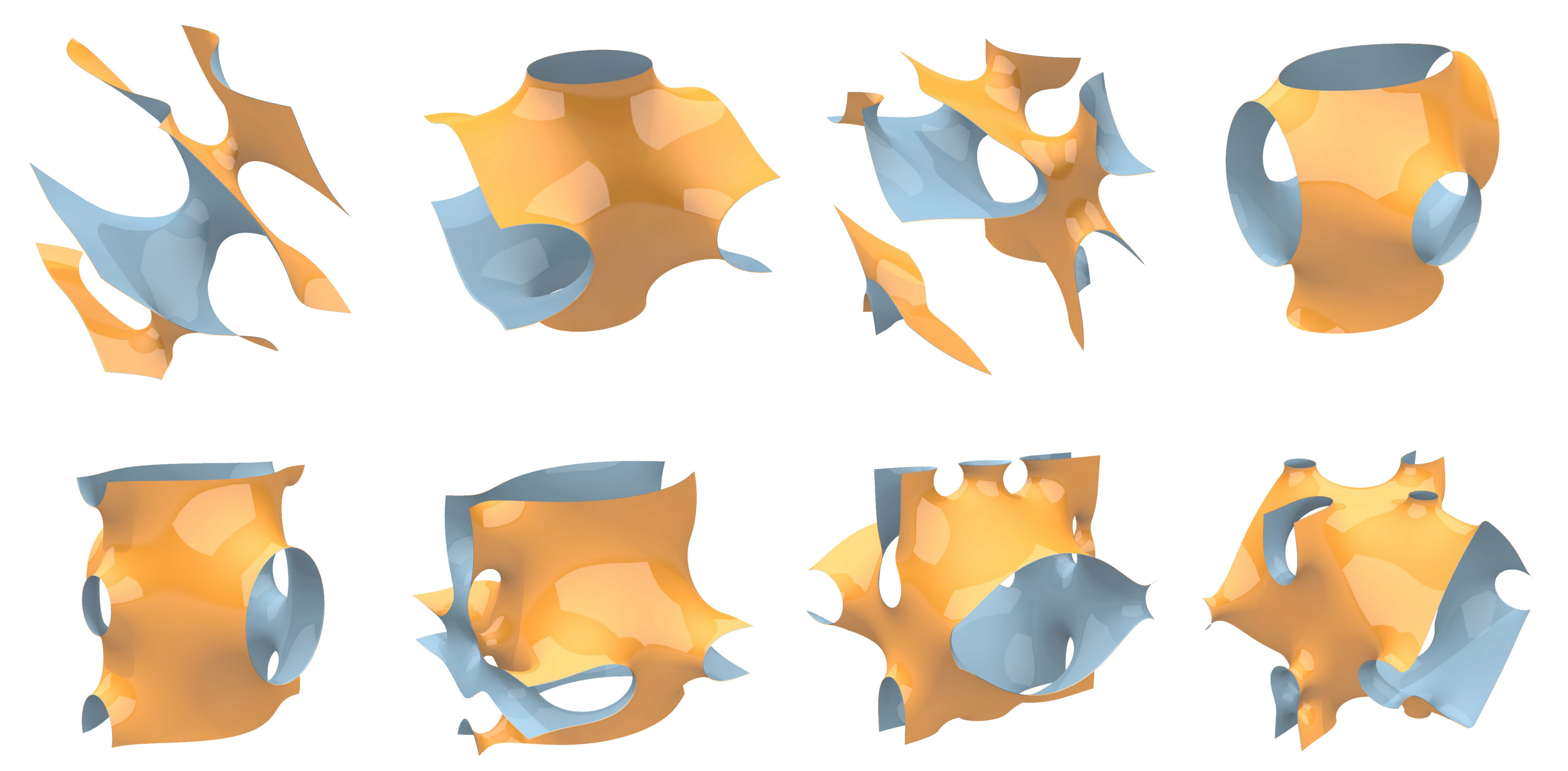}
					\end{overpic}
				\caption{
		Gallery of TPMS generated by optimizing $\ka_A^{Avg}$.
					}
				\label{fig:tpms-gal}
			\end{figure}

\begin{figure}[t]
	\centering
	\begin{overpic}[width=0.7\textwidth,keepaspectratio]{mcf-vs-apac-c.pdf}
		\put(-4,0){\small $A$}
		\put(2,0){\small $4.0010$}
		\put(20,0){\small $3.4650$}
		\put(37,0){\small $2.6184$}
		\put(53,0){\small $1.8346$}
		\put(70,0){\small $1.0042$}
		\put(87,0){\small $0.7908$}
		\put(-6,20){\small $\ka_A^{Avg}$}
		\put(2,19){\small $0.3889$}
		\put(19,19){\small $0.4489$}
		\put(35,19){\small $0.5159$}
		\put(52,19){\small $0.6012$}
		\put(69,19){\small $0.6632$}
		\put(86,19){\small $0.6666$}
	\end{overpic}
	\vspace{0mm}
	\caption{
		Comparison with mean curvature flow.
		Each row shows the surface evolution over iterations, from left to right.
		Top: optimizing AAC; Bottom: mean curvature flow.
		AAC values and areas are annotated below the corresponding surfaces.
	}
	\label{fig:mcf-vs-apac-c}
\end{figure}

\subsection{Generating TPMS by optimizing AAC}
\label{sec:gen-tpms-aac}
According to Theorem~\ref{thm:adc-upper-bound}, we can generate TPMS by optimizing AAC ($\ka_A^{Avg}$).
To validate this observation, we apply our optimization algorithm to a variety of input surfaces  with the objective of maximizing  $\ka_A^{Avg}$.
 As shown in Figure~\ref{fig:apac}, all input surfaces evolve into TPMS-like geometries, with AAC values significantly increased and approaching the theoretical upper bound $(\frac{2}{3})$.
 More results are shown in Figure~\ref{fig:tpms-gal}.
In contrast to classical conjugate surface construction methods~\cite{Makatura-proce}, which require manually specified boundary constraints that must satisfy certain compatibility conditions, our approach operates on arbitrary periodic surfaces without any additional prior geometric input.
This drawback also exhibits in the GMT based method~\cite{Stephanie2021,Palmer9879635}, where a boundary must be specified.

It is well known that applying MCF to generate TPMS often results in degenerate geometry (Figure~\ref{fig:mcf-vs-apac-c} bottom).
This is expected, as MCF minimizes surface area, whereas TPMS do not minimize area within a fixed topological class.
Consequently, the surface area may consistently decrease during the flow, eventually leading to collapse (Figure~\ref{fig:mcf-vs-apac-c} bottom).
In contrast, our method ensures that once the AAC reaches its optimal value (i.e., the upper bound $\kappa_A^*$), the resulting surface is guaranteed to be a TPMS (Figure~\ref{fig:mcf-vs-apac-c} top). This property provides both robustness and geometric fidelity, making our approach significantly more reliable than MCF-based alternatives.

Although optimizing the asymptotic bulk modulus also yields TPMS with theoretical guarantees~\cite{ads}, as does AAC, the associated variational equation results in  a linear system approximately three times larger than that arising from AAC optimization.
Moreover, the sensitivity analysis is more intricate to implement and compute, resulting in greater computational overhead.

\begin{figure}[h]
	\vspace{3mm}
	\centering
	\begin{overpic}[width=0.6\columnwidth,keepaspectratio]{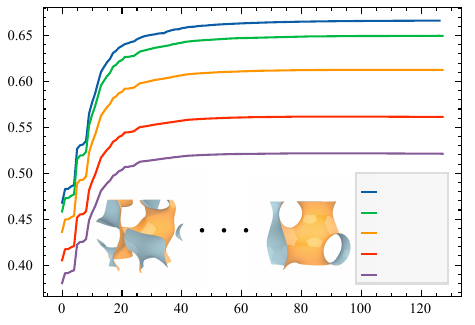}
		\put(80,27){\footnotesize $Asym.$}
		\put(80,23){\footnotesize $\eps=0.01$}
		\put(80,18.5){\footnotesize $\eps=0.025$}
		\put(80,14){\footnotesize $\eps=0.05$}
		\put(80,9.3){\footnotesize $\eps=0.075$}
		\put(50,-2){\footnotesize iteration}
		\put(25,28){\footnotesize Input}
		\put(60,28){\footnotesize Result}
		\put(-3,28){\footnotesize \rotatebox{90}{objective}}
	\end{overpic}
	\vspace{2mm}
	\caption{
		Variation of ADC (`Asym.') and normalized effective conductivity of shell lattice at different thicknesses $(\eps)$  during the optimization of AAC.
		The objective is evaluated as $\tr \k_A/3$ and $\tr \k_\eps/(3\rho_\eps)$, respectively.
	}
	\label{fig:effect-evolve-heat}
\end{figure}

\subsection{Evolution of effective conductivity during the optimization}
We compare the evolution of ADC and the effective conductivity of the corresponding shell lattices  at various thicknesses, throughout the optimization of $\ka_A^{\mathrm{Avg}}$.
The results (Figure~\ref{fig:effect-evolve-heat}) indicate that both quantities exhibit consistent trends over the course of the optimization.
Combined with the validation results in Section~\ref{sec:3order-acc}, this consistency  further supports  the effectiveness of our algorithm in the design of shell lattices with specified thermal conductivity.

\section{Conclusion}
\label{sec:concl}
We propose a novel theoretical framework for analyzing and designing shell-based thermal metamaterials by introducing the \emph{asymptotic directional conductivity} (ADC)—a metric that captures the directional conductivity of shell lattices as thickness approaches zero.
We prove a convergence theorem for computing ADC and establish its upper bound, which is attained when the surface is a TPMS. This provides the first rigorous justification for the superior thermal conduction of TPMS structures.
Although ADC is defined in the asymptotic sense, we show both theoretically and numerically that it approximates effective conductivity with third-order accuracy at low volume fractions. This makes ADC not only a theoretical tool but also a practical design objective. Our discretization scheme enables efficient evaluation of ADC on discrete triangular meshes, and our shape optimization algorithm allows tailoring shell lattices for various objectives.
Beyond thermal applications, the mathematical formulation of ADC extends naturally to other transport phenomena, such as electrical conductivity and diffusion, providing a new perspective for geometry-driven metamaterial design across disciplines.
In future work, we aim to explore the inverse design of shell lattice metamaterials under manufacturing constraints and extend our framework to coupled thermal-mechanical applications.

\section*{Acknowledgement}
This  work  is  supported  by  the  National Key  R\&D  Program  of  China  (2022YFB3303400)  and  the  National Natural Science Foundation of China (62025207).

\bibliographystyle{elsarticle-num} 
\bibliography{ref}       

	\appendix

\section{Auxiliary lemmas}
\begin{lemma}
	\label{lem:nec-cond-upp}
	If equality in~\eqref{eq:ka-upp-bound-p3} holds, then $\tw$ is either minimal or a cylinder parallel to $\p$.
\end{lemma}
\begin{proof}
	Assume~\eqref{eq:prove-equali-cond} equals zero, then we have $\nabla\bar u=0$, which implies that $\bar u$ is a constant.
	Thus we obtain $2Hp_3=\div \p_{\tw}=-\Delta\bar u=0$, indicating that either $H = 0$ or $p_3 = 0$ must hold on $\tw$.
	Let $F := \{ y \in \tw \mid H(y) \ne 0 \}$ denote the set of points with nonzero mean curvature. Due to the smoothness of $\tw$, $F$ is an open subset of $\tw$. We now consider two cases depending on whether $F$ is empty.
	
	\begin{itemize}
		\item 	Assume $F \ne \emptyset$. Let $\{F_k\}_{k \in I}$ be the collection of maximal connected open subsets of $F$, where $I$ is an index set. On each $F_k$, we have $H \ne 0$ and $H p_3 = 0$, thus $p_3 = 0$, meaning that each $F_k$ is a cylindrical surface parallel to $\p$.
		Note that on each ruling line (generator) of such a cylinder, the mean curvature is a nonzero constant that depends only on the curvature of the directrix. 
		
		We now prove that no generator of the cylinder intersects the boundary $\partial F_k$.
		Assume the contrary: there exists a point $x \in \partial F_k$ belongs to some ruling line such that $H(x) \ne 0$. Since $F_k$ is maximal, any neighborhood of $x$ must contain points where $H = 0$. By the smoothness of the surface, the mean curvature at $x$ must be continuous, and thus equal to 0 — a contradiction.
		Therefore, each generator of $F_k$ extends infinitely in the direction of $\p$, and so does $\partial F_k$.
		
		For the remainder of the surface, i.e., $\tw \setminus \bigcup_k F_k$, we have $H = 0$, implying that this region is a minimal surface with the boundaries $\{\partial F_k\}$, which are composed of straight lines parallel to $\p$.
		Since the normal vectors on the boundary of these minimal surfaces match those of the adjacent cylindrical regions $\overline{F_k} := F_k \cup \partial F_k$, the only possible surface satisfying both boundary position and boundary normal conditions is a plane. Hence, the entire surface $\tw$ must be a cylinder parallel to $\p$.
		\item 
		If $F = \emptyset$, then $H \equiv 0$, and thus $\tw$ is minimal.
	\end{itemize}
\end{proof}

\begin{lemma}
	\label{lm:cubic-srf-int}
	For a cubic symmetric periodic surface $\tw$ and any $\p\in \mathbb R^3$ with $\|\p\|=1$, we have 
	\begin{equation}
		\label{eq:lem-p3sq-13}
		\int_{\tw} p_3^2=\frac{1}{3}|\tw|
	\end{equation}
	where $p_3=\p\cdot\mathbf n$ denotes the normal component of $\p$ on surface $\tw$.
\end{lemma}
\begin{proof}
	The cubic symmetric group $G=\br{G_1,\cdots,G_{48}}$ is generated from 6 matrix
	\begin{equation}
		\begin{aligned}
			D_x=
			\begin{pmatrix}
				-1&&\\
				&1&\\
				& &1
			\end{pmatrix}
			D_y=
			\begin{pmatrix}
				1&&\\
				&-1&\\
				& &1
			\end{pmatrix}
			D_z=
			\begin{pmatrix}
				1&&\\
				&1&\\
				& &-1
			\end{pmatrix}\\
			S_{xy}=
			\begin{pmatrix}
				&1&\\
				1&&\\
				& &1
			\end{pmatrix}\quad
			S_{yz}=
			\begin{pmatrix}
				1&&\\
				&&1\\
				& 1&
			\end{pmatrix}
			S_{zx}=
			\begin{pmatrix}
				&&1\\
				&1&\\
				1& &
			\end{pmatrix}
		\end{aligned}.
	\end{equation}
	As all generators are symmetric matrices, the transpose of any element in $G$ is still in $G$.
	
	Now, we consider $\p_0=\pr{1,2,3}^\top$ and let $G$ act on it to get an orbit set $\br{\p_k},k=1,\cdots, 48$, where $\p_k=G_k\p_0$.  Hence, we have
	\begin{equation}
		\label{eq:qsol-condition}
		\int_{\tw}\pr{\p_k\cdot\n}^2=\int_{\tw}\pr{\p_0\cdot G_k^\top\n}^2=\int_{\tw}\pr{\p_0\cdot\n}^2.
	\end{equation}
	The last equality is because $G_k^\top\in G$ and $\tw$ is cubic symmetric.
	
	Note that the integral in~\eqref{eq:lem-p3sq-13} is a quadratic form on $\p$, i.e.,
	\begin{equation}
		\int_{\tw}p_3^2=\int_{\tw}\pr{\p\cdot\n}^2=\p^\top\mathbf{Q}\p,
	\end{equation}
	where $\mathbf Q$ is a symmetric matrix defined as
	\begin{equation}
		\mathbf{Q}= \int_{\tw}\n\n^\top.
	\end{equation}
	Since $\mathbf Q$ is symmetric, there are only 6 independent elements which can be solved from the following linear system derived from~\eqref{eq:qsol-condition}:
	\begin{equation}
		\label{eq:lin-sol-q}
		\p_k^\top\mathbf Q\p_k = c\|\p_0\|^2,\quad k=1,\cdots, 48.
	\end{equation}
	Here, we regard elements of $\mathbf Q$ as the unknown variables, and
	$c=\int_{\tw}\pr{\p_0\cdot\n}^2/\|\p_0\|^2$ is a constant.
	It can be verified through computation that the matrix form of the linear system~\eqref{eq:lin-sol-q} has rank $6$, and thus it has a unique solution.
	However, $\mathbf Q=c \mathbf I$ is a solution of~\eqref{eq:lin-sol-q}.
	Hence, the only solution to~\eqref{eq:lin-sol-q} is $\mathbf Q=c\mathbf I$.
	Consequently, any $\p\in\mathbb R^3$ with $\|\p\|=1$ satisfies
	\begin{equation}
		\label{eq:all-dir-same}
		\int_{\tw} \pr{\p\cdot\n}^2 = c.
	\end{equation}
	
	Extend $\br{\p}$ into a orthonormal basis $\br{\p,\p^1,\p^2}$, then we have
	\begin{equation}
		\int_{\tw} 1 = \int_{\tw} \pr{\p\cdot\n}^2+\pr{\p^1\cdot\n}^2 +\pr{\p^2\cdot\n}^2
	\end{equation}
	and 
	\begin{equation}
		\int_{\tw}\pr{\p\cdot\n}^2=
		\int_{\tw}\pr{\p^1\cdot\n}^2=
		\int_{\tw}\pr{\p^2\cdot\n}^2
	\end{equation}
	according to~\eqref{eq:all-dir-same}.
	It follows that
	\begin{equation}
		\int_{\tw}\pr{\p\cdot\n}^2=\frac{1}{3}|\tw|,
	\end{equation}
	and thus proves the lemma.
\end{proof}

\begin{lemma}
	\label{lem:div-tg-hp3}
	For a periodic surface $\tw$  and constant vector $\p\in\R^3$, we have the following identity
	\begin{equation}
		\label{eq:divpw-2hp3}
		\div \p_{\tw} = 2Hp_3
	\end{equation}
	where `$\,\div$' denotes the intrinsic divergence operator on $\tw$; $\p_{\tw}$ is the tangent component of $\p$ on $\tw$ and $H$ denotes the mean curvature.
\end{lemma}

\begin{proof}
	We can verify this through direct computation. By definition of divergence, we have
	\begin{equation}
		\div\p_{\tw}=a^{\a\b}\partial_\a p_\b -a^{\a\c}\Gamma_{\a\c}^\b p_\b
	\end{equation}
	under the local coordinate system of $\tw$.
	Note that
	\begin{equation}
		\begin{aligned}
			\partial_{\a} p_\b =\partial_\a \pr{\p\cdot \boldsymbol{a}_\b}=\p\cdot\partial_\a \boldsymbol{a}_\b&=\p\cdot\pr{\Gamma_{\a\b}^\c \boldsymbol a_\c + b_{\a\b}\n}\\
			&=\Gamma_{\a\b}^\c p_\c+b_{\a\b}p_3
		\end{aligned}
	\end{equation}
	Hence
	\begin{equation}
		\div\p_{\tw}= a^{\a\b}b_{\a\b} p_3= 2H p_3
	\end{equation}
\end{proof}

\begin{lemma}
	\label{lem:vol-remainder}
	The volume of $\toe$ is given by
	\begin{equation}
		\label{eq:close-volum-form}
		|\toe|=2\eps |\tw| +\frac{4\pi}{3}\eps^3\chi(\tw)
	\end{equation}
	where $\chi(\tw)$ is the Euler characteristic of middle surface $\tw$.
\end{lemma}
\begin{proof}
	Considering the local parameterization of $\toe$
	\begin{equation}
		\bs r^\eps:U\times (\eps,\eps)\to \toe,\, (y, x_3)\mapsto \bs r(y)+\eps x_3\bs n(y)
	\end{equation}
	where $\bs r(y)$ is the parameterization map of $\tw$, and $U\subset\R^2$ is the parameter domain.
	Let $g^\eps(y,x_3)$ and $\sqrt{a}$ be the volume forms along $\toe$ and $\tw$ induced by $\bs r^\eps$ and $\bs r$, respectively. 
	Then one can verify the following relation (see also Lemma~\ref{lem:vector-asymp-approx})
	\begin{equation}
		\sqrt{g^\eps(y,x_3)}=\pr{1-2H(y) x_3 + K(y)x_3^2}\sqrt{a(y)}
	\end{equation}
	Therefore, we have
	\begin{equation}
		\begin{aligned}
			|\toe|&=\int_U\int_{-\eps}^\eps\sqrt{g^\eps}dx_3 dy\\
			&=\int_U\int_{-\eps}^\eps\pr{1-2H(y) x_3 + K(y)x_3^2}\sqrt{a(y)}dx_3 dy\\
			&=\int_U \pr{2\eps +\frac{2}{3}\eps^3K(y)}\sqrt{a} dy\\
			&=2\eps|\tw|+\frac{2}{3}\eps^3\int_{\tw}K 
		\end{aligned}
	\end{equation}
	Applying Gauss-Bonnet theorem, we have 
	\begin{equation}
		\int_{\tw} K = 2\pi \chi(\tw),
	\end{equation}
	then we arrive at~\eqref{eq:close-volum-form}.
\end{proof}

\section{Effective conductivity of lattice metamaterial}
\label{sec:eff-cond}
Numerical homogenization is a well-established tool to determine the effective property of lattice metamaterials~\cite{Allaire2002}.
Given a RVE $Y=[-1,1]^3$ of the lattice metamaterial, with the solid region denoted as $\Omega\subset Y$,
the homogenization method computes the effective conductivity tensor $\k_H:=k_H^{ij}\e_i\otimes\e_j$ as 
\begin{equation}
	\label{eq:CH}
	k_H^{ij}=\frac{1}{|Y|}\int_\Omega\k(\nabla 
	u^{i}+\e^{i})\cdot(\nabla u^{j}+\e^{j})\itd\Omega
\end{equation}
where the temperature field $u^{i}$ is the solution of the following  \emph{cell problem} when $\p=\e^i$:
\begin{equation}
	\label{eq:cell-prop}
	\left\{
	\begin{aligned}
		&\nabla\cdot\left(\k\left(\nabla u+\p\right)\right)=\mathbf 0,\quad &x\in\Omega\qquad\\
		&\k\left(\nabla u+\p\right)\cdot\n=\mathbf 0,\quad &x\in\partial\Omega\backslash\partial Y\\
		&u \quad\text{Y-periodic}\\
		&\int_{\Omega}u= 0
	\end{aligned}
	\right.
\end{equation}
Here, the matrix
$
\k=\kappa\mathbf I
$
is the thermal conductivity matrix of the isotropic base material within $\Omega$;  $\kappa$ is the conductivity and $\mathbf I$ is the identity matrix. 
The notation $\p\in\R^3$ denotes the \emph{macroscopic temperature gradient}.
The constraint $\int_\Omega u=\mathbf 0$ is used to remove the constant shift and hence guarantee the unique solution.

The above cell problem is typically converted into a variational equation, also known as the weak form of~\eqref{eq:cell-prop}.
Considering the periodic boundary conditions, we identify $Y$ with the unit 3-dimensional flat torus $\T^3:=\R^3/(2\mathbb Z)^3$ via isometric covering map~\cite{Allaire2002}
\begin{equation}
	\label{eq:R-2-T}
	\tau:\R^3\to \T^3, (x,y,z)\mapsto (\tilde{x},\tilde{y},\tilde{z})
\end{equation}
where $\tilde{x}:=x+2\mathbb Z$, same as $\tilde{y},\tilde{z}$. We denote $\tilde{\Omega}:={\tau(\overline{{\Omega}})}$, which is a periodic domain.
Through integration by parts, the cell problem is then converted into the following variational equation
\begin{problem}
	\label{prob:vari-form}
	Find $u\in V_\#(\t\Omega)$, such that
	\begin{equation}
		\label{eq:vari-form}
		\int_{\t\Omega}\k\nabla u\cdot \nabla v=-\int_{\t\Omega} \k\p\cdot\nabla v,\quad\forall v\in V_\#(\t\Omega),
	\end{equation}
	where the function space $V_\#(\t\Omega)$ is 
	\begin{equation}
		\label{eq:func-sp-shell}
		V_\#(\t\Omega):=\left\{v\in H^1(\t\Omega):  \int_{\t\Omega} v= 0 \right\}.
	\end{equation}
	and $H^1(\tO)$ denotes the first-order Sobolev space.
\end{problem}
\begin{remark}
	The cell problem is also equivalent to the following minimization problem
	\begin{equation}
		\label{eq:vari-min}
		\argmin_{u\in V_\#(\t\Omega)} I_{\t\Omega}(u;\p)
	\end{equation}
	where 
	\begin{equation*}
		I_{\t\Omega}(u;\p):=\frac{1}{|Y|}\int_{\t\Omega} \k(\nabla u+\p)\cdot(\nabla u+\p).
	\end{equation*}
\end{remark}
\begin{remark}
	\label{rm:homogenization}
	The effective directional conductivity is computed as
	\begin{equation}
		\ka_{\p}=\k_H\p\cdot\p=k^{ij}_H p_ip_j
	\end{equation}
	for a direction $\p$ with $\|\p\|=1$.
	According to~\eqref{eq:CH} and noting that the solution $u$ of cell problem is linearly dependent on $\p$, we have
	\begin{equation}
		\ka_{\p} =\frac{1}{|Y|}\int_{\tO}\k\pr{\nabla u_{\p}+\p}\cdot\pr{\nabla u_{\p}+\p}=\min_{u\in V_\#(\tO)} I_{\tO}(u;\p)
	\end{equation}
	where $u_{\p}$ is the solution of cell problem when the macroscopic temperature gradient is $\p$.
	In the main text, we are considering the  shell lattice $\toe$ and use notation $\ka_\eps(\tw;\p)$ to denote the effective conductivity. 
	Applying above formulas, we have
	\begin{equation}
		\label{eq:min-energy}
		\ka_\eps(\tw;\p) = \min_{u\in V_\#(\toe)} I_{\toe}(u;\p)
	\end{equation}
	
\end{remark}

\begin{remark}
	\label{eq:pois-shell}
	As the base material is isotropic with conductivity matrix $\k=\ka\mathbf I$, we have $\Delta u=0$ according to the first equation in~\eqref{eq:cell-prop}.
\end{remark}

\section{Mathematical model of shell lattice metamaterial}
\label{sec:shell-repr}

\emph{
The content of this section largely follows Section 2 of the supplementary material in~\cite{ads}, and is included here for completeness and ease of reference.
	}

\subsection{Representation}
We consider a closed smooth 2-manifold $\tw$ isometrically embedded in $\T^3$,
which serves as the middle surface of the shell. 
As $\tw$ is compact, there exist finite charts $\br{\pr{S_i,\boldsymbol\psi_i}}_{i\in I}$ covering $\tw$, where 
$
\boldsymbol{\psi}_i:= S_i\to \U_i
$
is a map from a open set $S_i\subset\tw$ to the parameter domain $\U_i\subset\R^2$ and
$I$ is a finite index set.
We suppose all $\U_i$ are bounded in $\R^2$ and $\boldsymbol{\psi}_i$ are  continuous to boundary of $\U_i$, w.l.o.g.
Then we construct the shell as a 3-manifold embedded in $\T^3$ by offsetting  $\tw$ normally on both sides by a small distance.
This can be formally stated as the following lemma
\begin{lemma}
	\def\bx{\boldsymbol x}
	\def\by{\boldsymbol y}
	\label{lem:shell-is-mani}
	For sufficiently small\footnote{ Hereafter, sufficiently small $\eps$ means there exist $\eps_0$, such that the conclusion is true when $0<\eps\le\eps_0$. }
	$\eps$, the shell lattice  is the subset 
	\begin{equation}
		\label{eq:toe-uni-thick}
		\tilde{\Omega}^\eps:=\left\{{\bx}\in\T^3: {\bx}={\by}+x_3\n({\by}), x_3\in(-\eps,\eps),{\by}\in\tw\right\},
	\end{equation}
	where $\n({\boldsymbol{y}})$ denotes the normal vector of $\tw$ at ${\boldsymbol{y}}\in\tw$. 
	Moreover, it is a smooth manifold embedded in $\T^3$ with charts $(B_i^\eps,\boldsymbol{\Phi}_i^\eps)_{i\in I}$ given by
	\begin{equation}
		\label{eq:shell-chart-eps-map}
		\begin{aligned}
			B_i^\eps=\left\{{\bx}\in\T^3: {\bx}={\by}+x_3\n({\by}), x_3\in(-\eps,\eps),{\by}\in S_i\right\}\\
			\P_i^\eps: B_i^\eps\to \mathcal{V}_i^\eps,\, {\by}+x_3\n({\by})\mapsto\left(\boldsymbol\psi_i(\by),x_3\right)
		\end{aligned}
	\end{equation}
	where $\mathcal{V}_i^\eps:=\U_i\times(-\eps,\eps)$.
\end{lemma}
\noindent As  above lemma is clear, the proof is omitted.

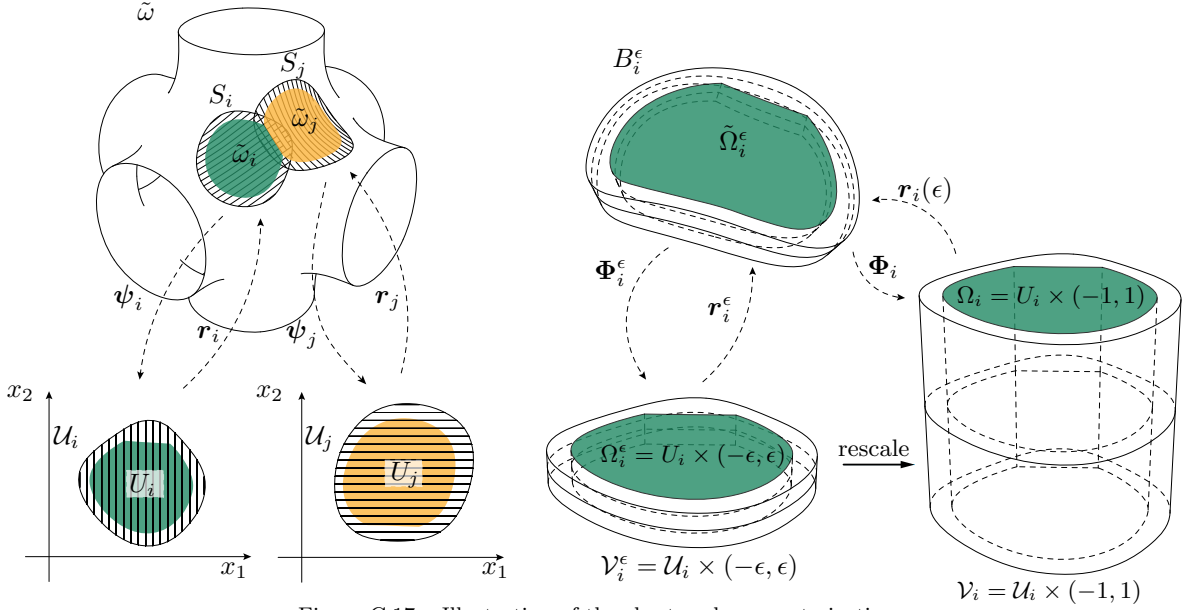
\begin{figure*}[h]
	\centering
	\begin{overpic}[width=0.95\linewidth]{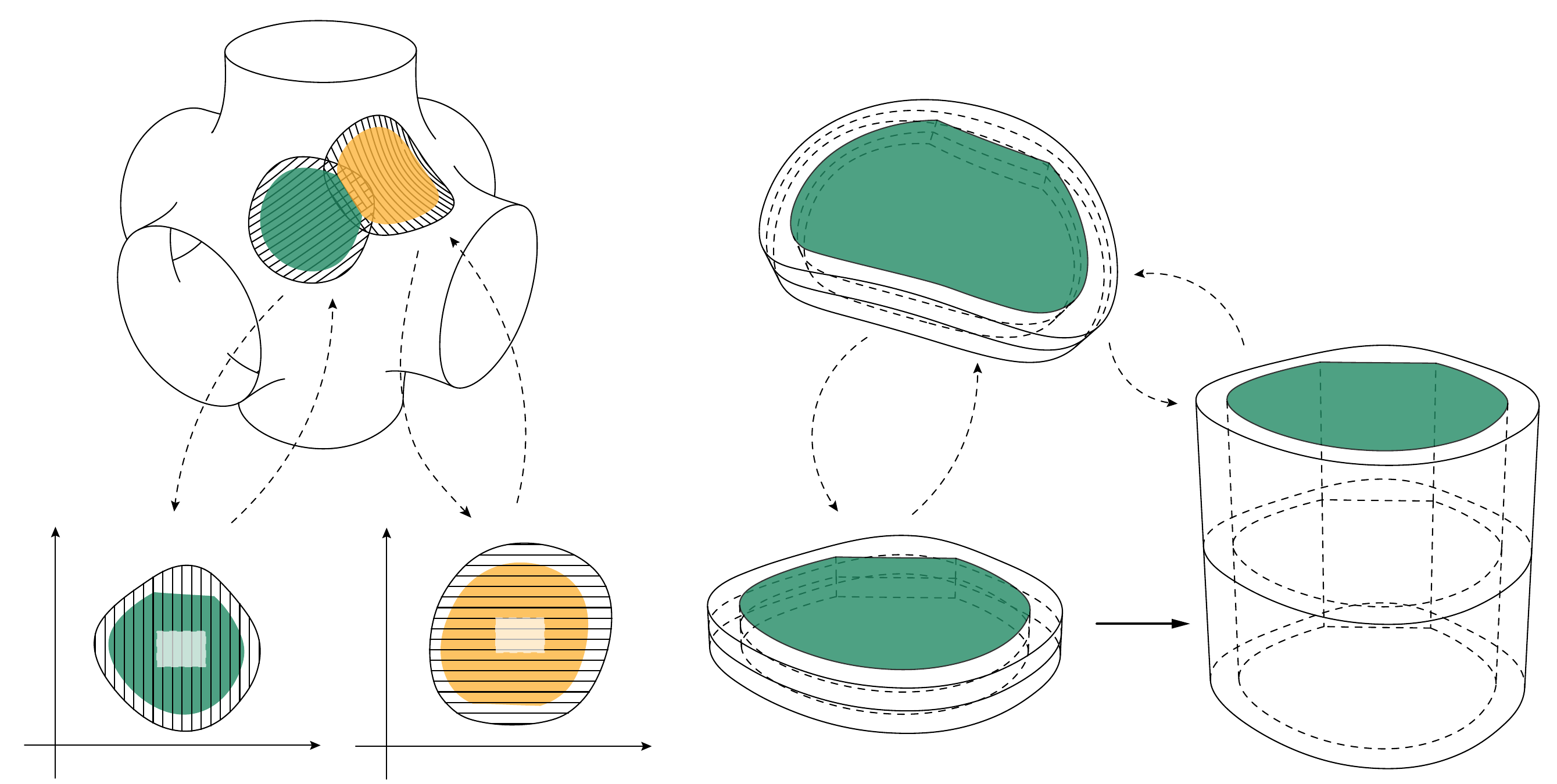}
		{
			\put(11,48){$\tw$}
			\put(19,36){$\tw_i$}
			\put(24,39){$\tw_j$}
			\put(17,41){$S_i$}
			\put(23,43.7){$S_j$}
			\put(10.3,8){$U_i$}
			\put(4,12){$\U_i$}
			\put(32.3,9){$U_j$}
			\put(25.2,12){$\U_j$}
			\put(18,1){$x_1$}
			\put(40,1){$x_1$}
			\put(0,16){$x_2$}
			\put(21,16){$x_2$}
			\put(9,24){$\bs\psi_i$}
			\put(16,21){$\bs r_i$}
			\put(23.5,21){$\bs\psi_j$}
			\put(31,24){$\bs r_j$}
			\put(50,1){$\V_i^\eps=\U_i\times(-\eps,\eps)$}
			\put(80,-1){\small$\V_i=\U_i\times(-1,1)$}
			\put(50,10.5){\small$\Omega_i^\eps=U_i\times(-\eps,\eps)$}
			\put(80,24){\small$\Omega_i=U_i\times(-1,1)$}
			\put(51,44){$B_i^\eps$}
			\put(60,37){$\tO_i^\eps$}
			\put(49.5,26){$\P_i^\eps$}
			\put(59,23){$\bs r_i^\eps$}
			\put(72.5,26.5){$\P_i$}
			\put(75,33){$\bs r_i(\eps)$}
			\put(70,11){\small rescale}
		}
	\end{overpic}
	\vspace{-2mm}
	\caption{
		Illustration of  the chart and parameterization.
	}
	\label{fig:parameter}
\end{figure*}

To conduct the asymptotic analysis as $\epsilon$ approaches zero, we utilize the native curvilinear coordinates on shell.
These curvilinear coordinates rely on the chart of the middle surface.
However, covering the middle surface $\tilde\omega$ with a single chart is infeasible due to its  non-trivial topology. Therefore, we partition the middle surface $\tilde{\omega}$ into several disjoint patches $\{\tw_i\}$.
We suppose $\tw_i$ is covered by chart $(S_i,\bs\psi_i)$. 
The parametric domain for $\tw_i$ is denoted as $U_i:=\boldsymbol\psi_i(\tw_i)$.
It follows that $\tw_i\subset S_i$ and $U_i\subset \U_i$.
Subsequently, we split the shell along normal direction from the cut locus on the mid-surface. 
The split patches constitutes the partition of $\toe$, denoted as $\{\t{\Omega}_i^\eps\}$, and the corresponding chart of each $\toe_i$ is $(B_i^\eps,\boldsymbol\Phi^\eps_i)$.
We denote  $\Omega_i^\eps:=\boldsymbol\Phi^\eps_i(\t{\Omega}_i^\eps)$.
Then we have $\t{\Omega}_i^\eps\subset B_i^\eps$ and $\Omega_i^\eps=\P^\eps_i(\t{\Omega}_i^\eps)=U_i\times(-\eps,\eps)\subset\mathcal{V}_i^\eps$.
See the illustration in Figure~\ref{fig:parameter}.

\subsection{Symbol definition}

For each chart of the shell, we define the related geometry symbols, which are mostly borrowed from \cite{ads,CiarMem1996,CiarletFlex1996}. 
For brevity, we omit the index of the chart in certain cases, and the symbols are defined and referenced on the chart indicated by the context. 
We adopt the following conventions throughout this article:
\begin{itemize}
	\item  {Greek indices (except subscript of chart) and exponents (excluding $\epsilon$) are drawn from the set $\{1, 2\}$.}
	\item {Latin indices and exponents (except when explicitly specified, such as when used for indexing sequences) are selected from the set $\{1, 2, 3\}$.}
\end{itemize}

For each chart $(B_m^\eps,\P_m^\eps)$,  its parameterization $\bs{r}^\eps_m:=(\P_m^\eps)^{-1}$ is 
\begin{equation}
	\label{eq:shell-solid-parameter}
	\bs{r}^\eps_m: \Omega_m^\eps\to \t{\Omega}_m^\eps,\, (x_1,x_2,x_3)\mapsto \r_m(x_1,x_2)+x_3\n(x_1,x_2)
\end{equation}
where $\Omega^\eps_m=U_m\times(-\eps,\eps)$ and  $\r_m:=\boldsymbol{\psi}_m^{-1}$.
We denote the tangent vector of the surface as
\begin{equation}
	\label{eq:surf-tangent}
	\boldsymbol{a}_\alpha=\partial_{\alpha}\r_m,\quad\alpha=1,2
\end{equation}
which are linearly independent on $U_m$.
They form the covariant basis of the tangent plane to the surface $ \tw$.
The two vectors $\ba^\alpha,\,\alpha=1,2$  of the same tangent plane defined by the relation
\begin{equation}
	\label{eq:contr-tangent}
	\ba^\alpha\cdot\ba_\beta=\delta^\alpha_\beta
\end{equation}
constitute its contravariant basis.
We also define the unit vector 
\begin{equation}
	\label{eq:unit-3}
	\boldsymbol{a}_3(y)=\boldsymbol{a}^3(y):=\boldsymbol n(y)
\end{equation}
which is normal to $\tw$ at the point $y\in\tw$.

We define the \emph{first fundamental form}, also known as the metric tensor, $(a_{\a\b})$ or $(a^{\a\b})$ (in covariant or contravariant components); the \emph{second fundamental form}, also known as the \emph{curvature tensor}, $(b_{\a\b})$ or $(b^{\a\b})$ (in covariant or contravariant components), and the \emph{Christoffel symbols} $\Gamma_{\alpha \beta}^\sigma$, of the surface $\t\omega$ by letting:
\begin{equation}
	\label{eq:tensors}
	\begin{gathered}
		a_{\alpha \beta}:=\boldsymbol{a}_\alpha \cdot \boldsymbol{a}_\beta, \quad a^{\alpha \beta}:=\boldsymbol{a}^\alpha \cdot \boldsymbol{a}^\beta, \\
		b_{\alpha \beta}:=\boldsymbol{a}^3 \cdot \partial_\beta \boldsymbol{a}_\alpha, \quad b_\alpha^\beta:=a^{\beta \sigma} b_{\sigma \alpha}, \\
		\Gamma_{\alpha \beta}^\sigma:=\boldsymbol{a}^\sigma \cdot \partial_\beta \boldsymbol{a}_\alpha .
	\end{gathered}
\end{equation}
Note the symmetries:
\begin{equation}
	\label{eq:sym-symmetry}
	a_{\alpha \beta}=a_{\beta \alpha}, \quad a^{\alpha \beta}=a^{\beta \alpha}, \quad b_{\alpha \beta}=b_{\beta \alpha}, \quad \Gamma_{\alpha \beta}^\sigma=\Gamma_{\beta \alpha}^\sigma.
\end{equation}
The \emph{volume form} along $\t\omega$ is  
\begin{equation}
	\label{eq:area-element}
	\itd\t{\omega} =\sqrt{a}\itd x^1\wedge \itd x^2, 
\end{equation}
where 
\begin{equation}
	\label{eq:area-det}
	a:=\det{(a_{\a\b})}.
\end{equation}

Similarly, the tangent vectors over $\tO_m^\eps$,
\begin{equation}
	\g^\epsilon_i:=\partial_i\bs{r}_m^\epsilon
\end{equation}
are linearly independent, which form the covariant basis of the tangent space locally, and the three vectors $\boldsymbol{g}^{i, \epsilon}$ defined by
\begin{equation}
	\g^{j,\epsilon}\cdot\g_{i}^\epsilon = \delta^j_i
\end{equation}
form the contravariant basis. 
We then define the \emph{metric tensor} $\left(g_{i j}^{\varepsilon}\right)$ or $\left(g^{i j, \varepsilon}\right)$ (in covariant or contravariant components) and the \emph{Christoffel symbols} of $\toe$ by letting 
\begin{equation}
	\label{eq:vol-tensor}
	\begin{gathered}
		g_{i j}^{\epsilon}:=\boldsymbol{g}_i^{\epsilon} \cdot \boldsymbol{g}_j^{\epsilon}, \quad g^{i j, \epsilon}:=\boldsymbol{g}^{i, \epsilon} \cdot \boldsymbol{g}^{j, \epsilon}, \\
		\Gamma_{i j}^{p, \epsilon}:=\boldsymbol{g}^{p, \epsilon} \cdot \partial_i \boldsymbol{g}_j^{\epsilon} .
	\end{gathered}
\end{equation}
Note the symmetries
\begin{equation}
	\label{eq:define-eps-g-Gamma}
	g_{i j}^{\epsilon}=g_{j i}^{\epsilon}, \quad g^{i j, \epsilon}=g^{j i, \epsilon}, \quad \Gamma_{i j}^{p, \epsilon}=\Gamma_{j i}^{p, \epsilon} .
\end{equation}
The \emph{volume form} on $\toe$ is 
\begin{equation}
	\label{eq:vol-element}
	\itd\toe:=\sqrt{g^{\epsilon}} \itd x^1\wedge \itd x^2\wedge \itd x^3,
\end{equation}
where
\begin{equation}
	\label{eq:vol-det-eps}
	g^{\epsilon}:=\operatorname{det}\left(g_{i j}^{\epsilon}\right).
\end{equation}
For a vector $\p\in\R^3$, we denote its covariant components or contravariant components as
\begin{equation}
	\begin{aligned}
		p_i^\eps = \p\cdot\g_i^\eps,\quad
		p^{i,\eps} = \p\cdot\g^{i,\eps}
	\end{aligned}
\end{equation}
under the local coordinate system of $\toe$ and
\begin{equation}
	p_i =\p\cdot \ba_i,\quad p^i = \p\cdot\ba^i 
\end{equation}
under the local coordinate system of $\tw$.
The conductivity matrix of base material is written as
\begin{equation}
	\k=\ka g^{ij,\eps} \g_i^\eps\otimes\g_j^\eps
\end{equation}
under the local coordinate system.

The gradient of a scalar field $v$ over $\toe$ is denoted by
\begin{equation}
	\label{eq:grad-bef}
	\nabla v=e_i^\eps(v) \g^{i,\eps}
\end{equation}
where $e_i^\eps(v)=\partial_i v$. This seemingly redundant notation is introduced to facilitate the corresponding definitions on the rescaled domain, as we discuss in next section.

\begin{remark}
	\label{rem:vol-form-bounded}
	As we have assumed that each chart map $\bs\psi_m$ is smoothly continuous to the boundary of $\U_m$,  the above functions defined in this section (except $e^\eps_i(v)$) are also continuous to the boundary of $\U_m$.
	Hence, they are all bounded on each chart.
\end{remark}

\subsection{Rescaled parameter domain}
\label{sec:rescale-prob}

To facilitate the asymptotic analysis, we eliminate the dependence of the parameter domain on $\epsilon$ by rescaling it along the surface normal direction.
This is achieved by introducing a new chart $(B_i^\eps, \P_i)$ defined as
\begin{equation}
	\P_i: B_i^\eps\to \mathcal{V}_i,\, {\bs{y}}+x_3\n({\bs{y}})\mapsto\left(\boldsymbol\psi_i({\bs{y}}),\frac{x_3}{\eps}\right),
\end{equation}
where the rescaled domain is given by
\begin{equation}
	\V_i:=\U_i\times(-1,1).
\end{equation}
After rescaling, the domain $\Omega_m^\eps=U_m\times(-\eps,\eps)$ becomes $\Omega_m:=U_m\times(-1,1)$ accordingly.
Then the parameterization on the rescaled domain, denoted as $\bs{r}_i(\eps)=\P_m^{-1}$, becomes
\begin{equation}
	\bs{r}_i(\eps) : \V_i\to B_i^\eps,\, (y,x_3)\mapsto \bs{r}^\eps_i(y,\eps x_3).
\end{equation}
The geometric interpretation of these mappings is illustrated in Figure~\ref{fig:parameter}.

We associate any scalar function $f^\eps$ defined on the original domain with a function $f(\eps)$ defined on the rescaled domain via
\begin{equation}
	\label{eq:rescale-define-scalar}
	f(\eps)(x_1,x_2,x_3):=f^\eps(x_1,x_2,\eps x_3).
\end{equation}
This convention is consistently used for all geometric and physical quantities under rescaling, such as $\Gamma^p_{ij}(\eps),\,g^{ij}(\eps)$, $\g^i(\eps)$, the vector component $p_i(\eps)$ or $p^i(\eps)$ for $\p\in\R^3$, etc.

We also introduce the following symbols on the rescaled domain to denote the components of the gradient $\nabla v$
\begin{equation}
	\label{eq:def-derivative-symb}
	\begin{aligned}
		e_{\a}(\eps)(v)=\partial_\a v(\eps),\quad
		& e_{3}(\eps)(v)=\frac{1}{\eps}\partial_3 v(\eps).
	\end{aligned}
\end{equation}
Then we have $ \nabla v = e_i(\eps)(v) \g^i(\eps) $ on $\Omega_m$.
For comparison, the unscaled definition is given in~\eqref{eq:grad-bef}.

\section{Asymptotic analysis}
\label{app:asym-ana}
We first introduce several preliminary results that are will be used in the subsequent analysis.
We denote the Sobolev norm~\cite{RiemSob} of a function $f$ over $\toe$ as 
\begin{equation}
	\|f\|_{1,\toe}=\left(\int_{\toe}|f|^2 +|\nabla f|^2 \itd\toe\right)^{1/2},
\end{equation}
where $|\nabla f|$ is the length of the gradient of $f$, which can be computed on local coordinates as 
\begin{equation}
	|\nabla f|^2 = g^{ij}\partial_i f\partial_j f.
\end{equation}
Here, $g^{ij}$ is the metric tensor on local chart.
Similar definition holds for $\|\cdot\|_{1,\tw}$, which represents $H^1(\tw)$ norm.
The Sobolev norm of a function $h$ defined on a parameter domain $\Omega_m\subset\R^3$ is denoted as
\begin{equation}
	\|h\|_{0,\Omega_m} =\pr{\int_{\Omega_m}h^2}^{1/2},\quad\|h\|_{1,\Omega_m}=\pr{\int_{\Omega_m} h^2+|\nabla h|^2}^{1/2},
\end{equation}
which are the $L^2(\Omega_m)$ and $H^1(\Omega_m)$ norm, respectively.
The analogical definition also holds for $\|\cdot\|_{0,U_m}$ and $\|\cdot\|_{1,U_m}$.
When consider the convergence of functions in Sobolev space, we adopt the following convention
\begin{itemize}
	\item {Symbol `$\lesssim$' signifies `less than or equal up to a constant'.}
	\item {Symbols `$\wc$' and `$\to$' denote weak and strong convergence, respectively}.
	\item {Symbols $\mathscr{C}^\infty(\Omega_m)$ and $\mathscr{D}(\Omega_m)$ represent the space of smooth functions and smooth functions with compact support, respectively.}
\end{itemize}
\subsection{Preliminary}

We begin by deriving several fundamental results concerning the geometric quantities  defined in  Section~\ref{sec:shell-repr}.
Based on their definitions, we can immediately get the following relations
\begin{lemma}
	\label{lem:vector-asymp-approx}
	On each chart $\pr{B_m^\eps,\P_m}$, define the symbols on the rescaled domain according to~\eqref{eq:rescale-define-scalar}, then we have
	\begin{align}
		\label{eq:tgvec-low-alpha}
		\g_\a(\eps) = \ba_\a - \eps x_3 b_\a^\b \ba_\b\\
		\label{eq:tgvec-upp-alpha}
		\g^\a(\eps) =\ba^\a +\eps x_3 b_\s^\a \ba^\s + R^\a_\c(\eps)\ba^\c\\
		\label{eq:tgvec-x3}
		\boldsymbol g^3(\eps)=\boldsymbol g_3(\eps) =\boldsymbol n\\
		\label{eq:g-low-a-b}
		g_{\a\b}(\eps) = a_{\a\b}-2\eps x_3 b_{\a\b}+\eps^2 x_3^2 b_{\a}^\c b_{\c\b}\\
		\label{eq:g-upp-a-b}
		g^{\a\b}(\eps) = a^{\a\b}+2\eps x_3 b^{\a\b}+ O(\eps^2)\\
		\label{eq:g-i3}
		g^{i3}(\eps)=g_{i3}(\eps)=0\\
		\label{eq:g-33}
		g^{33}(\eps) = g_{33}(\eps)=1\\
		\label{eq:p-eps-esti}
		p_\a(\eps) = p_\a -\eps x_3 b_\a^\b p_\b,\quad
		p_3(\eps) = p_3\\
		\label{eq:g-eps-det}
		\sqrt{g(\eps)} = 1-2H\eps x_3 + K\eps^2 x_3^2\\
		\label{eq:g-det-bound}
		0<a_0\le a, \quad 0<a_0\le g(\eps)\\
		\label{eq:g-spd}
		a^{\a\b}t_\a t_\b\ge C \sum_\a |t_\a|^2,\quad g^{ij}(\eps)t_i t_j\ge C \sum_i |t_i|^2
	\end{align}
	for sufficiently small $\eps$ and some constant $a_0,\,C > 0$ independent of $\eps$. 
	Here, $H$ and $K$ denote the mean curvature and Gauss curvature, respectively; $R^\a_\c(\eps)$ is the solution of 
	\begin{equation}
		\label{eq:g-a-error-sol}
		(\delta_{\a}^\b-\eps x_3 b^\b_\a) R^\c_\b(\eps)=\eps^2 x_3^2 b_\a^\s b_\s^\c ,
	\end{equation}
	satisfying $R^\a_\c(\eps)\sim \eps^2$.
\end{lemma}
These relations are derived through direct computation, 
we refer the reader to the supplementary material of~\cite{ads} for details.

In the following discussion, we frequently discusses the normal average of a function.  
We enumerate its preliminary properties as the following lemma, which has been established in the supplementary material of~\cite{ads}.
\begin{lemma}
	\label{lem:x3-average}
	For a scalar function $f\in L^2(\Omega_m)$, we define the normal average 
	$
	\bar{f}:=\frac{1}{2}\int_{-1}^1 f \itd x_3.
	$
	
	(a) The average $\bar f$ is finite for almost all $y\in U_m$ and 
	$
	\|\bar f\|_{0,U_m}\le\frac{1}{\sqrt{2}}\|f\|_{0,\Omega_m}.
	$
	
	(b) If $f\in H^1(\Omega_m)$, then $\bar f\in H^1(U_m)$, $\partial_\a\bar f=\overline{\partial_\a f}$ and 
	$
	\|\bar f\|_{1,U_m}\le \frac{1}{\sqrt{2}}\|f\|_{1,\Omega_m}.
	$
	
	(c) If $f\in H^1(\Omega_m)$, then 
	$
	\overline{x_3 f}=\frac{1}{2}\overline{(1-x_3^2)\partial_3 f},\, \text{ a.e. } U_m.
	$
	
	(d) If $f\in H^1(\toe)$, then the normal average on each $\V_m$ constitutes a global function $\bar f\in H^1(\tw)$ and 
	$\|\bar{f}\|_{H^1(\tw)}\le \frac{C}{\sqrt{\eps}}\|f\|_{H^1(\toe)}$ for some constant $C>0$.
\end{lemma}

The final lemma is used to establish a  bound on the solution, leading to the conclusion of convergence.
\begin{lemma}
	\label{lm:shell-korn-chart}
	For sufficiently small $\eps$, there exists a constant $C>0$ independent of $\eps$, such that
	\begin{equation}
		\label{eq:chart-pointcare}
		\sum_{m}\|v(\eps)\|_{1,\V_m}^2\le C\sum_m\sum_{i}\left\|e_{i}(\eps)(v)\right\|^2_{0,\V_m}, \quad\forall v\in V_\#(\toe) 
	\end{equation}
\end{lemma}
\begin{proof}
	We prove this lemma by contradiction.
	Suppose~\eqref{eq:chart-pointcare} does not hold, then there exist  sequences $\eps_k$ and $v^{k}\in V_\#(\tO^{\eps_k})$, $k=1,2,\cdots$, such that 
	\begin{align}
		&\eps_k\to 0\quad as\,\,k\to \infty\\
		\label{eq:chart-poincare-cond-1}
		&\sum_m \|v_m^k(\eps_k)\|_{1,\V_m}^2=1\quad as\,\,k\to \infty\\
		\label{eq:chart-poincare-cond-2}
		&\|e_i(\eps_k)(v_m^{k})\|_{0,\V_m}\to 0\quad as\,\,k\to \infty
	\end{align}
	where $v^{k}_m$ denotes the coordinate representation of $v^{k}$ on chart $(B_m^\eps,\P^\eps_m)$ and the notation $v^k_m(\eps_k)$ follows from the definition in~\eqref{eq:rescale-define-scalar} (regarding $v_m^k$ as $v_m^{k,\eps}$).
	From~\eqref{eq:chart-poincare-cond-1}, we know there  exist a subsequence of $v^k(\eps_k)$, denoted as $v^n(\eps_n)$, such that
	\begin{equation}
		\begin{aligned}
			\label{eq:subsq-connv}
			v_m^n(\eps_n)\wc v_m\ \text{ in } H^1(\V_m) \quad as\,\,n\to \infty\\
			v_m^n(\eps_n)\to v_m\ \text{ in } L^2(\V_m)\quad as\,\,n\to \infty\\
		\end{aligned}
	\end{equation}
	for some $v_m\in H^1(\V_m)$ on each $\V_m$, where the second convergence is derived from the first one by Rellich–Kondrachov theorem.
	
	Next, we apply Poincar\'e inequality on each $\V_m$ and obtain
	\begin{equation}
		\label{eq:chart-ineq-appl-poinc}
		\begin{aligned}
			\|v_m^n(\eps_n)-\BB{v_m^n(\eps_n)}\|_{0,\V_m}
			&\le C\br{\sum_i\|\partial_i v^n_m(\eps_n)\|_{0,\V_m}^2}^{1/2}\\
			&= C\br{\sum_\a\|e_\a(\eps_n)(v_m^n)\|^2+\|\eps_n e_3(\eps_n)(v_m^n)\|_{0,\V_m}^2}^{1/2}\\ &\to 0
		\end{aligned}
	\end{equation}
	where
	\begin{equation}
		\BB{v_m^n(\eps_n)}:=\frac{1}{|\V_m|}\int_{\V_m} v_m^n(\eps_n)
	\end{equation}
	denotes the average over $\V_m$.
	The last term in~\eqref{eq:chart-ineq-appl-poinc} vanishes as $n\to \infty$ because of~\eqref{eq:chart-poincare-cond-2}.
	Since $v_m^n(\eps_n)$ is bounded in $L^2(\V_m)$ due to~\eqref{eq:subsq-connv}, $\BB{v_m^n(\eps_n)}$  is also uniformly bounded with respect to $n$.
	By a diagonalization argument, there exist a subsequence of $v^n(\eps_n)$, still denoted as $v^n(\eps_n)$ for brevity, such that
	\begin{equation}
		\label{eq:constant-conv}
		\BB{v_m^n(\eps_n)}\to v_m^0 \quad as\,\,n\to\infty
	\end{equation}
	for some constant $v_m^0\in \R$.
	It then follows from~\eqref{eq:chart-ineq-appl-poinc} that
	\begin{equation}
		\label{eq:v-conv}
		v^n_m(\eps_n)\to v_m^0\text{ in } L^2(\V_m)\quad as\,\, n\to\infty
	\end{equation}
	
	For any incident chart $B_{m_1},B_{m_2}$ with $B_{m_1}\cap B_{m_2}\ne\emptyset$, the coordinate transition function $\phi_{m_1m_2}:=\P_{m_2}\circ\P_{m_1}^{-1}$ is a continuous function from $\P_{m_1}(B_{m_1}\cap B_{m_2})\subset\V_{m_1}$ to $\P_{m_2}(B_{m_1}\cap B_{m_2})\subset\V_{m_2}$.
	It can be verified that $\phi_{m_1m_2}$ is independent of $\eps$ and satisfies
	\begin{equation}
		v^n_{m_2}(\eps_n)\circ\phi_{m_1m_2} = v^n_{m_1}(\eps_n),\quad \aeon\, \P_{m_1}(B_{m_1}\cap B_{m_2})
	\end{equation}
	Let $n\to\infty$, then we have
	\begin{equation}
		v^0_{m_1}\circ\phi_{m_1m_2}=v_{m_2}^0\quad\aeon \P_{m_1}(B_{m_1}\cap B_{m_2})
	\end{equation}
	which implies $v^0_{m_1} = v^0_{m_2}$.
	Therefore, $v^0_m$ equals to the same value for all $m$, denoted as $v^0$.
	
	From~\eqref{eq:subsq-connv} and~\eqref{eq:v-conv}, we know $v_m=v^0$ for all $m$.
	Then we obtain
	\begin{equation}
		\begin{aligned}
			\lim_{n\to\infty}0=\lim_{n\to\infty}\frac{1}{\eps_n}\int_{\tO^{\eps_n}}v^n=\lim_{n\to\infty}\sum_{m}\int_{\Omega_m} v^n_m(\eps_n)\sqrt{g(\eps_n)}\itd x\\ = \sum_m\int_{\Omega_m}v^0\sqrt{a}\itd x=2v^0|\omega|.
		\end{aligned}
	\end{equation}
	Consequently
	\begin{equation}
		\label{eq:zero-conv}
		v_m=v^0=0\quad \aeon \V_m
	\end{equation}
	Thus, we have 
	\begin{equation}
		\begin{aligned}
			\sum_{m}\|v_m^n(\eps_n)\|_{1,\V_m}^2&=\sum_{m}\left(\|v^n_m(\eps_n)\|_{0,\V_m}^2+\sum_i\|\partial_i v^n_m(\eps_n)\|_{0,\V_m}^2\right)\\
			&\le\sum_m\| v^n_m(\eps_n)\|_{0,\V_m}^2+\sum_m\left(\sum_\a\|e_\a(\eps_n)(v^n_m) \|_{0,\V_m}^2+\|\eps_n e_3(\eps_n)(v^n_m)\|_{0,\V_m}^2\right)\\ 
			&\to 0 \quad as\ n\ \to \infty, 
		\end{aligned}
	\end{equation}
	due to~\eqref{eq:subsq-connv},~\eqref{eq:chart-poincare-cond-2} and~\eqref{eq:zero-conv}.
	However, this contradicts~\eqref{eq:chart-poincare-cond-1}.
\end{proof}
\subsection{Proof of convergence theorem} 
\label{sec:deri-limit-equation}
This section elaborates the detailed proof of the convergence theorem in the main text.
Several preliminary lemmas are proved in last section.
We first introduce some notations for brevity.
Since $|\toe|\sim 2\eps |\tw|$ as $\eps\to0$, we let
\begin{equation}
	J_\eps(\tw;u,\p) := \frac{1}{2\eps|\tw|}\int_{\toe}\k\pr{\nb u+\p}\cdot\pr{\nb u+\p}
\end{equation}
and 
\begin{equation}
	\label{eq:j-eps-star}
	J^*_\eps(\tw;\p):=\min_{u\in V_\#(\toe)} J_\eps(\tw;u,\p).
\end{equation}
Therefore, we have
\begin{equation}
	\ka_A(\tw;\p) = \lim_{\eps\to0}J_\eps^*(\tw;\p).
\end{equation}
We let $u^\eps$ be the solution of~\eqref{eq:j-eps-star}:
\begin{equation}
	u^\eps:=\argmin_{u\in V_\#(\toe)} J_\eps(\tw;u,\p),
\end{equation}
thus $J_\eps^*(\tw;\p)=J_\eps(\tw;u^\eps,\p)$.
With above definitions, our convergence theorem is stated as
\begin{theorem}
	\label{thm:adc-theorem}
	The limit of $J_\eps^*(\tw;\p)$ satisfies
	\begin{equation}
		\label{eq:j-eps-limit-expr}
		\lim_{\eps\to0} J_\eps^*(\tw;\p) = \frac{1}{|\tw|}\int_{\tw}\ka\pr{\nabla\bar u+\p_{\tw}}\cdot\pr{\nabla\bar u+\p_{\tw}},
	\end{equation}
	where $\bar u$ is the solution of the following equation on $\tw$
	\begin{equation}
		\label{eq:avg-u-poisson}
		\Delta \bar{u} = -\div \p_{\tw},
	\end{equation}
	and $\p_{\tw}$ is the tangent component of $\p$ on $\tw$.
\end{theorem}

We divide our proof into 5 steps (\textbf{Step i} to \textbf{Step v}).
The relevant geometry symbols are defined in Section~\ref{sec:shell-repr}.

\emph{
	For brevity, we denote $e_i(\eps)(u^\eps)$ as $e_i(\eps)$ (see definition in \eqref{eq:def-derivative-symb}), and simplify the notations of integral on multiple charts `$\sum_m\int_{\Omega_m}$' and `$\sum_m\int_{U_m}$' to notations `$\int_{\Omega}$' and `$\int_{U}$', respectively. 
}

\paragraph*{\textbf{Step i}}
\itshape{
	On each chart $(B_m^\eps,\P_m)$, the following inequality holds 
	\begin{equation}
		\label{eq:strain-bound}
		\sum_{i}\|e_{i}(\eps)\|_{0,\V_m}\le C
	\end{equation}
	for some constant $C>0$.
	Consequently, we also have $\sum_{i}\|e_{i}(\eps)\|_{0,\Omega_m}\le C$ since $\Omega_m\subset\V_m$.
} 
\begin{proof}
	This can be seen from the following inequalities
	\begin{equation}
		\begin{aligned}
			\sum_{i}\|e_{i}(\epsilon)\|_{0,\V_m}^2
			 &\le  2\sum_{i}\|e_{i}(\epsilon)+p_{i}(\epsilon)\|_{0,\V_m}^2 +2\sum_i \|p_i(\eps)\|^2_{0,\V_m}\\
			&\le C_2\int_{\V_m} \ka g^{ij}(\eps)\left(e_{i}(\epsilon)+p_{i}(\epsilon)\right)\left(e_{j}(\epsilon)+p_{j}(\epsilon)\right) \sqrt{g(\eps)} \itd x + C_1\\
			&= \frac{C_2}{\eps}\int_{B_m^\eps} \k\cdot\left(\nabla u^\eps+\p\right)\cdot\left(\nabla u^\eps+\p\right) \itd \tO^\eps + C_1\\
			&\le 2{C_2} |\tw|J_\epsilon^*(\t\omega;\p)+C_1\\
			&\le 2{C_2 }|\tw|J_\epsilon(\t\omega; 0,\p)+C_1\\
			&\le C,
		\end{aligned}
	\end{equation}
	where $C_1,C_2>0$ are constants independent of $\eps$.
	The first inequality holds due to the triangle inequality.
	For the second inequality, we have used the positive definiteness of the thermal conductivity matrix and the positive bound on $g(\eps)$, 
	and $|p_i(\eps)|^2=|\p\cdot\g_i(\eps)|^2\le |\p|^2 |\g_i(\eps)|^2 $ which is also bounded independent of $\eps$.
	The third inequality (fourth line) comes from the definition of the $J_\eps^*$ and the positive integrand.
	The fourth inequality (fifth line) is because $J_\eps^*$ is the minimum of $J_\eps$ and thus is less than the value of $J_\eps$ when $u=0$.
	The last inequality is because  $J_\eps(\omega;0,\p)=\ka\|\p\|^2|\toe|/(2\eps |\tw|)\sim \ka\|\p\|^2=\ka$ as $\eps\to 0$ and thus bounded for sufficiently small $\eps$.
\end{proof}

\paragraph*{\textbf{Step ii}}There exist a subsequence of $\eps$, still denoted as $\eps$, such that as $\eps\to 0$,
\begin{align}
	\label{eq:global-u-weak-conv-h1}
	\overline{u^\eps}&\wc \bar u\quad\inspace H^1(\tw)\\
	\label{eq:global-u-strong-conv-l2}
	\overline{u^\eps}&\to \bar u\quad\inspace L^2(\tw)
\end{align}
for some $\bar u\in H^1(\tw)$, where ` $\overline{u^\eps}$' signifies the normal average defined in Lemma~\ref{lem:x3-average}(d),
and on each $\Omega_m$, we have
\begin{align}
	\label{eq:u-eps-wc-bar-u-H1}
	u(\eps)&\wc \bar u\inspace H^1(\Omega_m)\\
	\label{eq:u-eps-to-bar-u-L2}
	u(\eps)&\to \bar u\inspace L^2(\Omega_m)\\
	\label{eq:strain-wc}
	{e_i(\eps)} &\wc e_i \,\,\inspace L^2(\Omega_m)
\end{align}
for some $e_i\in L^2(\Omega_m)$. Moreover, $e_i$ is independent of $x_3$ and  we have
\begin{align}
	\label{eq:avg-ea-eq-d-u}
	\overline{e_\a} =\partial_\a \bar u \quad\aeon U_m \\
	\label{eq:strain-limit-e3}
	e_3 = -p_3 \quad\aeon \Omega_m
\end{align}
\begin{proof}
	First, we note that~\eqref{eq:strain-wc} is the immediate conclusion from the bound in \textbf{Step i}.
	
	For~\eqref{eq:global-u-weak-conv-h1} and~\eqref{eq:global-u-strong-conv-l2},  we consider the norm of $\overline{u^\eps}$
	\begin{equation}
		\label{eq:h1-norm-bound-by-step-i}
		\begin{aligned}
			\|\overline{u^\eps}\|_{1,\tw}^2&=\int_{\tw}\pr{\overline{u^\eps}}^2+\int_{\tw}|\nabla\overline{u^\eps}|^2\\
			&=\int_{U} \pr{\overline{u(\eps)}^2+ a^{\a\b}\partial_\a\overline{u(\eps)}\partial_\b\overline{u(\eps)}}\sqrt{a}\itd s\\
			&\lesssim \int_{U}\pr{\overline{u(\eps)}^2+\sum_\a\left|\partial_\a\overline{u(\eps)}\right|^2}\itd s \\
			&\lesssim\sum_m\|u(\eps)\|^2_{0,\Omega_m}+\sum_\a\|\partial_\a u(\eps)\|_{0,\Omega_m}^2\\
			&=\sum_m\|u(\eps)\|_{1,\Omega_m}^2 \lesssim\sum_m\sum_i\|e_i(\eps)\|_{0,\V_m}^2
		\end{aligned}
	\end{equation}
	Here, the first inequality is because the  metric tensor is bounded (see Remark~\ref{rem:vol-form-bounded} and Lemma~\ref{lem:vector-asymp-approx}).
	The second inequality comes from Lemma~\ref{lem:x3-average}(a)(b).
	The last inequality is due to Lemma~\ref{lm:shell-korn-chart}.
	From the result in \textbf{Step i}, there is a constant upper bound for the last term.
	Therefore, there must exist a subsequence of $\eps$ such that~\eqref{eq:global-u-weak-conv-h1} is true and~\eqref{eq:global-u-strong-conv-l2} follows from~\eqref{eq:global-u-weak-conv-h1} by Rellich–Kondrachov theorem.
	
	Note that the upper bound on~\eqref{eq:h1-norm-bound-by-step-i} also implies a weak limit $u^0\in H^1(\Omega_m)$ such that
	\begin{align}
		\label{eq:u0-converge-h1}
		u(\eps)\wc u^0\inspace H^1(\Omega_m)\\
		u(\eps)\to u^0\inspace L^2(\Omega_m)
	\end{align}
	Then Lemma~\ref{lem:x3-average} shows that $\overline{u(\eps)}\to \overline{u^0}$ in $L^2(\Omega_m)$.
	Hence, $\overline{u^0} = \bar u$ a.e. $\Omega_m$.
	As $e_3(\eps)$ weakly converges, it is bounded in $L^2(\Omega_m)$. 
	Therefore, $\partial_3 u(\eps)=\eps e_3(\eps)\to 0$ in $L^2(\Omega_m)$.
	Then for $\forall \phi\in\mathscr{D}(\Omega_m)$, we have
	\begin{equation}
		\begin{aligned}
			\int_{\Omega_m}\phi\pd_3 u^0=\lim_{\eps\to0} \int_{\Omega_m}\phi\pd_3\pr{ u^0- u(\eps)}=\lim_{\eps\to0}\int_{\Omega_m}\pr{  u(\eps)-u^0}\partial_3\phi=0,\\
			\int_{\Omega_m}e_\a\pd_3\phi= \lim_{\eps\to0}\int_{\Omega_m}e_\a(\eps)\partial_3\phi=-\lim_{\eps\to0}\int_{\Omega_m}u(\eps)\pd_{\a3}\phi=\lim_{\eps\to0}\int_{\Omega_m}\partial_3 u(\eps)\partial_\a\phi=0,
		\end{aligned}
	\end{equation}
	where the first equation holds because of~\eqref{eq:u0-converge-h1}.
	This implies that $u^0,e_\a$ are both independent of $x_3$ (Lemma 3.2 in~\cite{CiarMem1996}).
	Consequently, $\overline{u^0} = u^0$ and thus~\eqref{eq:u-eps-wc-bar-u-H1} and~\eqref{eq:u-eps-to-bar-u-L2} are proved.
	
	We can extract a diagonal subsequence such that the limits in~\eqref{eq:global-u-weak-conv-h1},~\eqref{eq:global-u-strong-conv-l2}, 
	\eqref{eq:u-eps-wc-bar-u-H1}, \eqref{eq:u-eps-to-bar-u-L2} and~\eqref{eq:strain-wc} hold simultaneously.
	
	To prove~\eqref{eq:strain-limit-e3}, we construct a field $v^\eps$ whose coodinate repsentation is
	\begin{equation}
		v(\eps)(x) = \int_{-1}^{x_3}\varphi(x) dx  \quad on\,each\,\, \Omega_n
	\end{equation}
	where $\varphi\in \mathscr{D}(\Omega_m)$ and vanishes outside $\Omega_m$.
	By construction, $v{(\epsilon)}$ is zero on all $\Omega_n$ with $n \ne m$.
	Substituting $v^\eps$ into the weak form equation~\eqref{eq:vari-form} yields
	\begin{equation}
		\begin{aligned}
			\eps\int_{\Omega_m}g^{\a\b}(\eps)\pr{e_{\a}(\eps)+p_\a(\eps)}\partial_\b v(\eps)\sqrt{g(\eps)}\itd x\\
			+ \int_{\Omega_m} (e_3(\eps) + p_3)\varphi \sqrt{g(\eps)} \itd x = 0
		\end{aligned}
	\end{equation}
	Let $\eps\to 0$, we have
	\begin{equation}
		\int_{\Omega} (e_3 + p_3)\varphi\sqrt{a}\itd x=0
	\end{equation}
	As $\varphi$ is arbitrary, above equality implies~\eqref{eq:strain-limit-e3}.

	Next, we show that~\eqref{eq:avg-ea-eq-d-u} is true.
	For $\forall\phi\in\mathscr{D}(U_m)$, we have
	\begin{equation}
		\label{eq:der-limit-eq-limit-der}
		\begin{aligned}
			&\int_{U_m}\pr{\overline{e_\a}-\partial_\a\bar u}\phi\itd s\\ &=\int_{U_m}\Big\{\pr{\overline{e_\a}-\overline{e_\a(\eps)}}+\pr{\overline{e_\a(\eps)}-\partial_\a\overline{u(\eps)}}+\pr{\partial_\a\overline{u(\eps)}-\partial_a\bar{u}}\Big\}\phi\itd s\\
			&=\int_{U_m}\pr{\overline{e_\a}-\overline{e_\a(\eps)}}\phi\itd s+\int_{U_m}\pr{\partial_\a\overline{u(\eps)}-\partial_a\bar{u}}\phi\itd s\\
			&=\frac{1}{2}\int_{\Omega_m}\pr{e_{\a}-e_\a(\eps)}\phi\itd s+\int_{U_m}\pr{\overline{u(\eps)}-\bar u}\partial_\a\phi\itd s\quad\\
		\end{aligned}
	\end{equation}
	Here, we eliminate the term $\pr{\overline{e_\a(\eps)}-\partial_\a\overline{u(\eps)}}$ in the second equality because $\overline{e_\a(\eps)}=\overline{\partial_a u(\eps)}=\partial_\a\overline{u(\eps)}$ (Lemma~\ref{lem:x3-average}(b)).
	Note that the last term in~\eqref{eq:der-limit-eq-limit-der} vanishes as $\eps\to 0$ due to~\eqref{eq:global-u-strong-conv-l2} and~\eqref{eq:strain-wc}.
	Then we have 
	\begin{equation}
		\int_{U_m}\pr{\overline{e_\a}-\partial_\a\bar u}\phi\itd s=0,\quad \forall \phi\in \mathscr{D}(U_m)
	\end{equation}
	and thus~\eqref{eq:avg-ea-eq-d-u} is proved.
\end{proof}

\paragraph*{\textbf{Step iii}}
The limit  $\bar{u}$ in \textbf{Step ii} can be solved from the following equation
\begin{equation}
	\label{eq:limit-strong-form}
	\Delta \bar{u} = -\div \p_{\tw}
\end{equation}

\begin{proof}
	We prove a equivalent result.  We show that the limit $\bar u$ satisfies the weak from of~\eqref{eq:limit-strong-form}, i.e.,  $\bar u\in V_\#(\tw)$ and satisfies
	\begin{equation}
		\label{eq:u-avr-vari-form}
		\int_{\tw} \nabla {\bar u}\cdot\nabla v =-\int_{\tw} \p_{\tw} \cdot\nabla v \quad for\, all\,\, v\in V_\#(\tw),
	\end{equation}
	where
	\begin{equation}
		V_\#(\tw):=\br{v\in H^1(\tw) : \int_{\tw} v=0}.
	\end{equation}
	
	As $u^\eps$ satisfies~\eqref{eq:vari-form}, we can see that for $\forall v\in V_\#(\tw)$
	\begin{equation}
		\int_{\toe}\nabla u^\eps\cdot\nabla v=-\int_{\toe}\p\cdot\nabla v       
	\end{equation}
	Here, we regard $v$ as a function defined in $\toe$, which is constant along normal direction.
	Hence, the above formula is still valid\footnote{Although $\int_{\toe} v$ may be nonzero, namely, $v$ may not belong to $V_\#(\toe)$, we can remove the average of $v$ over $\toe$ without changing its gradient.}.
	We  expand the integrals onto charts, resulting in
	\begin{equation}
		\int_{\Omega}g^{ij}(\eps)e_i(\eps) e_j(\eps) (v)\sqrt{g(\eps)}=-\int_{\Omega}g^{ij}(\eps)p_i(\eps)e_j(\eps)( v)\sqrt{g(\eps)}
	\end{equation}
	As $\partial_3 v(\eps)=0$ and $g^{\a3}(\eps)=0$, we have
	\begin{equation}
		\int_{\Omega} g^{\a\b}(\eps)e_\a(\eps)e_\b(\eps)(v)\sqrt{g(\eps)}=-\int_{\Omega}g^{\a\b}(\eps)p_\a(\eps)e_\b(\eps)(v)\sqrt{g(\eps)}
	\end{equation}
	Since $v$ is independent on $x_3$, $e_\b(\eps)(v)=\partial_\b v$ is also constant w.r.t. $x_3$ and independent of $\eps$.
	Therefore, it follows by letting $\eps\to0$ that
	\begin{equation}
		\int_{\Omega} a^{\a\b}e_\a\partial_\b v \sqrt{a}=-\int_{\Omega} a^{\a\b} p_\a \partial_\b v\sqrt{a}
	\end{equation}
	By Fubini theorem, we obtain
	\begin{equation}
		\int_{U}a^{\a\b}\overline{e_\a}\partial_\b v\sqrt{a}=-\int_{U}a^{\a\b}p_\a\partial_\b v\sqrt{a}
	\end{equation}
	From \textbf{Step ii}, we know  $\overline{e_\a}=\partial_\a \bar{u}$.
	Hence 
	\begin{equation}
		\begin{aligned}
			\int_{\tw}\nabla\overline{u}\cdot\nabla v&=\int_{U}a^{\a\b}\partial_\a\bar{u}\partial_\b v\sqrt{a}\\
			&=-\int_{U}a^{\a\b}p_\a\partial_\b v\sqrt{a}=-\int_{\tw}\p_{\tw}\cdot\nabla v
		\end{aligned}
	\end{equation}
	
	Next, we show $\bar u\in V_\#(\tw)$.
	Given that $\bar u\in H^1(\tw)$,  it remains to show that $\int_{\tw}\bar u=0$.
	This follows from the following relations:
	\begin{equation}
		\begin{aligned}
			0=\int_{\toe} u^\eps=\lim_{\eps\to 0}\frac{1}{\eps}\int_{\toe} u^\eps&=\lim_{\eps\to0}\int_{\Omega} u(\eps)\sqrt{g(\eps)}\\
			&=\int_{\Omega} \bar{u}\sqrt{a}=2\int_{U}\bar{u}\sqrt{a}=2\int_{\tw} \bar{u}
		\end{aligned}
	\end{equation}
	Therefore, $\bar u\in V_\#(\tw)$.
\end{proof}

\paragraph*{\textbf{Step iv}} The limit in~\eqref{eq:j-eps-limit-expr} holds for the subsequence we found in \textbf{Step ii}.
\begin{proof}
	As $u^\eps$ satisfies the variational equation~\eqref{eq:vari-form}, we have
	\begin{equation}
		\int_{\toe}\k\left(\nabla u^\eps + \p\right)\cdot\left(\nabla u^\eps + \p\right)=\int_{\toe}\k\pr{\nabla u^\eps+\p}\cdot\p.
	\end{equation}
	Thus
	\begin{equation}
		\label{eq:heat-en-lim}
		\begin{aligned}
			\lim_{\eps\to 0}J_\eps^*(\tw;\p) =\lim_{\eps\to 0}\frac{1}{2\eps|\tw|}\int_{\toe}\k\pr{\nabla u^\eps+\p}\cdot\p \\
		\end{aligned}
	\end{equation}
	We split the integral onto charts as follows
	\begin{equation}
		\label{eq:ka-int-chart}
		\begin{aligned}
			&\frac{1}{2\eps|\tw|}\int_{\toe}\k\pr{\nabla u^\eps+\p}\cdot\p\\ &=\frac{1}{2|\tw|}\bigg\{\int_{\Omega}\ka g^{\a\b}(\eps)\pr{e_\a(\eps)+p_\a(\eps)}p_\b(\eps)\sqrt{g(\eps)}\\
			&+\int_{\Omega}\ka g^{33}(\eps)\pr{e_3(\eps)+p_3(\eps)}p_3(\eps)\sqrt{g(\eps)}\bigg\}
		\end{aligned}
	\end{equation}
	where we have utilized Lemma~\ref{lem:vector-asymp-approx} to eliminate the terms involving $g^{\a3}(\eps)$.
	
	We construct a function $v^\eps$ whose coordinate representation on each chart is
	\begin{equation}
		v(\eps) = \eps x_3 p_3(\eps)\quad\text{on } \Omega_m
	\end{equation}
	Though defined separately on each chart, it is clear that this function is globally well defined and is smooth over $\toe$.
	Note that $\pd_1 v(\eps)=\pd_2 v(\eps)=0$ and $\pd_3 v(\eps)=\eps p_3(\eps)$.
	Then we substitute $v^\eps$ into the variational equation~\eqref{eq:vari-form} and obtain
	\begin{equation}
		\label{eq:p3-term-vanishes}
		\int_\Omega \ka g^{33}(\eps)\pr{e_3(\eps)+p_3(\eps)}p_3(\eps)\sqrt{g(\eps)}=0
	\end{equation}
	Therefore, the second integral in~\eqref{eq:ka-int-chart} vanishes.
	
	For the first integral in~\eqref{eq:ka-int-chart}, we utilize~\eqref{eq:strain-wc},~\eqref{eq:avg-ea-eq-d-u} and Lemma~\ref{lem:vector-asymp-approx}, which leads to the following limit
	\begin{equation}
		\begin{aligned}
			&\lim_{\eps\to0}\int_{\Omega}\ka g^{\a\b}(\eps)\pr{e_\a(\eps)+p_\a(\eps)}p_\b(\eps)\sqrt{g(\eps)}\\
			&=\int_{\Omega} \ka a^{\a\bt}\pr{e_\a +p_\a} p_{\b}\sqrt{a}\\
			&=2\int_{U}\ka a^{\a\b}\pr{\overline{ e_\a} +p_\a}p_\b \sqrt{a}\\
			&=2\int_{U}\ka a^{\a\b}\pr{\partial_\a \bar u +p_\a}p_\b \sqrt{a}\\
			&=2\int_{\tw}\ka \pr{\nabla\bar u+\p_{\tw}}\cdot\p_{\tw},
		\end{aligned}
	\end{equation}
	Therefore, we have
	\begin{equation}
		\lim_{\eps\to 0}J_\eps^*(\tw;\p) =\frac{1}{|\tw|}\int_{\tw}\ka \pr{\nabla\bar u+\p_{\tw}}\cdot\p_{\tw}.
	\end{equation}
	According to the variational equation~\eqref{eq:u-avr-vari-form}, this limit is exactly~\eqref{eq:j-eps-limit-expr}.
\end{proof}

\paragraph*{\textbf{Step v}}
The convergence in \textbf{Step vi} holds for the entire sequence of $\eps$, not only for the subsequence in \textbf{Step ii}.
\begin{proof}
	For any subsequence of $\eps$, the upper bounds in \textbf{Step i} always holds.
	Therefore, following the process we have in \textbf{Step ii} and \textbf{Step iv},  there exists a  subsequence of this subsequence that satisfies~\eqref{eq:j-eps-limit-expr}.
	However, we know the limit $\bar u$ must be the  unique solution of~\eqref{eq:limit-strong-form} according \textbf{Step iii}.
	As a result, any subsequence contains a subsequence that converges to the same limit in \textbf{Step iv}, and hence the entire sequence converges.
\end{proof}

\upshape

\section{Third-order accuracy of ADC}
\label{app:third-oder-acc}
We prove the following theorem which implies the third-order accuracy of $\ka_A$.
\begin{theorem}
	\label{thm:second-order}
	When $\eps$ is sufficiently small, we have
	\begin{equation}
		\label{eq:second-order-accuracy}
		\ka_\eps(\tw;\p) -\frac{2\eps|\tw|}{|Y|}\ka_A(\tw;\p)\sim  O(\eps^3)
	\end{equation}
\end{theorem}

\begin{proof}
	We leverage part of the results in the proof of Theorem~\ref{thm:adc-theorem}\footnote{Although several results in the proof of Theorem~\ref{thm:adc-theorem} are valid only for a subsequence, we can see that~\eqref{eq:second-order-accuracy} applies to the entire sequence following the idea of \textbf{Step v}. }.
	We construct a scalar field $u^1$ whose coordinate representation is
	\begin{equation}
		\label{eq:u-1-def}
		u^1(\eps)=\bar u -\eps x_3 p_3\quad \text{ on each } \Omega_m,
	\end{equation}
	where $\bar u$ is the solution of~\eqref{eq:avg-u-poisson}.
	Then we have
	\begin{equation}
		\label{eq:u1-gradient}
		\begin{aligned}
			e_\a(\eps)( u^1) = \partial_\a \bar u+\eps x_3 b_\a^b p_\b,\quad
			e_3(\eps)( u^1) = - p_3.
		\end{aligned}
	\end{equation}
	The field $u^1$ can be considered as an approximation of the exact solution $u^\eps$.
	
	We let 
	\[
	u_t = u^\eps+ t(u^1-u^\eps)
	\] 
	and
	\begin{equation}
		E_t = \int_{\toe}\k\pr{\nabla u_t +\p}\cdot\pr{\nabla u_t+\p}
	\end{equation}
	Clearly, $E_t$ is a quadratic function w.r.t. $t$ satisfying $E_0=|Y|\ka_\eps(\tw;\p)$.
	Through computation,
	\begin{equation}
		\begin{aligned}
			\frac{dE_t}{dt}\Big|_{t=0}&=2\int_{\toe}\k\pr{\nabla u^\eps +\p}\cdot\nabla\pr{u^1-u^\eps}=0\\
			\frac{d^2E_t}{dt^2}&=2\int_{\toe}\k\nabla \pr{u^1-u^\eps}\cdot\nabla\pr{u^1-u^\eps}\\
			&= 2\ka \int_{\toe}|\nabla(u^1-u^\eps)|^2
		\end{aligned}
	\end{equation}
	Then we have
	\begin{equation}
		E_1 = E_0 + \frac{1}{2}\frac{d^2 E_t}{d t^2}
	\end{equation}
	Therefore, it suffices to prove the following relations
	\begin{align}
		\label{eq:E-1-eps3}
		E_1&= 2\eps|\tw| \ka_A(\tw;\p) +O(\eps^3)\\
		\label{eq:sec-der-eps3}
		\frac{d^2 E_t}{d t^2}&=O(\eps^3)
	\end{align}
	
	To prove~\eqref{eq:E-1-eps3}, we notice
	\begin{equation}
		\begin{aligned}
			E_1&=\int_{\toe}\k\cdot\pr{\nabla u^1+\p}\cdot\pr{\nabla u^1+\p}\\
&=\eps\int_{\Omega} \ka g^{\a\b}(\eps)\big(e_\a(\eps) (u^1)+p_\a(\eps)\big)\pr{e_\b(\eps) (u^1)+p_\b(\eps)}\sqrt{g(\eps)}\\
&=\eps\int_{\Omega} \ka g^{\a\b}(\eps)\big(\partial_\a\bar u+p_\a\big)\pr{\partial_\b \bar u+p_\b}\sqrt{g(\eps)}\\
&=\eps\int_{\Omega}\ka\bigg\{ \pr{a^{\a\b}+2\eps x_3 b^{\a\b}}\pr{\pd_\a\bar u+p_\a}\pr{\pd_\b\bar u+p_\b}(1-2\eps x_3H)+O(\eps^2)\bigg\}\sqrt{a},
		\end{aligned}
	\end{equation}
	Here, we have utilized Lemma~\ref{lem:vector-asymp-approx} and~\eqref{eq:u1-gradient}.
	Note that the last term contains a $O(\eps^3)$ term and two principle parts:
	\begin{equation}
		\begin{aligned}
			I_1(\eps) &= \ka\eps\int_{\Omega} a^{\a\b}\pr{\pd_\a\bar u+p_\a}\pr{\pd_\b\bar u+p_\b} \sqrt{a}\\
			I_2(\eps) &= 2\eps^2\int_{\Omega}\ka x_3\pr{b^{\a\b}-Ha^{\a\b}}\pr{\pd_\a\bar u+p_\a}\pr{\pd_\b\bar u+p_\b}\sqrt{a},
		\end{aligned}
	\end{equation}
	Since the integrand of $I_2$ is an odd function of $x_3$ and  $\Omega_m$ is symmetry w.r.t. $x_3$, we know $$I_2(\eps)=0.$$
	For $I_1(\eps)$, since the integrand is independent of $x_3$, we have
	\begin{equation}
		\begin{aligned}
			I_1(\eps) &= 2\eps\int_{U}\ka a^{\a\b}\pr{\pd_\a\bar u+p_\a}\pr{\pd_\b\bar u+p_\b}\sqrt{a}\\
			&=2\eps\ka \int_{\tw}\pr{\nabla \bar u + \p_{\tw}}\cdot\pr{\nabla \bar u + \p_{\tw}}
		\end{aligned}
	\end{equation}
	According to~\eqref{eq:ka-dir-expr}, we have $I_1(\eps) = 2\eps|\tw| \ka_A(\tw;\p)$.
	Hence,~\eqref{eq:E-1-eps3} is proved.
	
	To prove~\eqref{eq:sec-der-eps3}, consider the following relation
	\begin{equation}
		\label{eq:split-d2Et-to-parts}
		\begin{aligned}
			\frac{1}{2\ka}\frac{d^2E_t}{dt^2}&=\int_{\toe}\nabla\cdot\pr{(u^1-u^\eps)\nabla\pr{u^1-u^\eps}}\\
			&\quad -\int_{\toe}\pr{u^1-u^\eps}\Delta\pr{u^1-u^\eps}
		\end{aligned}
	\end{equation}
	
	For the first integral in~\eqref{eq:split-d2Et-to-parts}, we have
	\begin{equation}
		\label{eq:first-E6}
		\begin{aligned}
			\int_{\toe}\nabla\cdot\pr{(u^1-u^\eps)\nabla\pr{u^1-u^\eps}}
			=\int_{\partial\toe}(u^1-u^\eps)\nabla_{\n}\pr{u^1-u^\eps}
		\end{aligned}
	\end{equation}
	By construction, $u^1$ and $u^\eps$ satisfy the same boundary condition $\nabla_{\n} u=-p_3$  (top surface) or $p_3$ (bottom face) on $\partial\toe$ (see~\eqref{eq:cell-prop}).
	Therefore,~\eqref{eq:first-E6} equals zero.
	
	For the second integral in~\eqref{eq:split-d2Et-to-parts}, we have $\Delta u^\eps =0$ (Remark~\ref{eq:pois-shell}), thus it remains to consider the following integral
	\begin{equation}
		\begin{aligned}
			\int_{\toe}\pr{u^1 -u^\eps} \Delta u^1=\eps\int_{\Omega}(u^1-u^\eps)\Big(\pr{\Delta_{\tw} \bar u+2Hp_{3}} + \eps c_1(x_1,x_2) x_3+O(\eps^2)\Big)\sqrt{g(\eps)}
		\end{aligned}
	\end{equation}
	where $\Delta_{\tw}$ denotes the Laplace–Beltrami operator on $\tw$, and $c_1$ is 
	a continuous functions on each ${\U_m} $;  the relation
	\begin{equation}
		\Delta u^1 = (\Delta_{\tw} \bar u+2Hp_3) + \eps c_1(x_1,x_2) x_3+O(\eps^2)
	\end{equation}
	can be verified through direct computation based on definition~\eqref{eq:u-1-def}.
	Since $\Delta_{\tw}\bar u + 2 H p_{3}=0$ according to~\eqref{eq:avg-u-poisson} and~\eqref{eq:divpw-2hp3}, we get
	\begin{equation}
		\begin{aligned}
			&\int_{\toe}\pr{u^1-u^\eps}\Delta u^1 \\
			&=\eps\int_{\Omega}\pr{u^1-u^\eps}\eps c_1(x_1,x_2)x_3\sqrt{a}\itd s + O(\eps^3)\\
			&=2\eps^2\int_{U}\overline{(u^1-u^\eps)x_3}c_1(x_1,x_2)\sqrt{a}\itd s + O(\eps^3)\\
			&=\eps^2 \int_{U}\overline{\pr{1-x_3^2}\partial_3\pr{u^1-u^\eps}} c_1(x_1,x_2)\sqrt{a}\itd s + O(\eps^3)\\
			&=\eps^3\pr{\int_{U}\overline{\pr{x_3^2-1}\pr{e_3(\eps)+p_3}}c_1(x_1,x_2)\sqrt{a}\itd s} + O(\eps^3)
		\end{aligned}
	\end{equation}
	Here, the first equality applies~\eqref{eq:g-eps-det}; the third equality  is from Lemma~\ref{lem:x3-average}(c); the fourth equality is from the definition of $u^1$ and $e_3(\eps)$.
	As $e_3(\eps)+p_3\wc 0$ in $L^2(\Omega_m)$ according to~\eqref{eq:strain-wc} and~\eqref{eq:strain-limit-e3}, we have
	\begin{equation}
		\begin{aligned}
			\int_{U}\overline{\pr{x_3^2-1}\pr{e_3(\eps)+p_3}}c_1\sqrt{a}  \itd s=\frac{1}{2}\int_{\Omega} \pr{x_3^2-1}\pr{e_3(\eps)+p_3}c_1\sqrt{a}\itd x\to 0
		\end{aligned}
	\end{equation}
	which shows that
	\begin{equation}
		\int_{\toe}\pr{u^1-u^\eps}\Delta u^1\sim O(\eps^3),
	\end{equation}
	and thus~\eqref{eq:sec-der-eps3} is proved.
\end{proof}

\section{Upper bound of AAC for general shell lattice}
\label{sec:kavg-upp-general-shell}
In this section, we show that the upper bound $\frac{2}{3}\kappa$ for AAC in the main text holds for general shell lattice.
Specifically, the middle surface $\tw$ may not be closed, and the thickness of the shell may vary across $\tw$.
In this case, we assume the thickness distribution is given by two functions on $\tw$, denoted as $h^+_\eps$ and $h^-_\eps$, which represent the offset distances on both sides.
Note that we still require that $\tw$ is compact and the thickness of the shell tends to 0 as $\eps\to0$, i.e., the offset distances must satisfy $h^+_\eps,h^-_\eps\to 0$ as $\eps\to0$.

As $\tw$ remains compact, there still exist finite charts $(S_i,\boldsymbol{\psi}_i)_{i\in I}$ and partition $\br{\tw_m}$. 
Then we redefine 
\begin{equation}
	\label{eq:toe-vari-thick}
	\tO^{\eps}:=\left\{{\bx}\in\T^3: {\bx}={\by}+x_3\n({\by}), x_3\in\left(-h^-_\eps(\by),h^+_\eps(\by)\right),{\by}\in\tw\right\},
\end{equation}
which is still a 3-manifold when $\eps$ is sufficiently small.
The charts $(B_i^\eps,\P^\eps_i)_{i\in I}$ are redefined as
\begin{equation}
	\begin{aligned}
		B_i^{\eps}=\left\{{\bx}\in\T^3: {\bx}={\by}+x_3\n({\by}), x_3\in(-h^-_\eps(\by),h^+_\eps(\by)),{\by}\in S_i\right\}\\
		\label{eq:var-thick-shell-chart-map}
		\P_i^{\eps}: B_i^{\eps}\to \mathcal{V}_i^{\eps},\, {\by}+x_3\n({\by})\mapsto\left(\boldsymbol\psi_i(\by),x_3\right)
	\end{aligned}
\end{equation}
with $\mathcal{V}_i^{\eps}:=\br{(y,x_3):y\in \U_i, -h^-_\eps(y)<x_3<h^+_\eps(y)}$. The partition $\{\toe_m\}$ and the corresponding parameter domain $\{\Omega_m^\eps\}$ are defined accordingly, similar to the case of uniform thickness.

In this case, we denote the effective directional conductivity as $\ka(\toe;\p)$ and the corresponding volume fraction as $\rho(\toe):=|\toe|/|Y|$.
The asymptotic directional conductivity is still defined in a similar manner as in the case of uniform thickness:
\begin{equation}
	\ka_A(\p):=\lim_{\eps\to 0}\frac{\ka(\toe;\p)}{\rho(\toe)}
\end{equation}
where $\rho(\toe)$ denotes the volume fraction.
Let
\begin{equation}
	J_{\toe}(u,\p) := \frac{1}{|\toe|}\int_{\toe}\k\pr{\nb u+\p}\cdot\pr{\nb u+\p}
\end{equation}
and 
\begin{equation}
	\label{eq:j-Omega-eps-star}
	J_{\toe}^*(\p):=\min_{u\in V_\#(\toe)} J_{\toe}(u,\p).
\end{equation}
Then we define
\begin{equation}
	\ka_A(\p) = \lim_{\eps\to0}J_{\toe}^*(\p).
\end{equation} 

To show that the upper bound of $\ka_A^{Avg}$ also applies to $\ka_A(\p)$,
we construct a field $u^0$ whose coordinate representation is 
\begin{equation}
	\label{eq:u0-definition}
	u^0= -x_3 p_3\quad \text{ on each }\Omega_m^\eps.
\end{equation}
According to~\eqref{eq:j-Omega-eps-star}, we have the inequality
\begin{equation}
	J_{\toe}^*(\p)\le J_{\toe}(u^0,\p),
\end{equation}
which follows that
\begin{equation}
	\begin{aligned}
		J^*_{\toe}(\p)&\le \frac{1}{|\toe|}\int_{\toe}\k\pr{\nabla u^0+\p}\cdot\pr{\nabla u^0+\p}\\
		&=\frac{1}{|\toe|}\sum_m\int_{\Omega^\eps_m}\ka g^{ij,\eps}\pr{\partial_i u^0+p_i^\eps}\pr{\partial_j u^0+p_j^\eps}\sqrt{g^\eps}\\
		&=\frac{1}{|\toe|}\sum_m\int_{\Omega_m^\eps}\ka g^{\a\b,\eps}p_\a^\eps p_\b^\eps\sqrt{g^\eps},\\
	\end{aligned}
\end{equation}
where the last equality is from the definition of $u^0$.
Considering orthonormal basis $\{\p^i\}_{i=1,2,3}$, we have
\begin{equation}
	\begin{aligned}
		\frac{1}{3}\sum_i J_{\toe}^*(\p^i)&\le\frac{1}{|\toe|}\sum_m\int_{\Omega_m^\eps}\ka \frac{1}{3}\sum_ig^{\a\b,\eps} p_\a^{i,\eps} p_\b^{i,\eps}\sqrt{g^\eps}\\
		&= \frac{1}{|\toe|}\sum_m\int_{\Omega_m^\eps}\frac{2}{3}\ka \sqrt{g^\eps}= \frac{2}{3}\ka,
	\end{aligned}
\end{equation}
where the second line utilizes the following two identities
\begin{equation}
	\begin{aligned}
		g^{\a\b,\eps}p_\a^{i,\eps}p_\b^{i,\eps} &=\|\p^i \|^2 - \pr{\p^i\cdot\n}^2=1-\pr{\p^i\cdot\n}^2,\\
		\sum_i \pr{\p^i\cdot\n}^2&=\|\n\|^2=1.
	\end{aligned}
\end{equation}
Hence, we obtain
\begin{equation}
	\ka_A^{Avg}=\frac{1}{3}\sum_i\ka_A(\p^i)=\lim_{\eps\to0}\frac{1}{3}\sum_i J^*_{\toe}(\p^i)\le \frac{2}{3}\ka
\end{equation}
which is exactly the same bound as in main text.
Note that the above derivation does not require $\tw$ to be closed or the thickness to be uniform.

\section{Explicit form of ADC for revolution surfaces}
\label{app:exp-form-adc-revo}
First, we derive an explicit expression for ADC along the axial direction of  revolution surfaces:
\begin{theorem}
	\label{thm:revo-adc-ana}
	Let $\tw_R$ be a periodic surface of revolution about the $x$-axis, with radius profile given by $R(x)$. Then we have:
	\begin{equation}
		\label{eq:adc-ana-form}
		\ka_A(\tw_R;\e^1)=4 \left(\int_{-1}^1\frac{\sqrt{1+R'(x)^2}}{R(x)}\int_{-1}^1R(x)\sqrt{1+R'(x)^2}\right)^{-1}.
	\end{equation}
\end{theorem}
\begin{proof}
	Consider the following parameterization of $\tw_R$:
	\begin{equation}
		\label{eq:adc-rev-param}
		\bs{r}(x,\theta)=\pr{x,R(x)\cos\theta,R(x)\sin\theta},\quad x\in [-1,1],\theta\in[-\pi,\pi].
	\end{equation}
	Define
	\begin{equation}
		\begin{aligned}
			\bs a _1 &=\partial_x \bs r=(1,R'(x)\cos\theta,R'(x)\sin\theta),\\
			\bs a _2 &=\partial_\theta \bs r=(0,-R(x)\sin\theta,R(x)\cos\theta)\\
		\end{aligned}
	\end{equation}
	as the covariant basis vectors of the tangent space. In this basis, the covariant components of the tangential part of $\e^1$ are
	\[
	\e^1_{\tw}=(1,0).
	\]
	
	Due to the rotational symmetry of $\tw_R$, we can assume that the solution to \eqref{eq:avg-u-poisson} is of the form
	$s(x)$, i.e., independent of $\theta$.
	
	By the definition of divergence, we have:
	\begin{equation}
		\div \e^1_{\tw}=\frac{R'(x) \left(1-R(x) R''(x)+R'(x)^2\right)}{R(x) \left(1+R'(x)^2\right)^2}.
	\end{equation}
	
	Using the formula for the Laplace operator on a manifold, we obtain
	\begin{equation}
		\small
		\Delta \bar u = \frac{R(x) \pr{1+R'(x)^2} s''(x)+\left(R'(x)^3+R'(x) \left(1-R(x) R''(x)\right)\right)s'(x) }{R(x) \left(1+R'(x)^2\right)^2}.
	\end{equation}
	
	Then, the  equation \eqref{eq:avg-u-poisson} can be rewritten as the following second-order ordinary differential equation for $s(x)$:
	\begin{equation}
		s''(x)+w(x)s'(x)+w(x)=0,
	\end{equation}
	where
	\begin{equation}
		w(x)=\frac{R'(x) \left(1-R(x) R''(x)+R'(x)^2\right)}{R(x) \pr{1+R'(x)^2}}.
	\end{equation}
	
	The general solution is given by
	\begin{equation}
		\label{eq:gen-sol}
		\begin{aligned}
			s(x)=\int _1^x\frac{
				\sqrt{R'(\zeta )^2+1}}{R(\zeta )}\Bigg(c_1+
			\int _1^{\zeta }-\frac{R'(\xi ) \left(R'(\xi
				)^2-R(\xi ) R''(\xi )+1\right)}{\left(R'(\xi )^2+1\right)^{3/2}}d\xi \Bigg)d\zeta +c_2.
		\end{aligned}
	\end{equation}
	where $c_1,c_2$ are undetermined constants.
	Note that the integrand of the inner integral is precisely the derivative of $-\frac{R}{\sqrt{1+(R')^2}}$. Hence, the general solution~\eqref{eq:gen-sol} can be written as
	\begin{equation}
		\label{eq:rev-sol-c1}
		s(x)=\int _1^x \frac{
			\sqrt{R'(\zeta )^2+1}}{R(\zeta )}c_3-1\,d\zeta +c_2.
	\end{equation}
	
	According to  \eqref{eq:j-eps-limit-expr}, we only need the derivative of $s(x)$, thus it suffices to determine the constant $c_3$.
	
	Since $s(x)$ is periodic, we have
	\begin{equation}
		0=\int_{-1}^1 s'(x) dx =\int_{-1}^1 \frac{
			\sqrt{R'(x )^2+1}}{R(x )}c_1 -1\,dx,
	\end{equation}
	which gives
	\begin{equation}
		c_3=2\left(\int_{-1}^1\frac{
			\sqrt{R'(x )^2+1}}{R(x )} dx\right)^{-1}.
	\end{equation}
	
	Finally, by substituting $\bar u=s(x)$ and the area formula
	\begin{equation}
		|\tw_R|=2\pi\int_{-1}^1 R(x)\sqrt{1+R'(x)^2} dx
	\end{equation}
	into \eqref{eq:j-eps-limit-expr} and simplifying the result, we obtain \eqref{eq:adc-ana-form}.
\end{proof}

\section{Semi-analytical evaluation of effective conductivity}
\label{app:semi-any-eval}
To verify the third-order accuracy of ADC in approximating the effective thermal conductivity (Theorem~\ref{thm:second-order}), we thicken the revolution surface to form a shell lattice $\toe_R$ and consider the following parameterization:
\begin{equation}
	\begin{aligned}
		\bs{r^\eps}(x,\theta,z)&=\bs{r}(x,\theta) + z\n(x,\theta),\quad x\in [-1,1],\theta\in[-\pi,\pi], z\in [-\eps,\eps].
	\end{aligned}
\end{equation}
Here, $\bs{r}(x,\theta)$ is given by \eqref{eq:adc-rev-param}, and $\n(x,\theta)$ denotes the unit normal vector along $\tw_R$.

Define the following vectors:
\begin{equation}
	\g_1^\eps = \partial_x\bs r^\eps,\quad \g_2^\eps =\partial_\theta \bs r^\eps,\quad \g_3^\eps = \partial_z \bs r^\eps,
\end{equation}
as the covariant basis of the tangent space on the shell.

By symmetry, we can assume that the solution to the cell problem under $\p=\e^1$ takes the form
\begin{equation}
	u^\eps=s(x,z),
\end{equation}
i.e., $u^\eps$ is independent of $\theta$. Substituting it into \eqref{eq:min-energy} and simplifying yields:
\begin{equation}
	\label{eq:min-energy-discr}
	\kappa_\eps(\tw_R;\p)=\min_{s}\frac{1}{|Y|}\int_{-1}^12\pi\int_{-\eps}^\eps c_0 + \bs c_1^\top \nabla s + \nabla s^\top\bs C_2 \nabla s\, dz dx,
\end{equation}
where $c_0$, $\bs c_1$, and $\bs C_2$ are scalar, vector, and matrix functions, respectively,  of $x$ and $z$.

We uniformly discretize the domain $[-1,1]\times[-\eps,\eps]$ into $N\times M$ elements and denote the grid points by $(x_i,z_j)$. The value of $s$ at the grid point is denoted by $s_{i,j}$.

The gradient $\nabla s$ is approximated using finite difference scheme:
\begin{equation}
	\nabla s\approx \bs\zeta_{i,j}:=
	\left[
	\begin{aligned}
		\frac{s_{i+1,j}-s_{i-1,j}}{2\Delta x}\\
		\frac{s_{i,j+1}-s_{i,j-1}}{2\Delta z}
	\end{aligned}
	\right].
\end{equation}
Note that when $z_j$ lies on the boundary, the approximation of $\nabla_z s$ must be replaced with a one-sided finite difference, such as $\frac{s_{i,j}-s_{i,j-1}}{\Delta z}$ or $\frac{s_{i,j+1}-s_{i,j}}{\Delta z}$.
Since $s$ is periodic in $x$, we need to identify the two endpoints $x_0$ and $x_N$ as the same node.

After discretization, \eqref{eq:min-energy-discr} becomes:
\begin{equation}
	\label{eq:min-energy-discr-expand}
	\begin{aligned}
		\kappa_\eps(\tw_R;\p)=\min_{\bs s^h}\frac{1}{|Y|}\sum_{i=0}^{N-1}2\pi\Delta x\sum_{j=0}^{M}\Big(  c_0(x_i,z_j)
		+\bs c_1(x_i,z_j)^\top \bs \zeta_{i,j}+\bs\zeta_{i,j}^\top \bs C_2(x_i,z_j)\bs\zeta_{i,j}\Big)\Delta z_j.
	\end{aligned}
\end{equation}
Here, $\bs s^h$ denotes the vector of nodal values, and $\Delta z_j$ is computed as:
\begin{equation}
	\Delta z_j =\left\{
	\begin{aligned}
		&\Delta z/2, \quad& j=0, M\\
		&\Delta z, & \text{otherwise}
	\end{aligned}
	\right..
\end{equation}
Clearly, \eqref{eq:min-energy-discr-expand} is a quadratic minimization problem in terms of $\bs s^h$, which can be solved by a linear system.

To ensure sufficient accuracy, we set $N=100000$ and $M=20$.

\section{Time derivative of ADC in normal flow}
\label{app:dtime-flow}
Suppose the middle surface $\tw$ undergoes a normal velocity $v_n:\tw\to\mathbb R$, we can first verify the following time derivatives 
\begin{align}
	\label{eq:a-time-deri}
	\dt{a}^{\a\b}&=2v_n b^{\a\b}\\
	\dt{\sqrt{\det a}}&=-2v_n H\sqrt{\det a}\\
	\label{eq:pa-rate}
	\dt{p}_\a &= p_3\partial_\a v_n - v_n b_\a^\b p_\b
\end{align}

For a time dependent integrand defined on $\tw(t)$, the time derivative of its integral is
\begin{equation}
	\label{eq:mov-w-f-der}
	\frac{d}{dt}\int_{\tw(t)} f(t)=\int_{\tw(t)}\dt f(t)-\int_{\tw(t)}2v_n H f(t)
\end{equation}
For the surface area $A:=|\tw|$, the derivative is $\dt{A}=-\int_{\tw}2v_nH$.

Now, we give the proof of Proposition 1 in the main text:
\begin{proof}
	Applying above formulas, we deduce the shape derivative of ADC as
	\begin{equation}
		\begin{aligned}
			\dt{I}_A^c&=\frac{d}{dt}\int_{U}\ka a^{\a\b}(\partial_\a\bar u+p_\a)(\partial_\b\bar u+p_\b)\sqrt{\det a} \itd s\\
			&=\int_{U}\ka 2v_n b^{\a\b}(\partial_\a\bar u+p_\a)(\partial_\b\bar u+p_\b)\sqrt{\det a} \itd s \\
			&+\int_{U}\ka a^{\a\b}\pr{\partial_\a\bar u+p_\a}\pr{\partial_\b\dt{\bar u}+\dt{p}_\b}\sqrt{\det a} \itd s\\
			&+\int_U \ka a^{\a\b}\pr{\partial_\a\bar u+p_\a}\pr{\partial_\b\bar u+p_\b}(-2v_n H\sqrt{\det a}) \itd s\\
			&=\int_{\tw} 2v_n\ka\pr{\bs b- H\mathbf I}\pr{\nabla\bar u+\p_{\tw}}\cdot\pr{\nabla\bar u+\p_{\tw}} \itd \tw \\
			&+\int_{\tw}\ka\pr{\nabla\bar u+\p_{\tw}}\cdot\pr{\nabla\dt{\bar u}+\dt\p_{\tw}}\itd\tw
		\end{aligned}
	\end{equation}
	Substituting $v=\dt{\bar u}$  into the variational equation~\eqref{eq:u-avr-vari-form}, we have 
	\begin{equation}
		\int_{\tw}\ka\pr{\nabla\bar u+\p_{\tw}}\cdot\nabla\dt{\bar u} =0
	\end{equation}
	Thus it suffices to show that
	\begin{equation}
		\int_{\tw}\ka\pr{\nabla\bar u+\p_{\tw}}\cdot\dt{\p}_{\tw} = 0.
	\end{equation}
	From~\eqref{eq:pa-rate}, we get
	\begin{equation}
		\int_{\tw}\ka\pr{\nabla\bar u +\p_{\tw}}\cdot\dt{\p}_{\tw}=\int_{\tw}\ka\pr{\nabla\bar u+\p_{\tw}}\cdot\pr{p_3\nabla v_n-v_n\bs b\cdot\p_{\tw}}
	\end{equation}
	Note that
	\begin{equation}
		\label{eq:dot-p-eq-zero}
		\begin{aligned}
			p_3\nabla v_n &= \nabla(p_3 v_n) - v_n\nabla p_3\\
			\nabla p_3 &= -\bs b\cdot\p_{\tw}.
		\end{aligned}
	\end{equation}
	Hence, by using~\eqref{eq:u-avr-vari-form} again, the term involving $\nabla(p_3 v_n)$ is eliminated, and the remaining terms cancel out.
\end{proof}

\end{document}